\documentclass[10pt]{amsart}
\usepackage{}
\usepackage{amsfonts}
\usepackage{amsfonts,latexsym,rawfonts,amsmath,amssymb,amsthm}
\usepackage[plainpages=false]{hyperref}
\usepackage{mathrsfs}
\RequirePackage{amsmath} \RequirePackage{amssymb}
\usepackage{enumitem}
\usepackage{fancyhdr}
\usepackage{amscd,latexsym,amsthm,amsfonts,amssymb,amsmath,amsxtra}
\usepackage{amscd,mathrsfs,latexsym,amsmath,amsthm,amssymb,amsxtra,}
\usepackage{latexsym,amsmath,amsthm,amssymb,amsxtra,mathrsfs}
\usepackage[all]{xy}
\usepackage{longtable}
\usepackage{verbatim}
\usepackage{bm}
\usepackage{mathtools}
\usepackage{color}



\numberwithin{equation}{section}

\renewcommand{\theequation}{\arabic{section}.\arabic{equation}}

\newcommand{\al}{\alpha}

\newcommand{\la}{\lambda}

\def\Ga{\Gamma}
\def\ga{\gamma}

\newcommand{\si}{\sigma}

\newcommand{\beq}{\begin{equation}}
\newcommand{\eeq}{\end{equation}}
\newcommand{\beqs}{\begin{eqnarray*}}
\newcommand{\eeqs}{\end{eqnarray*}}
\newcommand{\beqn}{\begin{eqnarray}}
\newcommand{\eeqn}{\end{eqnarray}}
\newcommand{\beqa}{\begin{array}}
\newcommand{\eeqa}{\end{array}}

 \DeclareMathOperator{\sym}{Sym}

\def\var{\varepsilon}
\def\lra{\longrightarrow}

\def\la{\lambda}

\def\na{\nabla}

\def\lan{\langle}
\def\ran{\rangle}

\def\bc{\begin{center}}
\def\ec{\end{center}}

\def\o{\omega}

\def\S{\Sigma}

\def\pr{\prime}

\def\Z{{\mathbb Z}}

\def\tr{\rm tr}

\def\Id{{\rm Id}}

\def\begeq{\begin{equation}}
\def\endeq{\end{equation}}
\def\and{\quad{\rm and}\quad}

\def\Del{\Delta}

\let\lra=\longrightarrow

\def\mapright\#1{\,\smash{\mathop{\lra}\limits^{\#1}}\,}

\def\si{\sigma}

\def\Sph{\mathbb{S}}
\def\IDW{\mathcal{IDW}}
\def\md{\mathrm{d}}
\def\Rm{\mathrm{Rm}}
\def\Ric{\mathrm{Ric}}
\def\Vol{\mathrm{Vol}}
\def\ad{\mathrm{ad}}
\def\tr{\mathrm{tr}}
\def\Scal{\mathrm{Scal}}
\newtheorem{prop}{Proposition}[section]
\newtheorem{theo}[prop]{Theorem}
\newtheorem{lem}[prop]{Lemma}
\newtheorem{claim}[prop]{Claim}
\newtheorem{cor}[prop]{Corollary}
\newtheorem{rem}[prop]{Remark}

\newtheorem{conj}[prop]{Conjecture}

\begin{document}
\bibliographystyle{plain}

\title  {Classification of cohomogeneity one steady gradient Ricci solitons in dimensions four and five}

\date{}

\author {Yuxing $\text{DENG}^*$ }
\author { Guohao $\text{WANG}$ }

    \thanks {* Supported by National Key R$\&$D Program of China 2022YFA1007600 and NSFC No. 12571059.}

\address{ Yuxing Deng\\School of Mathematics and Statistics, Beijing Institute of Technology,
Beijing, 100081, China\\
6120180026@bit.edu.cn}

\address{ Guohao Wang\\School of Mathematics and Statistics, Beijing Institute of Technology,
Beijing, 100081, China\\
3120245981@bit.edu.cn}

\begin{abstract}
   In this paper, we first prove Appleton's conjecture on the uniqueness of noncollapsed steady solitons constructed in \cite{Appleton}. We also classify non-Ricci-flat $\kappa$-noncollapsed cohomogeneity one steady gradient Ricci solitons in dimensions four and five. After passing to a finite cover, in dimension four, they must be the Bryant soliton or Appleton's noncollapsed 4d steady soliton; in dimension five, they must be the Bryant soliton or solitons constructed by Ivey--Dancer--Wang.
\end{abstract}

\maketitle


\section{introduction}

A complete Riemannian manifold $(M^n,g)$ is called a steady gradient Ricci soliton if there is a smooth function $f$ such that
\begin{align}\label{eq:steady}
    \Ric+\na^2 f=0.
\end{align}
A steady gradient Ricci soliton generates an eternal Ricci flow by pullback under the one-parameter family of diffeomorphisms generated by $\na f$. 


Steady gradient Ricci solitons are important since they may arise as the type II singularity models of the Ricci flow.  In dimension three, Brendle proved that every complete nonflat $\kappa$-noncollapsed steady gradient Ricci soliton is isometric to the Bryant soliton, which solving a conjecture of Perelman\cite{Per02, Bre13}.  Brendle proved rotational symmetry for steady solitons with positive sectional curvature and asymptotically cylindrical condition\cite{Bre14}. Later on, Zhu and the first named author generalized Brendle's work and classified steady solitons with linear scalar curvature decay\cite{DZ19, DZ2020, DZ20b}. For further generalizations, see \cite{ZZ2022, Law}. There are also many works on the rigidity of steady gradient Ricci solitons, see \cite{CC12,CX,Lai25,MMS,Fang} and references therein. 

In higher dimensions, the classification of $\kappa$-noncollapsed steady gradient Ricci solitons remains widely open. There are quite a lot of $\kappa$-noncollapsed steady gradient Ricci solitons for $n\ge 4$. Bryant constructed steady gradient Ricci solitons with positive curvature operator and $O(n)$-symmetry \cite{Bryant}. In \cite{Appleton}, Appleton constructed  steady gradient Ricci solitons with $U(n)$-symmetry which is non-K\"ahler. Ivey also constructed complete non-K\"ahler steady solitons using doubly warped products\cite{Ivey}. Dancer-Wang extended this construction to multiply warped products with positive Einstein factors\cite{DancerWangExamples}. Buzano-Dancer-Wang constructed further families of steady solitons and Ricci-flat metrics with different asymptotic behavior\cite{BDW2015}. These examples are all cohomogeneity one steady Ricci solitons. Other cohomogeneity one steady Ricci soliton examples include \cite{Cao1996, DancerWang, Yang2012, BDW2016, Stolarski2015, Wink2021, Wink2023, Wang2024, Chi2024, Chi2026}.
 Moreover, there are also some fly wing type $\kappa$-noncollapsed steady Ricci solitons. Lai constructed  steady Ricci solitons with positive curvature operator and $O(n-1)\times \mathbb{Z}_2$-symmetry. These solitons are noncollapsed for $n\ge 4$ \cite{Lai24}. See \cite{CLL} for recent progress on constructing fly wing type steady Ricci solitons. We also remark that for examples of steady K\"ahler Ricci solitons, one may see \cite{Cao1996,CD,AC,CCL} and references therein.

In dimension four, Bamler-Chow-Ma-Zhang and the first named author classified the tangent flows of steady gradient Ricci solitons which are singularity models \cite{BCDMZ2022}.  An important type of these solitons have tangent flow $\mathbb{S}^3/\Gamma\times \mathbb{R}$. As a special case of the main theorem in \cite{Bre14},  steady gradient Ricci solitons with positive sectional curvature  and tangent flows $\mathbb{S}^3/\Gamma\times \mathbb{R}$ must be isometric to the Bryant soliton. The first named author generalized Brendle's results to steady gradient Ricci solitons on orbifolds \cite{Deng2025}, i.e.  steady gradient Ricci solitons with positive sectional curvature  and tangent flows $\mathbb{S}^3/\Gamma\times \mathbb{R}$ must be isometric to ${\rm Bry}/\Gamma$. In dimension $4$, Appleton also has constructed steady gradient Ricci solitons with tangent flow $\mathbb{S}^3/\mathbb{Z}_k\times \mathbb{R}$ for $k>2$. There is a folklore conjecture that \textit{any four-dimensional steady gradient Ricci soliton with tangent flow $\mathbb{S}^3/\Gamma\times \mathbb R$ should be a cohomogeneity one steady gradient Ricci soliton if $\mathbb{S}^3/\Gamma$ is homogeneous.} It remains interesting to classify $\kappa$-noncollapsed cohomogeneity one steady gradient Ricci solitons in dimension at least four. In  \cite{RT}, 
Ribeiro and Tran classified complete four-dimensional steady gradient Ricci solitons that are conformally Kähler with isoparametric potential. They showed that these hypotheses force cohomogeneity one $U(2)$-symmetry, and their classification contains both the collapsed and the
noncollapsed solitons. In particular, their results imply a classification
of the cohomogeneity one steady gradient Ricci solitons with $U(2)$-symmetry. However, Buttsworth \cite{B2025} has constructed  $SU(2)$-symmetric cohomogeneity one steady gradient Ricci solitons which are not $U(2)$-symmetric. So the classification of four-dimensional $\kappa$-noncollapsed cohomogeneity one steady gradient Ricci solitons is still open.

In this paper, we first show that the universal cover of a complete, non-Ricci-flat and $\kappa$-noncollapsed steady gradient Ricci soliton with cohomogeneity one symmetry is either the Bryant soliton or the Appleton's soliton.

\begin{theo}\label{thm: four}
Let $(M^4,g,f)$ be a complete, non-Ricci-flat and $\kappa$-noncollapsed steady gradient Ricci soliton.  Suppose that a compact connected Lie group $G$ acts almost effectively, isometrically, and with cohomogeneity one.  Then there is a finite Riemannian covering
\[
  (\widehat M,\widehat g,\widehat f)\longrightarrow(M,g,f)
\]
such that $\mathrm{Iso}(\widehat M,\widehat g)$ contains a connected subgroup locally isomorphic to $U(2)$.

More precisely, there are only two possibilities:
\begin{enumerate}[label=\textup{(\roman*)}]
\item If the unique singular orbit is a point, then $g$ is rotationally symmetric and, after scaling, $(M,g,f)$ is the four-dimensional Bryant soliton.
\item If the unique singular orbit has codimension two, then for some integer $k\ge 3$, $\widehat{M}$ is diffeomorphic to $\mathcal O(-k)$, and after scaling, the lifted soliton is the noncollapsed 4d Appleton soliton.  Moreover, if $k$ is odd, then $(M,g,f)$ itself~(no covering) is the noncollapsed 4d Appleton soliton on $\mathcal O(-k)$. If $k$ is even, then $M$ is the noncollapsed 4d Appleton soliton on $\mathcal O(-k)$ or the free $\mathbb{Z}_2$-quotient of it.
\end{enumerate}
\end{theo}

As a corollary, we can see that the non-$U(2)$-invariant steady solitons constructed by Buttsworth in \cite{B2025} are collapsed.


 Of particular relevance in dimension four are Appleton's noncollapsed, non-K\"ahler steady solitons on the total space of complex line bundles $\mathcal{O}(-k)\to \mathbb{CP}^1$ for $k\geq 3$~\cite{Appleton}. In fact, Appleton's construction generalizes to the existence of noncollapsed steady soliton on the total space of complex line bundle over any compact positive K\"ahler-Einstein manifold.  
 
 Now, we recall Appleton's construction. Let $(B^{2n},J,\widehat g)$, $n\geq1$, be a compact connected
K\"ahler--Einstein manifold with
$\operatorname{Ric}_{\widehat g}=2(n+1)\widehat g$.
Let $p=p(B)$ be the positive divisibility of the first Chern class:
\begin{align}\label{eq:def-q-1}
    p:=\max\bigl\{q\in\mathbb N:
                 c_1(B)\in qH^2(B;\mathbb Z)\bigr\}.
\end{align} 
Then $c_1(B)=p\eta$ for some primitive $\eta\in H^2(B;\mathbb Z)$. For an integer $k>p$, let $\mathcal{L}_k$ be the total space of the complex line
bundle over $B$ with first Chern class $-k\eta$. On the complement of the zero section, metric $g$ is of the form
\begin{equation}\label{eq:ansatz-1}
 g=dt^2+a(t)^2\theta\otimes\theta
       +b(t)^2\pi^*\widehat g,
\end{equation}
where $\pi$ is the bundle projection and $t$ is distance from the zero
section. Appleton proved that there exist a smooth function $f=f(t)$ and a metric $g$ of the form (\ref{eq:ansatz-1}) such that $(\mathcal{L}_k,g,f)$ is a complete $\kappa$-noncollapsed steady gradient Ricci soliton. See Section~2 for the detailed description of Appleton's construction. However, it is unknown whether the $\kappa$-noncollapsed soliton on $\mathcal{L}_k$ constructed by Appleton is unique. Appleton made the following conjecture.

\begin{conj}[\cite{Appleton}]\label{conj:Appleton}
Let $\mathcal{L}_k$ be the total space of the complex line
bundle over $B$ as defined above. $p$ is defined by (\ref{eq:def-q-1}). For any fixed integer $k>p$, up to scaling and isometry, there is exactly one complete noncollapsed steady gradient Ricci soliton on $\mathcal{L}_k$ of the form
\eqref{eq:ansatz-1} that extends smoothly across the zero section.
It is the soliton constructed in \cite[Theorem~1.5]{Appleton}.
\end{conj}

 In Appleton's construction,  after scaling the metric and adding a constant to the potential, we impose
\begin{equation}\label{eq:initial-data-1}
\begin{gathered}
 a(0)=0,\qquad a'(0)=\frac{(n+1)k}{p},\qquad
 b(0)=1,\\ b'(0)=0, \qquad
 f(0)=0,\qquad f'(0)=0.
\end{gathered}
\end{equation}
Then the reduced soliton equations determine a unique smooth
local solution, which is analytic, once
\[
 \mu=f''(0)
\]
is prescribed~(cf.~Proposition~\ref{prop:prelim-local}), where $f$ is the soliton potential. The remaining problem is to determine which values of $\mu$ yield complete
metrics and which yield noncollapsed metrics.
Our next result solves the conjecture of Appleton.

\begin{theo}\label{thm:conjecture}
Let $\mathcal{L}_k$ be the total space of the complex line
bundle over $B$ in Conjecture \ref{conj:Appleton}. $p$ is defined by (\ref{eq:def-q-1}). For any fixed integer $k>p$,  up to scaling and isometry, there is exactly one complete
noncollapsed steady gradient Ricci soliton of the form
\eqref{eq:ansatz-1} that extends smoothly across the zero section.

More precisely, under the normalization
\eqref{eq:initial-data-1}, there is a unique $\mu_*<0$ such that the
maximal solution corresponding to $\mu=f''(0)$ is
\begin{equation}\label{miu}
 \begin{cases}
  \text{incomplete},&\mu>\mu_*,\\
  \text{complete and noncollapsed},&\mu=\mu_*,\\
  \text{complete and collapsed},&\mu<\mu_*.
 \end{cases}
\end{equation}
In the noncollapsed case, $a(t)/b(t)\to1$ as $t\to\infty$; in the
collapsed case, $a(t)/b(t)\to0$ as $t\to\infty$.
\end{theo}

In particular, in dimension 4, we get the uniqueness of Appleton's noncollapsed non-K\"ahler steady solitons on the total space of the complex line bundles  $\mathcal{O}(-k),k\geq 3$ over $\mathbb{CP}^1$.

\begin{cor}\label{cor:unique four noncollapsed soliton}
    For every integer $k\geq 3$, the total space of
$\mathcal O(-k)\to\mathbb{CP}^1$ carries, up to scaling and isometry,
exactly one complete noncollapsed steady gradient Ricci soliton whose
metric is invariant under the standard $U(2)$-action. This unique noncollapsed steady soliton is exactly the one constructed in~\cite[Theorem~1.2]{Appleton}, which we call the noncollapsed 4d Appleton soliton on $\mathcal{O}(-k)$.
\end{cor}

\begin{rem}
  In dimension four, Corollary 1.4 and the trichotomy (\ref{miu}) were also
obtained by Ribeiro and Tran \cite[Theorem 6.1]{RT}, with a different proof.
Theorem \ref{thm:conjecture} holds for Appleton's solitons in all even dimensions. 
\end{rem}

In dimension five, if we take the Einstein manifold as the round $\Sph^2$ in Ivey's construction\cite{Ivey}, then we get a one-parameter family of five dimensional cohomogeneity one steady gradient Ricci solitons on $\mathbb{R}^3\times \Sph^2$. Note that this family of solitons coincide with the examples constructed by Dancer--Wang\cite{DancerWangExamples} when there are two identical round $\Sph^2$ Einstein factors. We use $\IDW$ to denote such one-parameter family steady solitons. Under the same conditions as in the four dimensional setting, cohomogeneity one steady solitons in dimension five must be the Bryant soliton or belong to the family $\IDW$.

\begin{theo}\label{thm:five}
Let $(M^5,g,f)$ be complete,  non-Ricci-flat and
$\kappa$-noncollapsed steady gradient Ricci soliton.  Suppose that a compact connected Lie group acts
almost effectively, isometrically, and with cohomogeneity one.  Then, after a
finite equivariant Riemannian covering $(\widehat M, \widehat g)$ and a homothety, exactly one of the
following occurs:
\begin{enumerate}[label=\textup{(\roman*)}]
\item $(\widehat M,\widehat g)$ is the five-dimensional Bryant soliton on
      $\mathbb{R}^5$;
\item $\widehat M\cong\mathbb{R}^3\times\Sph^2$, and $\widehat g$ belongs to the
      one-parameter Ivey--Dancer--Wang family $\IDW$.
\end{enumerate}
\end{theo}

The paper is organized as follows. In Section~2, we recall some preliminaries of cohomogeneity one steady gradient Ricci solitons and some results proved by Appleton. In Section~3, we construct a quantity using the warping functions to establish the uniqueness of Appleton's solitons, hence proving Theorem \ref{thm:conjecture}. In Section~4, we establishes the linear decay of the curvature operator for cohomogeneity one steady solitons with only one singular orbit. Consequently, the entropy estimates and finiteness of the fundamental group follow directly. In Section~5 and 6, we prove the four dimensional results, where the asymptotic property and a ratio maximum principle are essential. Section~7 and 8 treat the five dimensional classification, in which we first establish a topology trichotomy for five dimensional cohomogeneity one manifolds with only one singular orbit, then use the noncollapsing condition to exclude the case with $\Sph^1\times \Sph^3$ as principal orbit. Appendix~A gives detailed computations used in Section~3. Appendix~B records the homogeneous space classification used in the argument.

\medskip
\noindent\textbf{Acknowledgement.} The authors would like to thank Professor Ribeiro and Professor Tran for
pointing out their work \cite{RT}, which proves Appleton's uniqueness conjecture
in dimension four and coincides with the four-dimensional case of Theorem \ref{thm:conjecture}. Our proof is different and holds in all even dimensions. The two works were
carried out independently.

\medskip
\noindent\textbf{Declaration of AI use.} The main ideas underlying the classification of cohomogeneity one steady gradient Ricci solitons and the proof of Appleton’s uniqueness conjecture were developed by the authors. ChatGPT was used to do some calculations in the proof of Appleton’s uniqueness conjecture and to improve the language of the paper. The authors independently verified the AI-assisted material and take full responsibility for the mathematical arguments, results, and final content of this paper.

\section{Preliminaries of cohomogeneity one actions and Appleton's construction}

In this section, we first recall some basics about steady gradient Ricci solitons and cohomogeneity one actions. Then we use these to reduce our classification problem to cohomogeneity one actions with only one singular orbit of codimension at least 2. Finally, for the proof of Appleton's uniqueness conjecture, we record some results proved by Appleton in \cite{Appleton} which will be used in the next section.

For the steady gradient Ricci soliton, throughout the paper we use the sign convention \eqref{eq:steady}.  Tracing and differentiating the soliton equation \eqref{eq:steady}, we can get
\begin{equation}\label{eq:soliton-identities}
  R+\Delta f=0,
  \qquad
  \nabla R=2\mathrm{Ric}(\nabla f,\cdot),
  \qquad
  R+|\nabla f|^2=C_0,
\end{equation}
where $C_0$ is a constant. By the result of Chen~(cf.~\cite{Chen}), $C_0>0$ for a complete non-Ricci-flat steady gradient Ricci soliton. 

In the following, we always consider steady gradient Ricci solitons invariant under a cohomogeneity one action. It means that the soliton metric is invariant under a compact Lie group action and such action is cohomogenity one. If so, a natural question is whether the potential function $f$ is invariant under this action. Fortunately, the following lemma tells us that we can always assume the potential invariance.

\begin{lem}\label{lem:average}
Let $G$ be a compact group acting isometrically on a steady gradient Ricci soliton $(M,g,f)$.  Then there is a $G$-invariant function $\bar f$ satisfying
\[
  \mathrm{Ric}+\nabla^2\bar{f}=0.
\]
\end{lem}

\begin{proof}
Let $d\mu_G$ be the Haar measure on $G$ normalized such that
\begin{align*}
    \int_Gd\mu_G=1.
\end{align*}
Then we can define a function $\bar{f}$ on $M$ by
\[
  \bar f(x):=\int_G f(a\cdot x)\,d\mu_G(a),\qquad \forall\,x\in M.
\]
For every $a\in G$, the isometry $a$ preserves $\mathrm{Ric}_g$ and commutes with the Hessian. Hence,
\[
  \nabla^2(f\circ a)=a^*(\nabla^2 f)=-\mathrm{Ric}.
\]
Integrating this identity with respect to $d\mu_G$ gives us that $\nabla^2\bar f=-\mathrm{Ric}$. For the $G$-invariance of $\bar f$, it follows from the definition of $\bar{f}$ and the invariance of Haar measure.
\end{proof}

Therefore, in the following, we can always assume the potential function $f$ is also $G$-invariant for a cohomogeneity one steady gradient Ricci soliton with action group $G$.

Next, we recall some basics about cohomogeneity one actions. Let $G$ be a compact and connected Lie group. If $(M,g)$ admits an almost effective, isometric and cohomogeneity one $G$-action, then we say $(M,g)$ is a $G$-invariant cohomogeneity one Riemannian manifold. For cohomogeneity one action, it means that the orbit space $M/G$ is one dimensional. In this case, the orbit space $M/G$ is one of
\[
  \mathbb{S}^1,\quad [0,1],\quad \mathbb{R},\quad [0,\infty).
\]
The first two cases imply that $M$ is a compact manifold. For a compact steady gradient Ricci soliton, we know that it is Ricci-flat.

If $M/G=\mathbb{R}$, all orbits are principal and $M$ has two ends. Such steady solitons are also Ricci-flat, due to the following theorem of Munteanu--Wang~\cite[Theorem~4.2]{MW2022}.

\begin{theo}\label{thm:two-end}
A complete steady gradient Ricci soliton with at least two ends is isometric to
\[
  (\mathbb{R}\times N,d s^2+h),
\]
where $(N,h)$ is compact and Ricci-flat.  In particular, the soliton metric is Ricci-flat.
\end{theo}

It follows immediately that a non-Ricci-flat soliton in our setting has orbit space $[0,\infty)$. In this case, there is only one boundary orbit and the cohomogeneity one action is described by a group diagram
\[
  H\subset K\subset G,
\]
where $G/H$ is a principal orbit, $G/K$ is the boundary orbit and
\[
  K/H\cong \mathbb{S}^\ell
\]
for some $\ell\ge0$.  If $\ell=0$, then the boundary orbit has the same dimension as the principal orbit and is usually called exceptional rather than singular.  The following proposition exclude this case in our classification settings.

\begin{prop}\label{prop:double}
Let $G$ be a compact connected Lie group and $(M,g)$ be a complete $G$-invariant cohomogeneity one manifold with orbit space $[0,\infty)$. If this cohomogeneity one action has group diagram $H\subset K\subset G$ with $K/H\cong\mathbb{S}^0$. Then $(M,g)$ admits a two-fold Riemannian covering $\widehat M\to M$ such that $\widehat M/G\cong\mathbb{R}$.  If $M$ carries a $G$-invariant gradient steady soliton structure, the lifted structure on $\widehat M$ is complete and has two ends.
\end{prop}

\begin{proof}
Since $K/H\cong \mathbb S^0$, the slice representation of $K$ on the one-dimensional normal space has kernel $H$ and image $\{\pm1\}$.  Thus a neighborhood of the exceptional orbit is
\[
  G\times_K\mathbb{R},
\]
where the nontrivial element of $K/H$ acts on $\mathbb{R}$ by reflection.  The natural map
\[
  G\times_H\mathbb{R}=G/H\times\mathbb{R}\longrightarrow G\times_K\mathbb{R},
  \qquad [g,t]_H\longmapsto[g,t]_K,
\]
is a two-fold covering.  Since there is no second boundary orbit, this local normal cover extends along the principal part and gives a global two-fold cover $\widehat M\to M$.  Its orbit space is the double of $[0,\infty)$ across its boundary, hence is $\mathbb{R}$.

The metric and invariant potential lift smoothly and completeness is also preserved by a finite covering. Since $G$ is a compact Lie group, the orbits under $G$ action are all compact. Using this and the orbit space $\widehat M/G=\mathbb R$, we can see $\widehat M$ has two ends.
\end{proof}

\begin{cor}\label{cor:no-exceptional}
A complete non-Ricci-flat cohomogeneity one steady gradient Ricci soliton has exactly one singular orbit. In particular, this singular orbit has codimension at least two.
\end{cor}

\begin{proof}
By the discussion above, if the orbit space is a closed interval or a circle, then the manifold is compact. For a steady gradient Ricci soliton, it is Ricci-flat. If the orbit space is a line, then by the argument in the proof of Proposition \ref{prop:double} we can see that it has two ends, which is also Ricci-flat by Theorem \ref{thm:two-end}. Hence the orbit space must be $[0,+\infty)$, i.e., there is only one non-principal orbit. Proposition \ref{prop:double} tells us that if the non-principal orbit is not singular, then the universal cover of it has two ends and hence Ricci-flat.  In conclusion, such a steady soliton in Corollary \ref{cor:no-exceptional} must have exactly one singular orbit with codimension at least two.
\end{proof}

For cohomogeneity one non-Ricci-flat steady solitons, the soliton equation \eqref{eq:steady} can be reduced to a system of ODEs. Let $\S=G/K$ be the singular orbit and $t=d(\,\cdot\,,\S)$.  On the
principal part,
\begin{equation}\label{eq:general-cohom1}
  g=\md t^2+g_t,\qquad f=f(t).
\end{equation}
Let
\begin{align}\label{eq:def-LHxi}
    L=\frac12g_t^{-1}\dot g_t,\qquad
  H={\tr} L,\qquad
  \xi=H-\dot f,
\end{align}
where a dot denotes differentiation with respect to $t$. And let $r_t$ and $s_t$ denote the intrinsic Ricci endomorphism and scalar
curvature of $(G/H,g_t)$.  The radial, mixed, and tangential components of
the soliton equation are:(cf.~\cite[Section 2]{BDW2015})
\begin{align}
  -{\tr}\dot L-{\tr}(L^2)+\ddot f&=0,\label{eq:cohom-radial}\\
  \delta^{g_t}L&=0,\label{eq:cohom-mixed}\\
  \dot L&=r_t-\xi L.\label{eq:cohom-tangent}
\end{align}
In particular,
\begin{equation}\label{eq:xi-H}
  \dot\xi=-{\tr}(L^2),\qquad
  \dot H=s_t-\xi H.
\end{equation}
The volume $V(t)=\Vol(G/H,g_t)$ satisfies
\begin{equation}\label{eq:Vprime}
  \dot V=HV.
\end{equation}

The following theorem record some standard cohomogeneity one estimates.

\begin{theo}\label{thm:cohom-est}
Let $(M^{m+1},g,f)$ be a complete non-Ricci-flat cohomogeneity one steady
gradient Ricci soliton with one singular orbit and the normalization \eqref{eq:soliton-identities}. Then
\[
  \dot \xi(t)<0,\qquad \dot f(t)<0,\qquad 0<H(t)\le\frac{m}{t}, \qquad \forall\,t>0,
\]
and
\[
  \lim_{t\to+\infty}-\dot f(t)= \sqrt{C_0},\qquad
  \lim_{t\to+\infty}\xi(t)= \sqrt{C_0}.
\]
Moreover, $f$ is proper, $R(t)$ is strictly decreasing and there is a constant $C$ such that
\begin{equation}\label{eq:scalar-linear}
  0\le R(t)\le\frac{C}{1+t},
\end{equation}
 where $t(x)=d(x,\S)$ and  $\S=G/K$ denotes the singular orbit. 
\end{theo}

\begin{proof}
These are directly follow from Proposition 2.3, 2.4 and Corollary 2.6 in \cite{BDW2015}.
\end{proof}

For Appleton's construction, he considered a one variable line bundle ansatz on the total space of complex line bundles over positive compact K\"ahler-Einstein manifold. More precisely, let $(B^{2n},J,\widehat g)$, $n\geq1$, be a compact connected
K\"ahler--Einstein manifold with K\"ahler form $\widehat\omega$, normalized by
\[
 \operatorname{Ric}_{\widehat g}=2(n+1)\widehat g.
\]
Let $p=p(B)$ be the positive divisibility of the first Chern class defined by
\begin{align}\label{eq:def-q}
    p:=\max\bigl\{q\in\mathbb N:
                 c_1(B)\in qH^2(B;\mathbb Z)\bigr\}.
\end{align}
Then $c_1(B)=p\eta$ for some primitive $\eta\in H^2(B;\mathbb Z)$. For $k\in\mathbb N$, let $\mathcal L_k$ be the total space of
the complex line bundle over $B$ with first Chern class $-k\eta$.
In particular, if $B=\mathbb{CP}^n$, then $\eta$ is the hyperplane class, $p=n+1$ and $\mathcal L_k=\mathcal O(-k)$.

As in Appleton's construction, on the complement of the zero section, we consider $(g,f)$ of the form
\begin{equation}\label{eq:ansatz}
 g=dt^2+a(t)^2\theta\otimes\theta
       +b(t)^2\pi^*\widehat g,
 \qquad f=f(t),
\end{equation}
where $\pi$ is the bundle projection and $t$ is distance from the zero
section. With Appleton's normalization of the connection form,
locally $\theta=d\tau-2A$, where $dA=\widehat\omega$ and $\tau$ has period $2\pi p/((n+1)k)$. For the projective base, this is the standard $U(n+1)$-invariant ansatz.

After scaling the metric and adding a constant to the potential, we impose
\begin{equation}\label{eq:initial-data}
\begin{gathered}
 a(0)=0,\qquad a'(0)=\frac{(n+1)k}{p},\qquad
 b(0)=1,\qquad b'(0)=0,\\
 f(0)=0,\qquad f'(0)=0.
\end{gathered}
\end{equation}
We also require that $a$ extend smoothly and oddly across $t=0$, and that
$b$ and $f$ extend smoothly and evenly.

 In this case, the cohomogeneity one soliton system \eqref{eq:cohom-radial}--\eqref{eq:cohom-tangent} reduces to the following system of $a,b,f$~(cf.~\cite[(2.1)--(2.3)]{Appleton}).

\begin{subequations}\label{eq:prelim-ode}
\begin{align}
 f''&=\frac{a''}{a}+2n\frac{b''}{b},\label{eq:prelim-fode}\\
 a''&=2n\left(\frac{a^3}{b^4}-\frac{a'b'}{b}\right)+a'f',
 \label{eq:prelim-aode}\\
 b''&=\frac{2(n+1)}{b}-\frac{2a^2}{b^3}
       -\frac{a'b'}{a}-(2n-1)\frac{(b')^2}{b}+b'f'.
 \label{eq:prelim-bode}
\end{align}
\end{subequations}

For smooth extension across the zero section, we require \(a\) to extend oddly and $b,f$ to
extend evenly across $t=0$ and use the normalization in \eqref{eq:initial-data}. Now the  remaining initial parameter is $\mu=f^{\pr\pr}(0)$.
Appleton proved a local existence and uniqueness result of the system \eqref{eq:prelim-ode} with respect to the parameter $\mu$.

\begin{prop}[{\cite[Theorem~2.1, Appendix B]{Appleton}}]
\label{prop:prelim-local}
For each $\mu\in\mathbb R$, there exists $\varepsilon>0$ such that the following holds:

\begin{enumerate}
    \item  There exists a unique smooth solution $(f_\mu,a_\mu,b_\mu):(-\var,\var)\to \mathbb{R}^3$ to \eqref{eq:prelim-ode} with
\eqref{eq:initial-data} and $f_\mu''(0)=\mu$. Moreover, this solution is also analytic;
    \item  The function $a_\mu$ is an odd function and
$b_\mu,f_\mu$ are even functions;
    \item  The solution $(f_\mu,a_\mu,b_\mu)$ depends analytically on $\mu$.
\end{enumerate} 
\end{prop}

\begin{rem}\label{remark:smu def}
    We denote the maximal interval on which $a_\mu,b_\mu>0$ for $t>0$ by
$[0,t_\mu)$. By Proposition \ref{prop:prelim-local} we know that $t_\mu>0$. A solution on this interval defines a complete metric on
$\mathcal L_{k}$ precisely when $t_\mu=\infty$.
When the parameter is fixed, we omit the subscript $\mu$.
\end{rem}

To find parameters producing complete solutions, the following lemma tells us that it suffices to consider $\mu\leq 0$, which is a consequence of soliton identities.

\begin{lem}[{\cite[Lemmas~3.3--3.6]{Appleton}}]
\label{lem:prelim-identities}
Let $R$ denote the scalar curvature of a soliton solution of \eqref{eq:prelim-ode} with normalization \eqref{eq:initial-data} and $f^{\pr\pr}(0)=\mu$. Then
\begin{equation}\label{eq:prelim-hamilton}
 R+(f')^2=R(0)=-2\mu,
\end{equation}
and
\begin{equation}\label{eq:prelim-potential}
 f''+\left(\frac{a'}a+2n\frac{b'}b\right)f'-(f')^2=2\mu.
\end{equation}
Moreover,
\[
 \Delta_fR=-2|\operatorname{Ric}|^2,
 \qquad \Delta_f:=\Delta-\langle\nabla f,\nabla\cdot\rangle.
\]
If the solution is complete, then
\begin{equation}\label{eq:prelim-complete-bound}
 \mu\leq0,\qquad |f'|\leq\sqrt{-2\mu}.
\end{equation}
\end{lem}

Next, Appleton used soliton equations \eqref{eq:prelim-ode} to further analyze what parameters will give us complete solutions. First, combining the initial data \eqref{eq:initial-data}, ODE system \eqref{eq:prelim-ode} and the parity of $a,b,f$, we have the following Taylor expansions of $a,b,f$ around $t=0$. 
\begin{equation}\label{eq:prelim-taylor}
 \begin{aligned}
 a(t)&=\frac{(n+1)k}{p} t+
       \frac{(n+1)k}{6p}\bigl(\mu-2n(n+1)\bigr)t^3+O(t^5),\\
 b(t)&=1+\frac{n+1}{2}t^2+O(t^4),\\
 f(t)&=\frac{\mu}{2}t^2+O(t^4).
 \end{aligned}
\end{equation}
Indeed, the limit of \eqref{eq:prelim-bode} at $t=0$ gives
$b''(0)=n+1$, and \eqref{eq:prelim-fode} then gives
$a'''(0)={(n+1)k}\bigl(\mu-2n(n+1)\bigr)/p$.

Appleton introduced the ratio quantity:
\[
 Q:=\frac{a}{b}.
\]
Then by the initial data \eqref{eq:initial-data},
\begin{align}\label{eq:Q-initial}
    Q(0)=0, \qquad Q^{\pr}(0)=(n+1)k/p,   
\end{align}
and the soliton equations imply
\begin{equation}\label{eq:prelim-qode}
 Q''=\left(f'-(2n+1)\frac{b'}{b}\right)Q'
       +\frac{2(n+1)}{b^2}(Q^3-Q).
\end{equation}
In particular, every critical point of $Q$ with $0<Q<1$ is a strict local
maximum, whereas every critical point with $Q>1$ is a strict local minimum
\cite[equation~(3.1) and Lemma~3.1]{Appleton}. 

By analyzing the ODEs of $a,b,Q$, Appleton obtained the following monotonicity result.

\begin{prop}[{\cite[Lemmas~4.1, 4.2, 4.4 and Corollary~4.3]{Appleton}}]
\label{prop:prelim-monotonicity}
Suppose $(f,a,b)$ is a solution to the soliton equations \eqref{eq:prelim-ode} with \eqref{eq:initial-data} and $f^{\pr\pr}(0)=\mu\leq0$. Then $a'(t)>0$ for $0<t<t_\mu$~(see Remark \ref{remark:smu def} for the definition of $t_\mu$).
On any interval $[0,S)\subset[0,t_\mu)$ on which
$Q<\sqrt{n+1}$, we also have $b'(t)>0$ for $0<t<S$.
If $\mu<0$, then
\begin{equation}\label{eq:prelim-monotonicity}
 f'(t)<0,\qquad f''(t)<0,\qquad R'(t)<0
 \qquad(0<t<t_\mu).
\end{equation}
On an interval where $Q<\sqrt{n+1}$, we additionally have
$f'(t)\geq-\sqrt{-2\mu}$.
If $\mu=0$, then $f\equiv0$ and the corresponding metric is Ricci-flat
wherever it is defined.
\end{prop}

In particular, for a complete solution with $f^{\pr\pr}(0)=\mu<0$, the limit
\begin{equation}\label{eq:prelim-c}
 c:=-\lim_{t\to\infty}f'(t)
 \quad\text{exists and belongs to }(0,\sqrt{-2\mu}].
\end{equation}

The ratio $Q$ controls both continuation and the possible geometries at
infinity, as shown by the following two results of Appleton.

\begin{prop}[{\cite[Lemmas~5.3, 6.10 and Theorem~6.1]{Appleton}}]
\label{prop:prelim-continuation}
Suppose $(f,a,b)$ is a solution to the soliton equations \eqref{eq:prelim-ode} with \eqref{eq:initial-data} and $f^{\pr\pr}(0)=\mu\leq0$. Let $S<\infty$. If a solution is defined on $[0,S)$
and satisfies $Q<\sqrt{n+1}$ on $[0,S)$, then it extends past $S$.
Every complete solution satisfies $Q\leq1$.
Furthermore, if $Q(t_0)>1$ at some point of the maximal
interval of some solution with $\mu=f^{\pr\pr}(0)<0$, then $t_\mu<\infty$.
\end{prop}

\begin{rem}\label{rem:complete solution mu<0}
    For $f^{\pr\pr}(0)=\mu<0$, it follows that completeness is equivalent to the condition
$Q\leq1$ throughout the maximal interval. On a complete solution the
inequality is in fact strict for $t>0$: if $Q(t_0)=1$, then $Q'(t_0)=0$,
and uniqueness for \eqref{eq:prelim-qode} on $t>0$ would force
$Q\equiv1$, contrary to $Q(0)=0$. In particular, by Proposition \ref{prop:prelim-monotonicity}, $a',b'>0$ on every
complete non-Ricci-flat solution.
\end{rem}

\begin{theo}[{\cite[Theorem~6.1, Corollary~6.7 and Lemmas~6.4, 6.8, 6.9]{Appleton}}]
\label{thm:prelim-asymptotics}
Let $(f,a,b)$ be a complete solution to the soliton equations \eqref{eq:prelim-ode} with \eqref{eq:initial-data} and $f^{\pr\pr}(0)=\mu<0$.
Then
\begin{equation}\label{eq:prelim-dichotomy}
 Q_\infty:=\lim_{t\to\infty}Q(t)\in\{0,1\},
 \qquad \lim_{t\to\infty}Q'(t)=0.
\end{equation}
Let $c$ be the constant defined in \eqref{eq:prelim-c}. Then the two alternatives are as follows.
If $Q_\infty=0$, there exists $a_\infty\in(0,\infty)$ such that as $t\to+\infty$,
\begin{equation}\label{eq:prelim-collapsed}
 \begin{aligned}
 a(t)&\longrightarrow a_\infty,\\
 b(t)&\sim\sqrt{\frac{4(n+1)}c}\,\sqrt{t},\qquad
 b'(t)\sim\frac12\sqrt{\frac{4(n+1)}c}\,t^{-1/2}
 \end{aligned}
\end{equation}
 and the corresponding complete steady soliton is collapsed. If $Q_\infty=1$, then as $t\to+\infty$,
\begin{equation}\label{eq:prelim-noncollapsed}
 a(t)\sim b(t)\sim\sqrt{\frac{4n}c}\,\sqrt{t}
\end{equation}
and the corresponding steady soliton is $\kappa$-noncollapsed for some $\kappa>0$.
\end{theo}

For later reference, we record two elementary consequences of these
results. 

\begin{cor}
\label{cor:prelim-slope}
Every complete solution to the soliton equations \eqref{eq:prelim-ode} with \eqref{eq:initial-data} and $f^{\pr\pr}(0)=\mu<0$ satisfies
\begin{equation}\label{eq:prelim-slope}
 \lim_{t\to\infty}R(t)=0,\qquad
 c=\sqrt{-2\mu}.
\end{equation}
\end{cor}

\begin{proof}
Let $\mathcal J(t)=a(t)b^{2n}(t)$.
Theorem~\ref{thm:prelim-asymptotics} gives
$\log\mathcal J(t)=O(\log t)$. By Remark \ref{rem:complete solution mu<0}, $a^{\pr},b^\pr>0$ and hence $\mathcal J'>0$.
Integrating the trace equation
$R=-\Del f=-f''-(\log\mathcal J)'f'$ from a fixed $t_0>0$ to $t>t_0$ yields
\begin{align}\label{eq:integral R}
 \int_{t_0}^t R(s)\,ds
 =f'(t_0)-f'(t)
   -\int_{t_0}^t f'(s)(\log\mathcal J)'(s)\,ds
 =O(\log t),   
\end{align}
where we have used $\mathcal{J}^\pr>0$ and the bound $|f^{\pr}|\leq\sqrt{-2\mu}$ from \eqref{eq:prelim-complete-bound}. By \eqref{eq:prelim-monotonicity}, $R$ is decreasing. Since $R$ is nonnegative for a complete steady soliton~(cf.~\cite{Chen}), its limit must be zero by \eqref{eq:integral R}. Letting $t\to+\infty$ in \eqref{eq:prelim-hamilton} gives $c^2=-2\mu$.
\end{proof}

\begin{cor}\label{cor:prelim-ratio}
For every complete noncollapsed  solution to the soliton equations \eqref{eq:prelim-ode} with \eqref{eq:initial-data} and $f^{\pr\pr}(0)=\mu<0$,
\begin{equation}\label{eq:prelim-ratio-properties}
 0<Q(t)<1,\qquad Q'(t)>0,\qquad Q''(t)<0
 \qquad(t>0).
\end{equation}
In addition,
\begin{equation}\label{eq:prelim-ratio-bound}
 0<Q'(t)\leq\frac{(n+1)k}{p}\frac{e^{f(t)}}{b(t)^{2n+1}}.
\end{equation}
Consequently, $Q$ is a smooth global coordinate with range $(0,1)$ on the
radial interval $(0,\infty)$.
\end{cor}

\begin{proof}
By Remark \ref{rem:complete solution mu<0} and Theorem \ref{thm:prelim-asymptotics}, for a complete noncollapsed solution with $\mu<0$, we have $0<Q<1$ and $Q_\infty=\lim_{t\to\infty}Q(t)=1$.  Note that by \eqref{eq:prelim-qode} every critical point of $Q$ in $0<Q<1$ is a strict local maximum.
Since $Q'(0)=(n+1)k/p>0$, a first critical point would make $Q$ strictly
decreasing thereafter and a subsequent return to increasing behavior would
require a local minimum. This contradicts $Q_\infty=1$.
Thus $Q'>0$, and \eqref{eq:prelim-qode} gives $Q''<0$ because
$f'<0$, $b'>0$, and $0<Q<1$.
Finally, \eqref{eq:prelim-qode} can be written as
\begin{align}\label{eq:Q exp ode}
    \bigl(e^{-f}b^{2n+1}Q'\bigr)'
 =2(n+1)e^{-f}b^{2n-1}(Q^3-Q)<0.
\end{align}
Since $(e^{-f}b^{2n+1}Q^{\pr})(0)=(n+1)k/p$, integrating \eqref{eq:Q exp ode}
proves \eqref{eq:prelim-ratio-bound}.
\end{proof}

\begin{rem}\label{rem:k-range}
The preceding results, from Proposition \ref{prop:prelim-local} to Corollary \ref{cor:prelim-ratio},
hold for every positive integer $k$, subject to the hypotheses
of each statement. Indeed, the integer $k$ enters the reduced
ODE system only through the smoothness condition
$a'(0)=(n+1)k/p$, and no assumption $k>p$ is needed for
these results. The restriction $k>p$ is imposed in the
existence theorem below and in the uniqueness argument
of Section 3.
\end{rem}

Next, we recall Appleton’s construction of a noncollapsed parameter. Let
\[
 \mathcal C:=\{\mu\leq0:t_\mu=\infty\}
\]
be the set of complete normalized parameters.
Appleton proved the following theorem.

\begin{theo}[{\cite[Theorems~5.1, 7.2, Remark~5.5 and Lemma~7.1]{Appleton}}]
\label{thm:prelim-existence}
Fix the data of $(B^{2n},J,\widehat g)$ and $\mathcal{L}_k$ as before, with $k>p$ and the ansatz \eqref{eq:ansatz}. All sufficiently negative values of $\mu$ in Proposition \ref{prop:prelim-local} give complete collapsed solutions to the soliton equations \eqref{eq:prelim-ode} with \eqref{eq:initial-data} and $f^{\pr\pr}(0)=\mu$. The solution with $\mu=f^{\pr\pr}(0)=0$ satisfies $Q(s_1)>1$ for some finite $s_1>0$ and is incomplete. Let
\begin{equation}\label{eq:prelim-critical}
 \mu_{\mathrm A}:=
 \sup\bigl\{\nu\leq0:(-\infty,\nu]\subset\mathcal C\bigr\}.
\end{equation}
Then $-\infty<\mu_{\mathrm A}<0$,
$(-\infty,\mu_{\mathrm A}]\subset\mathcal C$ and the solution with
$\mu=f^{\pr\pr}(0)=\mu_{\mathrm A}$ satisfies $Q_\infty=1$. In particular, it defines a
complete noncollapsed steady gradient Ricci soliton on $\mathcal{L}_{k}$.
\end{theo}

The definition in \eqref{eq:prelim-critical} uses the endpoint of the initial
interval of complete parameters, rather than assuming in advance that
$\mathcal C$ is an interval. Likewise, the existence theorem does not
assert that every parameter below $\mu_{\mathrm A}$ is collapsed, or
that every parameter above $\mu_{\mathrm A}$ is incomplete; these are
part of the  trichotomy to be established in the next section.


\section{Proof of Theorem \ref{thm:conjecture}}\label{sec:comparison-quantities}

In this section we will give a complete proof of Theorem \ref{thm:conjecture}, solving a conjecture of Appleton. In particular, this will give us uniqueness of $U(2)$-invariant noncollapsed cohomogeneity one steady solitons on the total space of complex line bundles $\mathcal O(-k)$ of $\mathbb{CP}^1$ for $k\geq 3$. This uniqueness result will help us to classify four dimensional cohomogeneity one noncollapsed steady solitons.

The critical point for proving Theorem \ref{thm:conjecture} is the quantity $v$, defined in \eqref{eq:comparison-v-definition}. The rough idea is that we suppose that there are two parameters $\mu_1<\mu_2$ such that the soliton solutions determined by them by Proposition \ref{prop:prelim-local} are complete and noncollapsed. Then the corresponding $v_i,i=1,2$ must satisfy $v_1(0)=v_2(0)$ and $(v_2-v_1)_x>0$~(here $x$ is the variable defined in \eqref{eq:def-x} and subscript $x$ means differentiate with respect to $x$). Moreover, we can compute the second-order differential equation of $v_i$ to show that there must be $(v_2-v_1)_{x}>0$ on the entire domain of $x$. This will give us a positive lower bound of $v_2-v_1$. However, using the asymptotic analysis of Appleton's construction in Section 2, we can show that $v_i\to 0$ at infinity, which gives us a contradiction. This contradiction forces that there can only exist one noncollapsed complete parameter.

Throughout this section, $(f,a,b)$ denotes a complete noncollapsed solution in the class of Appleton's construction fixed in Section~2. Note that Appleton's existence theorem applies when $k>p$~(see Theorem \ref{thm:prelim-existence}). Thus
\begin{align}\label{eq:alpha n+1}
    \al:= a'(0)=\frac{(n+1)k}{p}>n+1
\end{align}
and Proposition \ref{prop:prelim-monotonicity}, Remark \ref{rem:complete solution mu<0} and Corollary \ref{cor:prelim-ratio} give
\begin{align}\label{eq:Qabf noncollapsed}
    0<Q<1,\qquad Q'>0,\qquad a',b'>0,\qquad f'<0
 \quad\text{on }(0,\infty).
\end{align}
We use the coordinate
\begin{align}\label{eq:def-x}
    x=Q^2\in(0,1),\qquad
 \partial_x=\frac{1}{2QQ'}\partial_t.
\end{align}
Primes continue to denote differentiation with respect to $t$; derivatives with respect to $x$ are indicated by subscripts. As in Section~2, we write $c=-\lim_{t\to\infty}f'(t)$. Then $c>0$ by \eqref{eq:prelim-c}.

Before using $v$ to get a contradiction, we need to introduce two auxiliary quantities, which we will denote by $m$ and $\mathcal{B}_n$. They will play an important role in the analysis of second-order differential equation of $v$. First, the quantity $m$ is defined by
\begin{equation}\label{eq:comparison-m-definition}
 m:=\frac{b'}{Q}=\frac{bb'}{a}.
\end{equation}
The parity of $a$ and $b$ implies that $m$ extends smoothly and evenly across $t=0$. Since $b''(0)=n+1$ and $Q'(0)=\alpha$, we have
\begin{align}\label{eq:m-initial}
    m(0)=\frac{n+1}{\alpha}=\frac pk<1.
\end{align}
Moreover, since $x=\al^2t^2+O(t^4)$, \(m\) extends smoothly as a function of \(x\).

By the computations in Appendix~A~(cf.~Lemma~\ref{lem:comparison-m-evolution}), $m$ satisfies the differential equations:
\begin{align}
    &m'=\frac{2(n+1-Q^2)}{bQ}
 -\left(2n\frac{b'}b-f'+2\frac{Q'}Q\right)m,\label{eq:comparison-m-first-order}\\
 &m''+\left((2n-1)\frac{b'}b-f'+3\frac{Q'}Q\right)m'
 =\frac{2(m^2-1)}{b^2}\left(nmQ^2+2bQ'\right).\label{eq:comparison-m-second-order}
\end{align}
and the Taylor expansion:
\begin{align}\label{eq:comparison-m-bolt}
    m(t)=\frac{n+1}{\alpha}
 +\frac{\alpha}{2}\left(\frac{(n+1)^2}{\alpha^2}-1\right)t^2
 +O(t^4).
\end{align}

Using the ODEs of $m$, we can show $m$ is strictly decreasing for $t\in (0,+\infty)$, hence giving an upper bound of $m$.

\begin{lem}\label{lem:comparison-m-bound}
For every $t>0$,
\begin{equation}\label{eq:comparison-m-bound}
 m'(t)<0,\qquad 0<m(t)<\frac{n+1}{\alpha}=\frac pk<1.
\end{equation}
In particular, if the base is $\mathbb{CP}^n$, then $\alpha=k$, this gives $m<(n+1)/k$.
\end{lem}

\begin{proof}
By \eqref{eq:Qabf noncollapsed}, $b^\pr>0$ for $t\in (0,+\infty)$. Then $m>0$ follows from the definition of $m$ in \eqref{eq:comparison-m-definition} and $Q=a/b>0$. Since $\alpha>n+1$~(see \eqref{eq:alpha n+1}), the coefficient of $t^2$ in \eqref{eq:comparison-m-bolt} is negative. Hence $m'<0$ for all sufficiently small positive $t$.

Note that by \eqref{eq:m-initial}, we know the initial value of $m$ satisfies $m(0)<1$. Suppose $t_0>0$ were the first point at which $m'(t_0)=0$. Then 
\[
 0<m(t_0)<m(0)<1,\qquad m''(t_0)\geq0.
\]
However, \eqref{eq:comparison-m-second-order} gives
\[
 m''(t_0)=\frac{2\bigl(m(t_0)^2-1\bigr)}{b(t_0)^2}
 \bigl(nm(t_0)Q(t_0)^2+2b(t_0)Q'(t_0)\bigr)<0,
\]
a contradiction. Thus $m'<0$ everywhere, and the asserted bound follows.
\end{proof}

Next, we introduce the quantity $\mathcal B_n$ defined by
\begin{equation}\label{eq:comparison-B-definition}
 \mathcal B_n
 :=-\frac{2(n+1)(1-Q^2)}{Q}(b^2Q')'
   +\bigl(nm^2-3n-2\bigr)b^2(Q')^2.
\end{equation}
This expression involves only the warping functions and their derivatives.  The quotient equation \eqref{eq:prelim-qode} of $Q$ gives the equivalent identity
\begin{equation}\label{eq:comparison-B-potential}
 \begin{aligned}
 \mathcal B_n
 ={}&2(n+1)(1-Q^2)
 \left[\frac{bQ'}Q\bigl((2n-1)b'-bf'\bigr)
       +2(n+1)(1-Q^2)\right]\\
 &+\bigl(nm^2-3n-2\bigr)b^2(Q')^2.
 \end{aligned}
\end{equation}
We will show that $\mathcal B_n$ is actually positive for $t\in(0,+\infty)$. The idea is by first determining the sign near infinity and then excluding a nonpositive value at finite distance with the help of ODE satisfied by $\mathcal B_n$.

\begin{lem}\label{lem:comparison-B-infinity}
As $t\to\infty$,
\begin{equation}\label{eq:comparison-end-limits}
 b'(t)\longrightarrow0,\qquad
 \frac{Q'(t)}{1-Q(t)}\longrightarrow c,\qquad
 \frac{\mathcal B_n(t)}{b(t)^2Q'(t)^2}\longrightarrow n+2.
\end{equation}
In particular, $\mathcal B_n>0$ for all sufficiently large $t$.
\end{lem}

\begin{proof}
For a complete noncollapsed solution, by Theorem \ref{thm:prelim-asymptotics} and \eqref{eq:prelim-c}, we know that as $t\to+\infty$,
\begin{align}\label{eq:asy abQf}
    a,b\to+\infty,\quad Q\to 1,\quad f^{\pr}\to -c<0. 
\end{align}
The equation for $b$ implies, for all sufficiently large $t$,
\[
 b''\leq\frac{2(n+1)}b-\frac c2b'.
\]
 Then for a sufficiently large $S$ and any $t\geq S$, we have
 \begin{align*}
     0<b^\pr(t)\leq e^{-\frac{c}{2}(t-S)}b^{\pr}(S)+2(n+1)\int_S^t e^{-\frac{c}{2}(t-\tau)}\frac{ d\tau}{b(\tau)}.
 \end{align*}
Since $b\to+\infty$ as $t\to+\infty$, this implies that
\begin{align}\label{eq:bprime-limit}
    \lim_{t\to+\infty}b^{\pr}(t)=0.
\end{align}

Let $\delta=1-Q$. By Corollary \ref{cor:prelim-ratio}, as $t\to+\infty$, $Q^{\pr}\to 0$. Integrating the quotient equation \eqref{eq:prelim-qode} from $t$ to infinity gives
\begin{equation}\label{eq:comparison-end-integral}
 Q'(t)=\int_t^\infty
 \left[\left((2n+1)\frac{b'}b-f'\right)Q'
       +\frac{2(n+1)Q(1+Q)}{b^2}\delta\right](\tau)\,d\tau.
\end{equation}
Moreover, the integral converges by the exponential estimates of Corollary~\ref{cor:prelim-ratio}. By \eqref{eq:asy abQf} and \eqref{eq:bprime-limit}, as $\tau\to+\infty$,
\begin{align*}
    (2n+1)\frac{b^\pr}{b}(\tau)-f^\pr(\tau) \longrightarrow c,\quad 0\leq \frac{2(n+1)Q(1+Q)}{b^2}(\tau)\longrightarrow 0
\end{align*}
 Since $Q^{\pr}>0$, $\delta$ is decreasing. Fix $0<\varepsilon<c$. For sufficiently large $t$, \eqref{eq:comparison-end-integral} therefore gives
\[
 -\delta^{\pr}(t)=Q'(t)\geq(c-\varepsilon)\delta(t).
\]
Consequently, for sufficiently large $t$,
\[
 \delta(s)\leq\delta(t)e^{-(c-\varepsilon)(s-t)}
 \quad\text{for }s\geq t.
\]
Using this estimate in \eqref{eq:comparison-end-integral}, and making $t$ larger if necessary, gives
\[
 Q'(t)\leq
 \left(c+\varepsilon+\frac{\varepsilon}{c-\varepsilon}\right)\delta(t).
\]
Letting $\varepsilon\downarrow0$ proves $Q'/\delta\to c$ as $t\to+\infty$.

Finally, divide \eqref{eq:comparison-B-potential} by $b^2(Q')^2$ to obtain
\[
 \begin{aligned}
 \frac{\mathcal B_n}{b^2(Q')^2}
 ={}&\frac{2(n+1)(1-Q^2)}{QQ'}
      \left((2n-1)\frac{b'}b-f'\right)\\
 &+\frac{4(n+1)^2(1-Q^2)^2}{b^2(Q')^2}
   +nm^2-(3n+2).
 \end{aligned}
\]
By the definition of $m$ in \eqref{eq:comparison-m-definition}, $m=b'/Q\to0$ as $t\to+\infty$. By the previous convergence $Q^{\pr}/(1-Q)\to c$ we have just established and \eqref{eq:asy abQf}--\eqref{eq:bprime-limit}, as $t\to +\infty$,
\begin{align*}
    \frac{2(n+1)(1-Q^2)}{QQ'}
      \left((2n-1)\frac{b'}b-f^{\pr} \right) &\longrightarrow 4(n+1),\\
      \frac{4(n+1)^2(1-Q^2)^2}{b^2(Q')^2}&\longrightarrow 0.
\end{align*}
This proves the last limit in \eqref{eq:comparison-end-limits}.
\end{proof}

The following elementary algebraic lemma will be used to determine the sign of $(\mathcal{B}_n)_x$ in Lemma \ref{lem:Bn-zero-crossing}.

\begin{lem}\label{lem:comparison-polynomial}
For every integer $n\geq1$, let $\mathscr P_n$ be the polynomial defined by
\begin{equation}\label{eq:comparison-polynomial}
 \begin{aligned}
 \mathscr P_n(s,z)={}&(n+2)z^2\\
 &+\left[2ns^2+(4n^2-6n-12)s
          +2(n+1)(5n+6)\right]z\\
 &+2n(2n-1)s^3+4(n^3-n^2-2n+1)s^2\\
 &+2(n+1)(10n^2+n-10)s+16(n+1)^3.
 \end{aligned}
\end{equation}
Then $\mathscr P_n(s,z)>0$ for any $s,z\geq 0$.
\end{lem}

\begin{proof}
For $n\geq3$, all displayed coefficients are nonnegative and the constant term is positive. For $n=1$, the polynomial is
\[
 \mathscr P_1(s,z)
 =3z^2+(2s^2-14s+44)z
   +2s\bigl((s-1)^2+1\bigr)+128.
\]
Its coefficient $2s^2-14s+44=2(s-7/2)^2+39/2$ is positive. For $n=2$,
\[
 \mathscr P_2(s,z)
 =4z^2+\bigl(4(s-1)^2+92\bigr)z
   +12s^3+4s^2+192s+432>0.
\]
This proves the claim in every dimension.
\end{proof}

Let
\[
 \vartheta(x)=\frac{bb'Q'}{Q(1-Q^2)},\qquad
 \zeta(x)=-\frac{b^2f'Q'}{Q(1-Q^2)},
\]
where $x=Q^2$. Both functions are positive on $(0,1)$ by \eqref{eq:Qabf noncollapsed}.
If $\mathcal B_n(x_0)=0$ for some $x_0\in(0,1)$, then at $x=x_0$, computations in Appendix~A~(cf.~Lemma~\ref{lem:Bn-zero-identities}) yield the following identity: 
\begin{equation}
 \frac{(3n+2)b^2(Q')^2}{(1-x)^2}
 =n\vartheta^2+2(n+1)(2n-1)\vartheta
   +2(n+1)\zeta+4(n+1)^2
\end{equation}
and a differential equation at $x_0$:
\begin{equation}\label{eq:Bn-zero-derivative}
\begin{aligned}
 (\mathcal B_n)_x
 ={}&2(1-x)(n\vartheta-n-1)
       \frac{n+(1-x)(1-\vartheta)}{x}\\
 &-\frac{(1-x)^3}{(3n+2)b^2(Q')^2}
   \Bigl[(n+1)\mathscr P_n(\vartheta,\zeta)
       +n(3n+2)\vartheta^2(2n\vartheta+\zeta)\Bigr],
\end{aligned}
\end{equation}
where $\mathscr P_n$ is the polynomial defined in~\eqref{eq:comparison-polynomial}. Moreover, Lemma~\ref{lem:Bn-zero-identities} also tells us that throughout $(0,1)$,
\begin{equation}\label{eq:mx-theta-zeta}
 (bQ')m_x
 =\frac{n+(1-x)(1-\vartheta)}{x}
 -\frac{(1-x)^2\vartheta(2n\vartheta+\zeta)}{2b^2(Q')^2}.
\end{equation}

With the above preparations, we can determine the direction in which $\mathcal{B}_n$ crosses at a zero. 

\begin{lem}\label{lem:Bn-zero-crossing}
If $\mathcal B_n(x_0)=0$ for some $x_0\in(0,1)$, then
\[
 (\mathcal B_n)_x(x_0)<0.
\]
\end{lem}

\begin{proof}
All quantities below are evaluated at $x=x_0$. We distinguish two cases.

Suppose first that $n\vartheta\leq n+1$. Then
\[
 n+(1-x)(1-\vartheta)
 \geq n-\frac{1-x}{n}>0.
\]
Thus the first term on the right-hand side of
\eqref{eq:Bn-zero-derivative} is nonpositive. Its second term is strictly
negative because $\mathscr P_n(\vartheta,\zeta)>0$ by Lemma~\ref{lem:comparison-polynomial}.
Consequently, $(\mathcal B_n)_x<0$.

Suppose instead that $n\vartheta>n+1$. Since $m_x<0$ by Lemma~\ref{lem:comparison-m-bound},
\eqref{eq:mx-theta-zeta} implies
\[
 \frac{n+(1-x)(1-\vartheta)}{x}
 <\frac{(1-x)^2\vartheta(2n\vartheta+\zeta)}{2b^2(Q')^2}.
\]
Multiplying this inequality by the positive number
$2(1-x)(n\vartheta-n-1)$ and substituting into
\eqref{eq:Bn-zero-derivative}, we obtain
\[
\begin{aligned}
 (\mathcal B_n)_x
 &<-\frac{(n+1)(1-x)^3}{(3n+2)b^2(Q')^2}
   \Bigl[\mathscr P_n(\vartheta,\zeta)
          +(3n+2)\vartheta(2n\vartheta+\zeta)\Bigr]\\
 &<0.
\end{aligned}
\]
This proves the assertion in both cases.
\end{proof}

Now we can prove the positivity of the quantity $\mathcal{B}_n$.

\begin{lem}\label{lem:comparison-B-positive}
Every complete noncollapsed soliton under consideration satisfies
\begin{equation}\label{eq:comparison-B-positive}
 \mathcal B_n(t)>0\qquad\text{for every }t>0.
\end{equation}
\end{lem}

\begin{proof}
By Lemma~\ref{lem:comparison-B-infinity}, $\mathcal B_n>0$ for $x$ sufficiently close to $1$. Suppose $\mathcal B_n(x_1)\leq0$ for some $x_1\in(0,1)$, and choose $x_2\in(x_1,1)$ with $\mathcal B_n(x_2)>0$. The last zero $x_0$ of $\mathcal B_n$ in $[x_1,x_2]$ has positive values immediately to its right, and hence $(\mathcal B_n)_x(x_0)\geq0$. This contradicts Lemma~\ref{lem:Bn-zero-crossing}. 
\end{proof}

We now define the third quantity $v$ by
\begin{equation}\label{eq:comparison-v-definition}
 v(x):=b(t)Q'(t),\qquad x=Q(t)^2.
\end{equation}
This is a positive function on $(0,1)$. Since \(bQ'\) is a smooth even function of \(t\) and \(x\) is a smooth function of \(t^2\) with derivative \(\alpha^2\) at zero, \(v\) extends smoothly to \(x=0\). The two quantities $m$ and $v$ determine the change of the warping functions in the $x$ coordinate through
\begin{equation}\label{eq:comparison-reconstruction}
 \frac{b_x}{b}=\frac{m}{2v},\qquad
 a=\sqrt{x}\,b,\qquad
 \frac{dt}{dx}=\frac{b}{2\sqrt{x}\,v}.
\end{equation}
Thus the normalization $b(0)=1$ fixes the scale in the reconstruction.

We can compute the second-order differential equation of $v$, which will be important for comparing two noncollapsed solutions.

\begin{lem}\label{lem:comparison-v-equation}
The function $v$ satisfies
\begin{equation}\label{eq:comparison-v-second-order}
 2v^2v_{xx}
 =\bigl(2(n+1)(1-x)+2mv\bigr)v_x
  +(nm^2+3n+2)v+2(n+1)(1-x)m.
\end{equation}
Moreover,
\begin{equation}\label{eq:comparison-v-potential}
 vv_x+nmv+(n+1)(1-x)
 =\frac{b^2f'Q'}{2Q}<0,
\end{equation}
and
\begin{equation}\label{eq:comparison-B-in-v}
 \mathcal B_n
 =-4(n+1)(1-x)vv_x-2(n+1)(1-x)mv
   +(nm^2-3n-2)v^2.
\end{equation}
In particular, $v_x<0$, and its endpoint data are
\begin{equation}\label{eq:comparison-v-endpoints}
 v(0)=\alpha,\qquad
 v_x(0)=\frac{\mu-2(n+1)^2}{2\alpha},\qquad
 \lim_{x\uparrow1}v(x)=0.
\end{equation}
\end{lem}

\begin{proof}
By equation \eqref{eq:prelim-qode}, we can compute
\begin{align}\label{eq:v2-diffx}
    (v^2)_x
 =\frac{bQ'}Q\bigl(bf'-2nb'\bigr)-2(n+1)(1-x).
\end{align}
Since $b'=mQ$, this proves \eqref{eq:comparison-v-potential}. Also,
\[
 \frac{(b^2Q')'}Q=(v^2)_x+mv,
\]
hence \eqref{eq:comparison-B-definition} becomes \eqref{eq:comparison-B-in-v}.

For the second-order equation, rewriting \eqref{eq:comparison-m-first-order} gives
\begin{equation}\label{eq:comparison-m-x}
 m_x=\frac{n+1-x}{xv}-\frac mx
      +\frac{m v_x}{v}+\frac{(n+1)(1-x)m}{v^2}.
\end{equation}
On the other hand, \eqref{eq:comparison-v-potential} is equivalent to
\begin{align}\label{eq:bfprime-eq}
    bf'=2\sqrt{x}\left(v_x+nm+\frac{(n+1)(1-x)}v\right).
\end{align}
Differentiating this identity with respect to $x$ gives us
\begin{equation}\label{eq:right-hand-side-diff}
    \begin{aligned}
        (bf^{\pr})_x=&\frac{1}{\sqrt{x}}\left( v_x+nm+\frac{(n+1)(1-x)}{v} \right)\\
    &+2\sqrt{x}\left( v_{xx}+nm_x-\frac{n+1}{v}-\frac{(n+1)(1-x)}{v^2}v_x \right).
    \end{aligned}
\end{equation}
On the other hand, using
\eqref{eq:prelim-fode} and \eqref{eq:comparison-reconstruction}, we can get
\begin{align}\label{eq:left-hand-side-diff}
    (bf')_x=\frac{m\,bf'}{2v}
         +\frac{v_x}{\sqrt{x}}
         +(2n+1)\sqrt{x}\,m_x
         +\frac{(n+1)m}{\sqrt{x}}.
\end{align}
Equating \eqref{eq:right-hand-side-diff} and \eqref{eq:left-hand-side-diff} while using \eqref{eq:bfprime-eq} to substitute $bf^\pr$ and \eqref{eq:comparison-m-x} to substitute $m_x$, we obtain
\[
 2v_{xx}
 =2\left(\frac{(n+1)(1-x)}{v^2}+\frac mv\right)v_x
  +\frac{nm^2+3n+2}{v}
  +\frac{2(n+1)(1-x)m}{v^2},
\]
which is \eqref{eq:comparison-v-second-order}.

Since $f'<0$, \eqref{eq:comparison-v-potential} implies $v_x<0$. By the Taylor expansion \eqref{eq:prelim-taylor}, we can get
\[
 v(t)=\alpha+\frac{\alpha}{2}
          \bigl(\mu-2(n+1)^2\bigr)t^2+O(t^4),
 \qquad x=\alpha^2t^2+O(t^4).
\]
This proves the initial data in \eqref{eq:comparison-v-endpoints}. Finally, the limit at $x=1$ in \eqref{eq:comparison-v-endpoints} follows from the exponential decay of $Q'$ in \eqref{eq:prelim-ratio-bound} and the square-root growth of $b$ established in \eqref{eq:prelim-noncollapsed}.
\end{proof}

Now we set about the proof of uniqueness of noncollapsed complete parameters. For two complete noncollapsed solutions with parameters $\mu_1,\mu_2$, we use subscripts $i=1,2$ for the corresponding quantities $v_i$, $m_i$, and $\mathcal B_{n,i}$. In the following, all comparisons are made at the same value of $x\in(0,1)$ while the corresponding radial coordinates on the two solitons need not agree. By \eqref{eq:Qabf noncollapsed}, Lemma \ref{lem:comparison-m-bound} and Lemma \ref{lem:comparison-B-positive}, 
\[
 v_i>0,\qquad 0<m_i<1,\qquad \mathcal B_{n,i}>0,
\]
and
\[
 v_i(v_i)_x+nm_iv_i+(n+1)(1-x)<0.
\]
Moreover, \eqref{eq:m-initial} and \eqref{eq:comparison-v-endpoints} give us
\begin{equation}\label{eq:uniqueness-profile-data}
 \begin{aligned}
 &v_i(0)=\alpha,\qquad
 (v_i)_x(0)=\frac{\mu_i-2(n+1)^2}{2\alpha},\\
 &m_i(0)=\frac{n+1}{\alpha},\qquad\lim_{x\uparrow1}v_i(x)=0.
 \end{aligned}
\end{equation}
By Lemma \ref{lem:comparison-v-equation}, we can write
\begin{align}\label{eq:v-second order-Phi}
    v_{xx}=\Phi_n(x,v,v_x,m),    
\end{align}
 where
\[
 \Phi_n(x,v,s,m)
 =\left(\frac{(n+1)(1-x)}{v^2}+\frac mv\right)s
  +\frac{nm^2+3n+2}{2v}
  +\frac{(n+1)(1-x)m}{v^2}.
\]
Holding the other arguments fixed, along the soliton, direct differentiation while using \eqref{eq:comparison-B-positive} and \eqref{eq:comparison-v-potential} gives
 \begin{align}
 \partial_m\Phi_n
 &=\frac{vv_x+nmv+(n+1)(1-x)}{v^2}<0,\label{eq:Phi partial m}\\
 \partial_v\Phi_n
 &=\frac{\mathcal B_n}{2v^4}
   -\frac{m}{v^3}\bigl(vv_x+nmv+(n+1)(1-x)\bigr)>0.\label{eq:Phi partial v}
 \end{align}

First, we derive an ordering of the ratios $m_i/v_i$.

\begin{lem}\label{lem:uniqueness-ratio}
Let $X\in(0,1)$. If $v_2(x)>v_1(x)$ for every $x\in(0,X]$, then
\[
 \frac{m_2(x)}{v_2(x)}<\frac{m_1(x)}{v_1(x)}
 \qquad\text{for every }x\in(0,X].
\]
\end{lem}

\begin{proof}
Equation~\eqref{eq:comparison-m-x} gives
\begin{equation}\label{eq:uniqueness-ratio-identity}
 \left(\frac{xm}{v}\right)_x
 =\frac{n+1-x+(n+1)x(1-x)m/v}{v^2}.
\end{equation}
Let
\[
 r(x):=x\left(\frac{m_2(x)}{v_2(x)}-\frac{m_1(x)}{v_1(x)}\right).
\]
Subtracting \eqref{eq:uniqueness-ratio-identity} for the two solutions while using $0<x<1$ and $v_2>v_1$ yields
\begin{equation}\label{eq:uniqueness-ratio-difference}
 \begin{aligned}
 r_x-\frac{(n+1)(1-x)}{v_2^2}r
 ={}&\left(n+1-x+(n+1)x(1-x)\frac{m_1}{v_1}\right)\\
 &\qquad\cdot\left(\frac1{v_2^2}-\frac1{v_1^2}\right)<0
 \qquad(0<x\leq X).
 \end{aligned}
\end{equation}
The coefficient of $r$ is continuous up to $x=0$, since $v_2(0)=\alpha>0$, and $r(0)=0$. Multiplying by a positive integrating factor and integrating from $0$ to $x$ gives $r(x)<0$. This is the desired ratio ordering.
\end{proof}

Next, we compare the second-order derivatives at a critical point of $v_2-v_1$.

\begin{lem}\label{lem:uniqueness-contact}
Suppose there exists $x_0\in(0,1)$ such that at $x_0$ we have
\[
 v_2>v_1,\qquad
 \frac{m_2}{v_2}\leq\frac{m_1}{v_1},\qquad
 (v_2)_x=(v_1)_x.
\]
Then $(v_2)_{xx}(x_0)>(v_1)_{xx}(x_0)$.
\end{lem}

\begin{proof}
Throughout this proof, all profile quantities are evaluated at $x_0$, which we abbreviate by $x$. Let
\begin{align}\label{eq:t h def}
    s=(v_1)_x=(v_2)_x,\qquad h=\frac{m_2}{v_2}.
\end{align}
By \eqref{eq:v-second order-Phi},
\begin{align*}
    (v_i)_{xx}=\Phi_n(x,v_i,s,m_i )
\end{align*}
We compare the right-hand sides by first decreasing $m$ at fixed $v_1$ and then increasing $v$ while keeping $m/v=h$ fixed.

Since $hv_1\leq m_1$, for every $m\in[hv_1,m_1]$ we can use \eqref{eq:Phi partial m} and \eqref{eq:comparison-v-potential} to get
\begin{align*}
    \partial_m\Phi_n(x,v_1,s,m)
 =&\frac{v_1s+nmv_1+(n+1)(1-x)}{v_1^2}\\
 \leq&\frac{v_1s+nm_1v_1+(n+1)(1-x)}{v_1^2}<0.
\end{align*}
It follows that
\begin{align}\label{eq:Phi hv1 m1}
    \Phi_n(x,v_1,s,hv_1)\geq\Phi_n(x,v_1,s,m_1).
\end{align}

For the second step, using \eqref{eq:Phi partial m}, \eqref{eq:Phi partial v} and a direct differentiation give
\begin{equation}\label{eq:uniqueness-path-derivative}
 \begin{aligned}
 2z^3\frac{d}{dz}\Phi_n(x,z,s,hz)
 ={}&-4(n+1)(1-x)s\\
 &-\bigl(2(n+1)(1-x)h+3n+2\bigr)z
   +nh^2z^3.
 \end{aligned}
\end{equation}
The derivative of the right-hand side of \eqref{eq:uniqueness-path-derivative} with respect to $z$ is
\[
 -2(n+1)(1-x)h-(3n+2)+3nh^2z^2<0
 \qquad(v_1\leq z\leq v_2),
\]
since $hz\leq hv_2=m_2<1$ by Lemma \ref{lem:comparison-m-bound}. Thus the right-hand side of \eqref{eq:uniqueness-path-derivative} is decreasing in $z$ on this interval. Its value at $z=v_2$ is
\[
 \frac{-4(n+1)(1-x)v_2s-2(n+1)(1-x)m_2v_2
          +(nm_2^2-3n-2)v_2^2}{v_2}
 =\frac{\mathcal B_{n,2}}{v_2}>0,
\]
where we have used \eqref{eq:comparison-B-in-v} to get the first equality and Lemma \ref{lem:comparison-B-positive} to get the positivity of $\mathcal{B}_n$.
Consequently, $z\mapsto\Phi_n(x,z,s,hz)$ is strictly increasing on $[v_1,v_2]$. Note that $hv_2=m_2$ by our definition of $h$ in \eqref{eq:t h def}, using this increasing property and \eqref{eq:Phi hv1 m1} we conclude that
\[
 \begin{aligned}
 (v_2)_{xx}
 &=\Phi_n(x,v_2,s,m_2)\\
 &>\Phi_n(x,v_1,s,hv_1)\\
 &\geq\Phi_n(x,v_1,s,m_1)
 =(v_1)_{xx}.
 \end{aligned}
\]
This proves the lemma.
\end{proof}

Now we can prove the uniqueness of the normalized noncollapsed parameter.

\begin{lem}\label{lem:uniqueness-parameter}
At most one parameter $\mu=f^{\pr\pr}(0)<0$ in the normalized  family in Appleton's construction produces a complete noncollapsed soliton.
\end{lem}

\begin{proof}
    Suppose there are two parameters $\mu_1<\mu_2<0$ such that the two solutions $(a_i,b_i,f_i)$ to the soliton equation \eqref{eq:prelim-ode} with \eqref{eq:initial-data} and $f_i^{\pr\pr}(0)=\mu_i$ are both complete and noncollapsed, where $i=1,2$. Then we can define $v_i$ for each soliton solution by \eqref{eq:comparison-v-definition}. By \eqref{eq:uniqueness-profile-data},
    \begin{align*}
        (v_2-v_1)(0)=0,\qquad
 (v_2-v_1)_x(0)=\frac{\mu_2-\mu_1}{2\alpha}>0.
    \end{align*}
    This implies that $(v_2-v_1)_x(x)>0$ for sufficiently small $x\in (0,1)$. Suppose $(v_2-v_1)_x$ attains zero at some point of $(0,1)$, then we can take
    \begin{align*}
        x_0:=\inf\lbrace x\in (0,1):(v_2-v_1)_x(x)=0 \rbrace.
    \end{align*}
    By the above discussion, we have $x_0\in (0,1)$ and
    \begin{equation}\label{eq:v-contra}
        \begin{aligned}
            &(v_2-v_1)_x(x)>0 \qquad \text{on}~(0,x_0),\\
        &(v_2-v_1)(x)>0\qquad\text{on}~(0,x_0),\\
        &(v_2-v_1)_x(x_0)=0.
        \end{aligned}
    \end{equation}
    From this we can get
    \begin{align*}
        (v_2-v_1)_{xx}(x_0)\leq  0.
    \end{align*}
    By Lemma~\ref{lem:uniqueness-ratio},  $m_2/v_2<m_1/v_1$ at $x_0$. However, combining this with \eqref{eq:v-contra}, we can see $(v_2-v_1)_{xx}(x_0)>0$ by Lemma \ref{lem:uniqueness-contact}. This is a contradiction. Hence, 
    \begin{align*}
        (v_2-v_1)_x>0\qquad \text{on}~(0,1).
    \end{align*}
    Now we get that $v_2-v_1$ is strictly increasing on $(0,1)$ with $(v_2-v_1)(0)=0$. This implies that
    \begin{align*}
        \liminf_{x\uparrow 1} (v_2(x)-v_1(x))\geq v_2(1/2)-v_1(1/2)>0.
    \end{align*}
    On the other hand, \eqref{eq:uniqueness-profile-data} tells us that $v_i(x)\to 0$ as $x\uparrow 1$ for each $i=1,2$. This gives us a contradiction. Therefore, such two complete noncollapsed parameters cannot exist. Finally, by the local uniqueness result in Proposition \ref{prop:prelim-local} and the ordinary differential equation uniqueness away from the zero section, we conclude that a complete noncollapsed parameter gives us a unique soliton metric. This completes the proof of the lemma.
\end{proof}

We now return to the entire normalized family $(f_\mu,a_\mu,b_\mu)$. For $\mu<0$, the results recalled in Section~2 imply that the maximal solution is either incomplete or complete with $Q_\infty\in \lbrace 0,1\rbrace$. We first show that the collapsed and incomplete alternatives are open conditions on the parameter. 

\begin{lem}\label{lem:uniqueness-collapsed-open}
The set of parameters $\mu=f^{\pr\pr}(0)<0$ producing complete collapsed solitons in Appleton's construction is open in $(-\infty,0)$.
\end{lem}

\begin{proof}
    Let $\mu_0=f^{\pr\pr}(0)<0$ be a complete collapsed parameter in Appleton's construction. By Remark \ref{rem:complete solution mu<0} and Theorem \ref{thm:prelim-asymptotics}, $Q(t)< 1$ on $(0,+\infty)$ and $\lim_{t\to\infty}Q_{\mu_0}(t)=0$. This implies that there exists $S\in (0,+\infty)$ such that 
    \begin{align*}
        0<Q_{\mu_0}(t)<1,\;\forall\,t\in [0,S],\qquad Q_{\mu_0}^\pr(S)<0.
    \end{align*}
    By the continuous dependence of the parameter, for parameters sufficiently close to $\mu_0$, the solutions of \eqref{eq:prelim-ode} also exist on $[0,S]$ and
    \begin{align*}
        0<Q_{\mu}(t)<1,\;\forall\,t\in [0,S],\qquad Q_{\mu}^\pr(S)<0.
    \end{align*}
    On the maximal existence interval, if $Q_\mu$ had a first critical point after $S$, then $Q_{\mu}<1$ at this point and hence \eqref{eq:prelim-qode} implies that $Q_\mu^{\pr\pr}<0$ at this point. This would give us a contradiction. Thus $Q_\mu<1$ on $(S,t_\mu)$. By Proposition \ref{prop:prelim-continuation}, we know that $t_\mu=\infty$ and such $\mu$ are complete parameters. Moreover, the above argument also implies that $Q_\mu^\pr<0$ on $(S,+\infty)$. Therefore,
    \begin{align*}
        \limsup_{t\to+\infty} Q_\mu(t)\leq Q_\mu(S)<1.
    \end{align*}
    Combining this with Theorem \ref{thm:prelim-asymptotics}, we can get that
    \begin{align*}
        \lim_{t\to+\infty} Q_\mu(t)=0.
    \end{align*}
    This shows that $\mu$ is also a complete collapsed parameter. Therefore, the set of complete collapsed parameters are open in $(-\infty,0)$.
\end{proof}

\begin{lem}\label{lem:uniqueness-incomplete-open}
The set of parameters $\mu=f^{\pr\pr}(0)<0$ producing incomplete solutions in Appleton's construction is open in $(-\infty,0)$. In addition, every negative parameter sufficiently close to zero produces an incomplete solution.
\end{lem}

\begin{proof}
    Let $\mu_0<0$ be an incomplete parameter in Appleton's construction. By Remark \ref{rem:complete solution mu<0}, there exists $S\in (0,+\infty)$ such that $Q_{\mu_0}(S)>1$. Using the continuous dependence of the parameter, this implies that $Q_\mu(S)>1$ for parameters $\mu$ sufficiently close to $\mu_0$. Now we can apply Remark \ref{rem:complete solution mu<0} again to get that solutions corresponding to such $\mu$ are incomplete. Therefore, the set of incomplete parameters is open in $(-\infty,0)$.

    For the zero parameter, by Theorem \ref{thm:prelim-existence}, the solution corresponding to this parameter satisfies $Q(s_1)>1$ for some finite $s_1>0$ and is incomplete. As above, we can use the continuous dependence property to get that every negative parameter sufficiently close to zero gives us an incomplete solution.
\end{proof}

At this stage, we can prove Theorem \ref{thm:conjecture}.

\begin{proof}[Proof of Theorem~\ref{thm:conjecture}]
By Theorem \ref{thm:prelim-existence} and Lemma \ref{lem:uniqueness-parameter}, there exists a unique complete noncollapsed parameter $\mu_*<0$. It suffices to show the trichotomy \eqref{miu}.

For this, we consider two connected open intervals $(-\infty,\mu_*)$ and $(\mu_*,0)$. By Lemma \ref{lem:uniqueness-collapsed-open} and \ref{lem:uniqueness-incomplete-open}, the sets of complete collapsed parameters and incomplete parameters are both open subsets of $(-\infty,\mu_*)$ and $(\mu_*,0)$. Moreover, since $\mu_*$ is the unique noncollapsed parameter, the complete collapsed parameters and incomplete parameters give us a partition of $(-\infty,\mu_*)$ and $(\mu_*,0)$. By connectedness, $(-\infty,\mu_*)$ contains only complete collapsed parameters or incomplete parameters, which is also true for $(\mu_*,0)$. By Theorem \ref{thm:prelim-existence}, there is a complete collapsed parameter lying in $(-\infty,\mu_*)$. This implies that $(-\infty,\mu_*)$ contains only complete and collapsed parameters. Using the same method, since there is an incomplete parameter in $(\mu_*,0)$ by Lemma \ref{lem:uniqueness-incomplete-open}, we can see $(\mu_*,0)$ contains only incomplete parameters. Finally, all $\mu\geq 0$ yield incomplete solution by Lemma \ref{lem:prelim-identities} and Theorem \ref{thm:prelim-existence}. This completes the proof of the trichotomy \eqref{miu}.
    
\end{proof}

\section{Linear full-curvature decay and finite fundamental group}
\label{sec:analytic}

Now we start the classification of cohomogeneity one $\kappa$-noncollapsed non-Ricci-flat steady gradient Ricci solitons in dimension four and five. In this section, we record two useful results in the classification of cohomogeneity one setting. The first result is an improvement of Theorem \ref{thm:cohom-est} from scalar curvature decay to the norm of the Riemann curvature tensor decay under the $\kappa$-noncollapsed condition.

\begin{prop}
\label{prop:full-linear}
Under the hypotheses of Theorem \ref{thm:cohom-est}, assume in addition that $(M,g)$ is $\kappa$-noncollapsed.  Then
\begin{equation}\label{eq:full-linear}
  |\Rm|(x)\le\frac{C}{1+t(x)},
\end{equation}
where $t(x)=d(x,Q)$ and  $Q=G/K$ is the singular orbit. In particular, the curvature is globally bounded.
\end{prop}

\begin{proof}
The norm $|\Rm|$ is constant on each homogeneous principal orbit.  If
\eqref{eq:full-linear} failed, there would be $t_i\to\infty$ with
$t_i|\Rm|(t_i)\to\infty$.  Let
\[
  I_i=\left[\frac{t_i}{2},\frac{3t_i}{2}\right],\qquad
  d_i(t)=\operatorname{dist}(t,\partial I_i),
\]
and choose $\bar t_i\in I_i$ such that
\begin{align*}
    |\Rm|(\bar{t}_i)d_i(\bar{t}_i)^2=\max_{t\in I_i} |\Rm|(t)d_i(t)^2    
\end{align*}
For simplicity, we denote
\[
  Q_i=|\Rm|(\bar t_i),\qquad D_i=d_i(\bar t_i).
\]
Since $|\Rm|(\bar{t}_i)d_i(\bar{t}_i)^2\geq |\Rm|(t_i)d_i(t_i)^2=|\Rm|(t_i)t_i^2/4$ and $D_i\leq t_i/2$,
\begin{equation}\label{eq:pointpick-scales}
    \begin{aligned}
        &Q_iD_i^2\to\infty, \qquad
  \sqrt{Q_i}\,D_i\to\infty,\\
  &Q_i D_i\geq \frac{|\Rm|(t_i)t_i^2}{4D_i}\geq \frac{1}{2}t_i|\Rm|(t_i)\to \infty,\\
  &Q_i D_i\leq Q_i \bar{t}_i\to \infty.
    \end{aligned}
\end{equation}
The maximality of choice of $\bar{t}_i$ gives
\begin{equation}\label{eq:pointpick-bound}
  |\Rm|(t)\le4Q_i
  \quad\text{whenever}\quad |t-\bar t_i|\le D_i/2.
\end{equation}
Choose $p_i$ on the orbit $t=\bar t_i$ and rescale the canonical flow by
\[
  g_i(s)=Q_i\,g(Q_i^{-1}s).
\]
The canonical diffeomorphisms move points in the static metric with speed at
most $\sqrt{C_0}$ by the third identity in \eqref{eq:soliton-identities}. Hence, for any
\begin{align*}
    (x,s)\in B_{g_i(s)}(p_i,A) \times [-A,A],
\end{align*}
we have
\begin{align*}
    |t(\Phi_{s/Q_i}(x))-\bar{t}_i|\leq &|t(\Phi_{s/Q_i}(x))-t(\Phi_{s/Q_i}(p_i))|+|t(\Phi_{s/Q_i}(p_i))-\bar{t}_i|\\
    \leq & d_g(\Phi_{s/Q_i}(x),\Phi_{s/Q_i}(p_i))+d_g(\Phi_{s/Q_i}(p_i),p_i)\\
    \leq & \frac{A}{\sqrt{Q_i}} +\frac{A \sqrt{C_0}}{Q_i}\leq \frac{D_i}{2}
\end{align*}
for sufficiently large $i$. Therefore, by \eqref{eq:pointpick-bound},
\[
  \sup_{|s|\le A}\;
  \sup_{x\in B_{g_i(s)}(p_i,A)}
  |\Rm_{g_i(s)}|_{g_i(s)}(x)
  \le 4
\]
for sufficiently large $i$.

The $\kappa$-noncollapsing hypothesis supplies a uniform injectivity-radius
lower bound at $p_i$.  By Hamilton compactness theorem, after passing to a subsequence, $(M,g_i(s),p_i)_{-\infty<s<+\infty}$ converges to an eternal Ricci flow
\[
  (M_\infty,g_\infty(s),p_\infty)
\]
in the Cheeger-Gromov sense. Moreover, by our construction, $|\Rm_{g_\infty(0)}|(p_\infty)=1$. For the scalar curvature at $p_\infty$, we can use \eqref{eq:scalar-linear} and \eqref{eq:pointpick-scales} to compute
\begin{align*}
    R_{g_\infty(0)}(p_\infty)=\lim_{i\to+\infty}R_{g_i(0)}(p_i)\leq \lim_{i\to+\infty}Q_i^{-1}\frac{C}{\bar{t}_i}=0.
\end{align*}
By the maximum principle~(cf.~\cite{Chen}), the ancient Ricci flow $(M_\infty,g_\infty(s))_{s\leq 0}$ is Ricci-flat.

Now we let $\gamma_i$ be the unit-speed normal geodesic with respect to $g_i(0)$ such that $\ga_i(0)=p_i$. Since the metric has the form $\md t^2+g_t$, every subsegment of $\ga_i$ is globally minimizing. By our choice of $p_i$, $d_{g_i(0)}(Q,p_i)=\sqrt{Q_i}\bar{t}_i$. This implies that the domain of $\ga_i$ is
\[
  \bigl(-\sqrt{Q_i}\,\bar t_i,\infty\bigr).
\]
By \eqref{eq:pointpick-scales}, the two sides of this interval tend to infinity. Hence, after passing to a subsequence, $\ga_i$ converges to a line $\ga_\infty$ through $p_\infty$. Since the limit is Ricci-flat, the Cheeger--Gromoll theorem produces an isometric translation
group along this line.

Now we claim that $M_\infty$ is homogeneous.  Let $x_\infty\in M_\infty$ and $x_i$ be a sequence of points converging to it. By this convergence, there exists a constant $C_1$ independent of $i$ such that
\begin{align*}
    d_{g_i(0)}(p_i,x_i)\leq C_1.
\end{align*}
In the global principal-part
coordinates we can write
\[
  x_i=(\bar t_i+\delta_iQ_i^{-1/2},y_i).
\]
 Let $q_i=(\bar t_i,y_i)$. Then $p_i$ and $q_i$ lie on the same homogeneous orbit, hence 
 \begin{align}\label{eq:piqi-bound}
     |\delta_i|=d_{g_i(0)}(p_i,q_i)=\sqrt{Q_i}|t(x_i)-\bar{t}_i|\leq d_{g_i(0)}(p_i,x_i)\leq C_1.
 \end{align}
Moreover, there is $a_i\in G$ with $a_i(p_i)=q_i$. Now we can get
\begin{align*}
    d_{g_i(0)}(p_i,a_i(p_i))=&d_{g_i(0)}(p_i,q_i)\\
    \leq& d_{g_i(0)}(p_i,x_i)+d_{g_i(0)}(x_i,q_i)\\
    \leq &2d_{g_i(0)}(p_i,x_i)\leq C_1.
\end{align*} 
This implies the isometries
$a_i$ converge to an isometry $a_\infty$ on $(M_\infty,g_\infty(0))$ after passing to a subsequence. Moreover, by \eqref{eq:piqi-bound} and after passing to another subsequence, we assume $\delta_i\to \delta_0$. Since the group
action preserves the normal coordinate, we can get
\[
  x_\infty=\lim_{i\to+\infty}x_i=\lim_{i\to+\infty}a_i(\ga_i(\delta_i))=a_\infty\!\left(\gamma_\infty(\delta_\infty)\right).
\]
This shows that a translation along the line followed by $a_\infty$ sends $p_\infty$ to
$x_\infty$. Therefore, the isometry group of $M_\infty$ is transitive.

By the theorem of
Alekseevskii--Kimel'fel'd \cite{AlekseevskiiKimelfeld}, every homogeneous Ricci-flat Riemannian manifold is flat.  This contradicts
$|\Rm|(p_\infty)=1$ and proves \eqref{eq:full-linear}.
\end{proof}

The next proposition converts the classical noncollapsing hypothesis into
the entropy condition.

\begin{prop}\label{prop:kappa-entropy}
Under the hypotheses of Proposition \ref{prop:full-linear}, the pointed Nash entropy of
the canonical flow is bounded from below.  Consequently condition
\begin{align}\label{NC condition}
    \nu(M,g)>-\infty
\end{align}
 holds and
\[
  |\pi_1(M)|<\infty.
\]
\end{prop}

\begin{proof}
   Let $\Phi_t$ be the flow generated by $\nabla f$ and $g(t):=\Phi_t^*g$. Then $(M,g(t))$ is the canonical Ricci flow of the steady soliton $(M,g)$. Fix $o\in M$ and let $d\nu_{o,0;-\tau}$ be the conjugate heat kernel measure.  For any $\tau>0$, let $(z_{\tau},-\tau)$ be an $H_{m+1}$-center of $(o,0)$, i.e.,
    \begin{align*}
        \int_M d_{g(-\tau)}^2(z_\tau,y)d\nu_{o,0;-\tau}(y)\leq H_{m+1}\tau.
    \end{align*}
Let $w_\tau=\Phi_{-\tau}(z_\tau)$, $\bar{\nu}_\tau=(\Phi_{-\tau})_*\nu_{o,0;-\tau}$. Since $\Phi_{-\tau}:(M,g(-\tau))\to (M,g)$ is an isometry,
\begin{align*}
    \int_Md_g^2(w_\tau,y)d\bar{\nu}_{\tau}(y)\leq H_{m+1}\tau.
\end{align*}
On the other hand, for any smooth function $\phi$ on $M$,
\begin{align*}
    \frac{d}{dt}\int_M\phi\circ \Phi_td\nu_{o,0;t}=&\int_M(\partial_t-\Del_{g(t)})\phi\circ\Phi_t d\nu_{o,0;t}\\    
    =&\int_M\lan \nabla f,\na \phi \ran\circ \Phi_t-\Del_{g(t)}(\phi\circ \Phi_t)d\nu_{o,0;t}.
\end{align*}
Hence, for $\tau=-t$ and $\phi=f$ we can get that
\begin{align*}
    \frac{d}{d\tau}\int_M fd\bar{\nu}_\tau=&-\int_M(|\na f|^2-\Del f )d\bar{\nu}_\tau\\
    =-&\int_M(|\na f|^2+R)d\bar{\nu}_\tau\\
    =&-1,
\end{align*}
where we have used the third equation of \eqref{eq:soliton-identities} in the last equality and we assume $C_0=1$ after a harmless normalization. This implies that
\begin{align*}
    \int_M f d\bar{\nu}_\tau=f(o)-\tau
\end{align*}
and hence
\begin{align*}
    |f(w_\tau)-f(o)+\tau|\leq& \int_M |f(y)-f(w_\tau)|d\bar{\nu}_\tau(y)\\
    \leq & \int_Md_g(y,w_\tau)d\bar{\nu}_\tau(y)\\
    \leq & \sqrt{H_{m+1}\tau}.
\end{align*}
Thus, $f(w_\tau)=-\tau+O(\sqrt{\tau})$ for large $\tau$. Since $|\na f|\leq 1$, we get that
\begin{align}\label{eq:bound f in ball}
    f(y)=-\tau+O(\sqrt{\tau}),\qquad\forall\,y\in B_g(w_\tau,\var \sqrt{\tau})
\end{align}
for sufficiently large $\tau$. By Theorem \ref{thm:cohom-est}, this implies that
\begin{align}\label{eq:distance ball}
    \frac{t(y)}{\tau}\longrightarrow 1,\qquad \forall\,y\in B_g(w_\tau,\var\sqrt{\tau}).
\end{align}

Now, fix a small $\varepsilon>0$, which will be specified below. Consider the ball $B_{g(-\tau)}(z_\tau,\varepsilon\sqrt{\tau})$, which is isometric to $B_g(w_\tau,\varepsilon\sqrt{\tau})$. By Proposition \ref{prop:full-linear} and \eqref{eq:distance ball}, there is a constant $C_1>0$ such that 
\begin{align*}
    \sup_{B_{g(-\tau)}(z_\tau,\varepsilon\sqrt{\tau})}|\mathrm{Rm}|_{g(-\tau)}\leq \frac{C_1}{\tau},\quad \forall\,0<\var<1.
\end{align*}
Choose $\var$ sufficiently small such that $C_1\var^2<1$. Then
\begin{align*}
    \frac{C_1}{\tau}\leq \frac{1}{\var^2\tau}=(\var \sqrt{\tau})^{-2}.
\end{align*}
Hence the $\kappa$-noncollapsing condition gives 
\begin{align*}
    \mathrm{Vol}_{g(-\tau)}B_{g(-\tau)}(z_\tau,\sqrt{
    \tau})\geq \mathrm{Vol}_{g(-\tau)}B_{g(-\tau)}(z_\tau,\var\sqrt{\tau})\geq \kappa \var^{m+1}\tau^{\frac{m+1}{2}}.
\end{align*}
By \cite[Corollary 1.11]{CMZ2023}, we have uniform Nash entropy bound $\mathcal{N}_{o,0}(\tau)\geq -\beta(\kappa,\var)$ for any $\tau>0$. Moreover, by \cite[Theorem 1.14]{CMZ2023}, this gives the bound of $\inf_{t\leq 0}\nu(g(t))$. Note that for a steady gradient Ricci soliton, $\nu(M,g)=\nu(g(t))$ for any $t\in \mathbb{R}$ and thus \eqref{NC condition} holds. Bamler's fundamental-group theorem for noncollapsed ancient
Ricci flows then yields $|\pi_1(M)|<\infty$
\cite{BamlerFundamentalGroup}.
\end{proof}

\section{Asymptotically cylindrical property in dimension four}\label{cylindrical}

In this section, we derive an asymptotically cylindrical property of $\kappa$-noncollapsed 4-dimensional cohomogeneity one steady gradient Ricci solitons. The asymptotic geometry will help us to analyze the asymptotic behavior of the metric components of an $SU(2)$-invariant metric in the next section.


The asymptotically cylindrical property is stated as follows.

\begin{prop}\label{prop: AC property}
    Let $(M,g,f)$ be a 4-dimensional complete non-Ricci-flat $\kappa$-noncollapsed steady gradient Ricci soliton. Suppose it is also a cohomogeneity one manifold invariant under a compact connected Lie group $G$-action. Then $(M,g)$ is asymptotically $\mathbb{S}^3/\Ga$-cylindrical. More precisely, for any sequence $p_i$ on $M$ tending to infinity,
    \begin{align*}
        (M,R(p_i)g,p_i)\longrightarrow (\mathcal{C}_\Ga,g_{cyl},p_\infty)
    \end{align*}
    in the smooth Cheeger-Gromov sense. Here $\mathcal C_\Gamma=(\mathbb{S}^3/\Gamma)\times\mathbb{R}$ with the finite group $\Gamma$ acting freely on $\mathbb{S}^3$ and $g_{cyl}$ is normalized to have scalar curvature 1.
\end{prop}

\begin{proof}
By Theorem \ref{thm:cohom-est} and \cite[Theorem~3]{Chan2019}, the norm of the curvature tensor is controlled by the scalar curvature, i.e.,
\begin{align}\label{eq:1}
    |{\rm Rm}|(x)\le CR(x),
\end{align}
where $C$ is independent of $x$. By Theorem \ref{thm:cohom-est},  $f$ is proper. By Proposition \ref{prop:full-linear}, the curvature decays linearly. Following Chan-Zhu's estimate in \cite{Chanzhu}, we get
\begin{align}\label{eq:2}
    |\nabla^k{\rm Rm}|(x)\cdot |f(x)|^{\frac{k+2}{2}}\le C(k),
\end{align}
where $C(k)$ is some constant may depend on $k$.
Let $\phi_s$ be the one-parameter diffeomorphism generated by $\nabla f$.  Then, $g(s)=\phi_{s}^{*}g$ satisfies the Ricci flow. By Theorem \ref{thm:cohom-est}, $\langle \nabla R, \nabla f\rangle\ge0$. Hence, $R(s)$ satisfies 
\begin{align}\label{eq:3}
    \frac{\partial R(s)}{\partial s}=\langle \nabla R, \nabla f\rangle\ge 0.
\end{align}
By (\ref{eq:1}),(\ref{eq:2}),(\ref{eq:3}) and the properness of $f$,  Deng-Zhu's dimension reduction arguments  in \cite{DZ20b} hold in cohomogeneity one case. Similar to   \cite[Lemma 6.2]{DZ20b}, for any $p_i\to \infty$, $(M,R(p_i)g(R^{-1}(p_i)s)),p_i)$ converge subsequentially to $(\Sigma\times \mathbb{R},g_{\Sigma}(s)+ds^2)$, where $(\Sigma,g_{\Sigma}(s))$ is a non-flat $\kappa$-noncollapsed ancient Ricci flow with backward type I estimate. Following the proof of   \cite[Theorem 5.1]{CDM}, $\Sigma=\mathbb{S}^3/\Gamma$ and $g_{\Sigma}(s)$ has constant sectional curvature for each $s$.
\end{proof}

\begin{rem}
 The argument in \cite[Proposition 3.8]{BCDMZ2022} can also be used to give another proof of Proposition \ref{prop: AC property} due to Proposition \ref{prop:kappa-entropy}. Although \cite[Proposition 3.8]{BCDMZ2022} is stated for steady soliton
singularity models, the neck argument used there depends only on local
properties of the canonical flow. 
\end{rem}

As a corollary, we can get the principal orbits are diffeomorphic to a spherical three-manifold.

\begin{cor}\label{coro: principal orbit}
    Let $(M,g,f)$ be a 4-dimensional complete, non-Ricci-flat, and $\kappa$-noncollapsed steady gradient Ricci soliton. Suppose it is also a cohomogeneity one manifold invariant under a compact connected Lie group $G$-action. Then the principal orbits are diffeomorphic to a spherical three-manifold.
\end{cor}

\begin{proof}
    Since $(M,g,f)$ has linear scalar curvature decay \eqref{eq:scalar-linear} and asymptotically $\mathbb{S}^3/\Ga$-cylindrical properties, one can show the level sets outside a compact set are diffeomorphic to a spherical three-manifold by Deng-Zhu's level set method(cf. \cite{DZ2020}). Note that the potential function $f$ is strictly decreasing(cf. \cite[Proposition 2.3]{BDW2015}) and $G$-invariant, each level set outside a compact set is exactly a principal orbit. Hence, the principal orbit is also diffeomorphic to a spherical three-manifold. 
\end{proof}

Since the principal orbit is diffeomorphic to a spherical three-manifold, the action group $G$ has a transitive action on this spherical three-manifold. We can use this property to reduce the action group to $SU(2)$. Before this, we need the following lifting lemma, which helps us to lift the action to the covering manifold.

\begin{lem}
\label{lemma: lift action}
Let $\widetilde M$ be connected and $\pi:(\widetilde M,\widetilde g)\longrightarrow (M,g)$ be a finite regular Riemannian covering. Let $G$ be a compact connected
Lie group acting smoothly, almost effectively, and isometrically on
$(M,g)$ with cohomogeneity $c$. Then there exist a compact connected Lie group $\widetilde G$ and a
finite covering homomorphism $\rho:\widetilde G\longrightarrow G$ such that $\widetilde G$ acts smoothly, almost effectively, and isometrically on $(\widetilde M,\widetilde g)$ with the same
cohomogeneity $c$. The actions are related by
\[
\pi(\widetilde g\cdot\widetilde x)
=
\rho(\widetilde g)\cdot\pi(\widetilde x)
\]
for every $\widetilde g\in\widetilde G$ and
$\widetilde x\in\widetilde M$. In particular, 

\begin{enumerate}[leftmargin=*]
\item if $G$ acts transitively on $M$, then $\widetilde G$ acts
transitively on $\widetilde M$;

\item if $G$ acts on $M$ with cohomogeneity one, then $\widetilde G$
acts on $\widetilde M$ with cohomogeneity one.
\end{enumerate}
\end{lem}

\begin{proof}
Let $\Gamma$ be the finite deck transformation group of $\pi$. Since the covering is regular,
\[
M\cong \widetilde M/\Gamma.
\]
Consider the group of lifts
\[
\widehat G
=
\left\{
(g,\widetilde g)\in
G\times\operatorname{Isom}(\widetilde M,\widetilde g):
\pi\circ\widetilde g=g\circ\pi
\right\}.
\]
By this definition, $\widehat G$ is a closed subgroup of $G\times \mathrm{Isom}(\widetilde M,\widetilde g)$, hence it is a Lie group. Since $G$ is connected and $\pi$ is regular, every element of $G$
admits a lift to $\widetilde M$ and two lifts differ by an element of
$\Gamma$. This implies that the projection
\[
q:\widehat G\longrightarrow G,
\qquad
q(g,\widetilde g)=g,
\]
is a smooth surjective homomorphism with finite kernel. In particular, $\widehat G$ is also a compact Lie group.

Let $\widetilde G:=\widehat G^{\,0}$ be the identity component. Since $G$ is connected, the image
$q(\widetilde G)$ is an open and closed subgroup of $G$. This forces that
\[
q(\widetilde G)=G.
\]
Thus
\[
\rho:=q|_{\widetilde G}:\widetilde G\longrightarrow G
\]
is a finite covering homomorphism. The action of $\widetilde G$ on
$\widetilde M$ is isometric, and the defining property of $\widehat G$
gives
\[
\pi(\widetilde g\cdot\widetilde x)
=
\rho(\widetilde g)\cdot\pi(\widetilde x).
\]

Let $\mathfrak{\widetilde G}, \mathfrak g$ be the Lie algebras of $\widetilde G$ and $G$. Let $\widetilde x\in\widetilde M$ and put $x=\pi(\widetilde x)$. If
$\widetilde X\in\mathfrak{\widetilde G}$ and
$X=d\rho_e(\widetilde X)\in\mathfrak g$, then
\[
d\pi_{\widetilde x}
\bigl(\widetilde X_{\widetilde x}^{\#}\bigr)
=
X_x^{\#}.
\]
Since $d\pi_{\widetilde x}$ and $d\rho_e$ are isomorphisms, it follows
that
\[
\dim\bigl(\widetilde G\cdot\widetilde x\bigr)
=
\dim\bigl(G\cdot x\bigr).
\]
Therefore the principal orbit dimensions of the two actions agree.
Since $\dim\widetilde M=\dim M$, the two actions have the same
cohomogeneity.

If an element of $\widetilde G$ acts trivially on $\widetilde M$, then
its image under $\rho$ acts trivially on $M$. Since the original action
is almost effective and $\rho$ has finite kernel, the kernel of the
lifted action is finite. Hence the lifted action is almost effective.
\end{proof}

\begin{lem}\label{lem:transitive}
Let $G$ be a compact connected Lie group acting almost effectively, transitively and isometrically on a spherical three-manifold  $\Sph^3/\Ga$. Then after passing to finite covers and finite ineffective quotients, the action is one of the following:
\[
  SU(2),\qquad U(2),\qquad SO(4).
\]
In particular, unless the acting group already contains a subgroup locally isomorphic to $U(2)$, it is locally isomorphic to $SU(2)$ and the  isotropy group is finite.
\end{lem}

\begin{proof}
 After quotienting by ineffective kernels, we may assume that $G$ acts transitively and effectively on $\mathbb{S}^3/\Ga$. By Lemma \ref{lemma: lift action}, there is a compact connected group $\tilde{G}$ acting transitively and almost effectively on $\mathbb{S}^3$ such that $\tilde{G}$ is a finite cover of $G$.  The classification of compact connected transitive transformation groups of spheres~\cite{MS1943} gives, up to finite ineffective kernels, the left $SU(2)$ action, the standard $U(2)$ action, and the full $SO(4)$ action.  Passing back to the spherical quotient changes only finite kernels and finite covering data.
\end{proof}

Moreover, we can reduce the isotropy group $K$ of the singular orbit to connected Lie groups. Note that in the following we use $\mathrm{Lie}(G)$ to denote the Lie algebra of a Lie group $G$

\begin{lem}\label{lem:connected-cover}
Let $(M^4,g)$ be an $SU(2)$-invariant cohomogeneity one manifold. Assume that the principal isotropy $H$ is finite and the singular orbit has codimension two. Let $K$ be its isotropy group and $K^0$ be the identity component of $K$. Then $K^0\cong S^1$, $K=K^0H$ and there is a finite Riemannian covering
\begin{equation}\label{eq:connected-cover}
  \widehat M=SU(2)\times_{K^0}D^2\longrightarrow
  M=SU(2)\times_KD^2
\end{equation}
such that $\widehat M$ is also an $SU(2)$-invariant cohomogeneity one manifold with principal isotropy group $H_0=H\cap K^0$.  If the lifted orbit metrics are biaxial with symmetry axis $\mathrm{Lie}(K^0)$, then the lifted metric admits an isometric action locally isomorphic to $U(2)$. 
\end{lem}

\begin{proof}
Since the singular orbit has codimension 2, $K/H\cong S^1$. Moreover, $H$ is finite, this implies that every component of $K$ meets $H$. Equivalently,
\[
  K=K^0H.
\]
Therefore, the natural map
\[
  K^0/(H\cap K^0)\longrightarrow K/H
\]
is a diffeomorphism.  This implies that the slice representation of $K$ on the normal disk restricts to the standard transitive circle representation of $K^0$ with kernel $H_0=H\cap K^0$. Consequently, the inclusion $K^0\subset K$ induces the finite covering of associated disk bundles in \eqref{eq:connected-cover}. On the principal part, it is the finite covering
\[
  SU(2)/H_0\longrightarrow SU(2)/H.
\]
Moreover, the invariant metric, the invariant potential lift and completeness are all preserved under this finite covering.

If the lifted orbit metrics are biaxial with symmetry axis $\mathrm{Lie}(K^0)$, then we can choose a Milnor frame such that $X_1$ generates $\mathrm{Lie}(K^0)$ and the lifted orbit metrics have the form
\[
  a(r)^2\omega_1^2+b(r)^2(\omega_2^2+\omega_3^2).
\]
From this we can see that such metrics are invariant under right multiplication by $K^0$.  Since $K^0$ is abelian, this right action is well defined on $SU(2)/H_0$ and extends over $SU(2)\times_{K^0}D^2$.  Together with left $SU(2)$ action it has Lie algebra
\[
  \mathfrak{su}(2)\oplus\mathbb R\cong\mathfrak u(2)
\]
and finite kernel.  If the original $K$ is connected, then $K=K^0$ and the same construction is defined directly on the original associated disk bundle and no covering is needed for the extra circle action.
\end{proof}

Consequently, to prove Theorem \ref{thm: four}, it suffices to treat the case in which $G$ is locally $SU(2)$.  We lift the acting group and, in the codimension-two case, pass to the finite connected-isotropy cover as in Lemma \ref{lem:connected-cover}.  On that cover the principal isotropy is finite and the singular isotropy is connected.

Let $K$ be the isotropy of the singular orbit on this cover.  Since $K/H$ is a positive-dimensional sphere, $K$ has positive-dimensional identity component. Moreover, since $K$ is the singular isotropy group, it must be closed. The next lemma shows that the only positive-dimensional closed connected proper subgroup of $SU(2)$ is a circle.

\begin{lem}\label{lemma: unique subgroup of su(2)}
    Let $K$ be a closed connected proper subgroup of $SU(2)$ with positive dimension. Then $K$ must be a circle. 
\end{lem}

\begin{proof}
    Let $\mathfrak{su}(2)$ be the Lie algebra of $SU(2)$ and as usual, $\mathrm{Lie}(K)$ be the Lie algebra of $K$.
    Since $\mathrm{dim}\,\mathfrak{su}(2)=3$, $\dim \mathrm{Lie}(K)=1$ or 2 by our assumption. We show that $\dim \mathrm{Lie}(K)=2$ is impossible.

    We use the identification of $\mathfrak{su}(2)$ with $\mathbb{R}^3$ to verify this. We choose a basis $E_{1},E_{2},E_3$ of $\mathfrak{su}(2)$ such that
    \begin{align*}
        [E_1,E_2]=2E_3,\quad[E_2,E_3]=2E_1,\quad [E_3,E_1]=2E_2.
    \end{align*}
    Under the linear identification
    \begin{align*}
        x_1E_1+x_2E_2+x_3E_3\longleftrightarrow (x_1,x_2,x_3)\in \mathbb{R}^3,
    \end{align*}
    the Lie bracket, is, up to the factor 2, the Euclidean cross product:
    \begin{align*}
        [X,Y]=2X\times Y.
    \end{align*}
    Now suppose $\mathrm{Lie}(K)$ is two dimensional and we choose linearly independent elements
    \begin{align*}
        X,Y\in\mathrm{Lie}(K).
    \end{align*}
    Since $\mathrm{Lie}(K)$ is a Lie subalgebra,
    \begin{align*}
        [X,Y]\in \mathrm{Lie}(K).
    \end{align*}
    However, $[X,Y]=2X\times Y$, which is nonzero and orthogonal to both $X$ and $Y$. Hence
    \begin{align*}
        [X,Y]\notin \mathrm{span}\, \lbrace X,Y \rbrace=\mathrm{Lie}(K),
    \end{align*}
    which is a contradiction. Consequently, $\dim K=1$ and since $K$ is compact, it must be a circle.
\end{proof}

In conclusion, to show Theorem \ref{thm: four}, there are only two possibilities need to be considered:
\begin{lem}\label{lem:four two cases}
    Let $(M^4,g,f)$ be an $SU(2)$-invariant non-Ricci-flat, $\kappa$-noncollpased cohomogeneity one steady gradient Ricci soliton. Then there is a finite covering  $(\widehat M,\widehat g,\widehat f)$ which is also an $SU(2)$ cohomogeneity one steady gradient Ricci soliton. Moreover, one of the following must hold for the singular isotropy group $K$ of $\widehat M$:
    \begin{enumerate}[label=\textup{(\alph*)}]
        \item $K=SU(2)$, so the singular orbit is a point and $K/H\cong\mathbb{S}^3$;
        \item $K\cong S^1$, so $K/H\cong S^1$ and the singular orbit has codimension two.
    \end{enumerate}
\end{lem}

\begin{proof}
    By Lemma \ref{lem:connected-cover}, we can lift the $SU(2)$ action to a finite covering $(\widehat M,\widehat g,\widehat f)$ such that the singular isotropy group $K$ is connected. By Proposition \ref{prop:double}, $K$ must have positive dimension. Then Lemma \ref{lemma: unique subgroup of su(2)} forces that the connected singular isotropy group $K$ must be either $SU(2)$ or the circle.
\end{proof}

\section{Proof of Theorem \ref{thm: four}}

In this section, we will assume $(M^4,g,f)$ is an $SU(2)$-invariant cohomogeneity one steady gradient Ricci soliton, which is complete, non-Ricci-flat and $\kappa$-noncollapsed. Using \cite[Proposition 2]{B2025}, by choosing a suitable Milnor frame $X_1,X_2,X_3$ on $SU(2)$ and after passing to a finite cover, the metric $g$ on the principal part can be written in the form: 

\begin{align}\label{su(2) invariant metric}
    g=dt^2+f_1^2(t)\omega_1\otimes \omega_1+f_2^2(t)\omega_2\otimes \omega_2+f_3^2(t)\omega_3\otimes\omega_3,
\end{align}
where $f_i:(0,+\infty)\to (0,+\infty)$ are smooth functions and $\omega_i$ are one forms dual to $X_i$. In the following, after a normalization, we may assume the Milnor frame $X_1,X_2,X_3$ satisfies
\begin{align}\label{eq:milnor}
    [X_1,X_2]=X_3,\quad [X_2,X_3]=X_1,\quad [X_3,X_1]=X_2.    
\end{align}

As in \cite{B2025}, we can define the following quantities:

\begin{align*}
    L_i=\frac{f_i'}{f_i},
  \qquad
  R_i=\frac{f_i}{f_jf_k},
  \qquad
  \{i,j,k\}=\{1,2,3\},
\end{align*}
and
\begin{equation}\label{eq:xidef}
  \xi=L_1+L_2+L_3-f'.
\end{equation}
Then the steady soliton equation is equivalent to the following ODE system~(cf. \cite{B2025}).

\begin{align}
  \xi^\prime&=-L_1^2-L_2^2-L_3^2 ,\label{eq:xi}\\
    L_i'&=-\xi L_i+\frac12R_i^2-\frac12(R_j-R_k)^2,
  \label{eq:Li}\\
  R_i'&=R_i(L_i-L_j-L_k),\qquad \{i,j,k\}=\{1,2,3\}.
  \label{eq:Ri}
\end{align}

For our use, for any $i,j\in\lbrace 1,2,3 \rbrace$, we define
\begin{align}\label{eq:def-Bij-qij}
    B_{ij}=L_i-L_j,
  \qquad
  q_{ij}=\frac{R_i}{R_j}=\frac{f_i^2}{f_j^2},
\end{align}
then
\begin{align}
    B_{ij}'&=-\xi B_{ij}+(R_i-R_j)(R_i+R_j-R_k),\label{eq:B}\\
  q_{ij}'&=2q_{ij}B_{ij}.
  \label{eq:q}
\end{align}

In the remainder of this section, $i,j,k$ will always denote indices such that $\lbrace i,j,k \rbrace=\lbrace 1,2,3 \rbrace$.

The proof of Theorem \ref{thm: four} relies on a series of ratio maximum principles, which start with the following lemma.

\begin{lem}\label{lem:max-implication}
Let $q_{ij},q_{kj}$ be defined as in \eqref{eq:def-Bij-qij}. Suppose $q_{ij}$ has an interior local maximum at $t_0$ and $q_{ij}(t_0)>1$.  Then
\begin{equation}\label{eq:triangle-reverse}
  R_k(t_0)\ge R_i(t_0)+R_j(t_0).
\end{equation}
Consequently,
\begin{equation}\label{eq:larger-ratio}
  q_{kj}(t_0)\ge q_{ij}(t_0)+1.
\end{equation}
\end{lem}

\begin{proof}
Since $q_{ij}$ has an interior local maximum at $t_0$, $q_{ij}^{\pr}(t_0)=0$.  Then equation \eqref{eq:q} and $q_{ij}(t_0)>1$ imply that $B_{ij}(t_0)=0$.  Differentiating \eqref{eq:q} yields
\[
  q_{ij}''=2q_{ij}'B_{ij}+2q_{ij}B_{ij}'.
\]
At $t_0$ we have
\[
  0\ge q_{ij}''(t_0)=2q_{ij}(t_0)B_{ij}'(t_0),
\]
and hence $B_{ij}'(t_0)\le0$.  Since $B_{ij}(t_0)=0$, equation \eqref{eq:B} becomes
\[
  0\geq B_{ij}'(t_0)
  =(R_i-R_j)(R_i+R_j-R_k)(t_0).
\]
Since $q_{ij}(t_0)>1$, we have $R_i(t_0)>R_j(t_0)$ by the definition of $q_{ij}$. This implies  $R_i(t_0)+R_j(t_0)-R_k(t_0)\le0$, which gives us \eqref{eq:triangle-reverse}.  Division by $R_j>0$ gives \eqref{eq:larger-ratio}.
\end{proof}

The cyclic implication yields a convenient global lemma.

\begin{lem}\label{lem:two-cycle}
Let $q_{ij},q_{kj}$ be defined as in \eqref{eq:def-Bij-qij}.  Suppose 
\begin{equation}\label{eq:end-limsup}
  \limsup_{t\to0}q_{ij}(t)\le1,
  \quad
  \limsup_{t\to\infty}q_{ij}(t)\le1,
\end{equation}
with the same two inequalities for $q_{kj}$.  Then
\[
  q_{ij}(t)\le1
  \qquad\text{for every }t>0.
\]
\end{lem}

\begin{proof}
For contradiction, we assume that
\[
  M_{ij}=\sup_{t>0}q_{ij}(t)>1.
\]
By \eqref{eq:end-limsup}, the supremum is attained at an interior point $t_{ij}$.  Then Lemma \ref{lem:max-implication} tells us that
\[
  q_{kj}(t_{ij})\ge M_{ij}+1.
\]
This implies that
\[
  M_{kj}:=\sup_{t>0}q_{kj}(t)\ge M_{ij}+1>1.
\]
Since $q_{kj}$ has the same endpoint assumption as in \eqref{eq:end-limsup}, $M_{kj}$ is attained at an interior point $t_{kj}$.  Using Lemma \ref{lem:max-implication} to the pair $(k,j)$ with remaining index $i$, we obtain
\[
  q_{ij}(t_{kj})\ge M_{kj}+1.
\]
Hence,
\[
  M_{ij}\ge M_{kj}+1\ge M_{ij}+2,
\]
which is a contradiction.
\end{proof}

\begin{lem}\label{lem:bolt-rigidity}
Let $q_{ij},q_{kj}$ be defined as in \eqref{eq:def-Bij-qij}. Suppose
\begin{align}
  &\lim_{t\to0}q_{23}(t)=1,\qquad
  \lim_{t\to0}q_{12}(t)=\lim_{t\to0}q_{13}(t)=0,
  \label{eq:bolt-left}\\
  &\lim_{t\to\infty}q_{ij}(t)=1\qquad
  \text{for every }i\ne j.
  \label{eq:bolt-right}
\end{align}
Then
\[
  R_2\equiv R_3,
  \qquad
  f_2\equiv f_3.
\]
\end{lem}

\begin{proof}
Since $q_{13}$ and $q_{23}$ satisfy the boundary conditions \eqref{eq:bolt-left}--\eqref{eq:bolt-right}, we can apply Lemma \ref{lem:two-cycle} to get $q_{23}(t)\le1$ for every $t>0$. By the definition of $q_{ij}$ in \eqref{eq:def-Bij-qij},
\begin{align*}
    \lim_{t\to 0}q_{32}(t)=\lim_{t\to 0}q_{23}^{-1}(t)=1.
\end{align*}
Combining this with \eqref{eq:bolt-left} and \eqref{eq:bolt-right}, we can also apply Lemma \ref{lem:two-cycle} to $q_{32}$ and $q_{12}$ to obtain $q_{32}(t)\le1$ for any $t>0$. Since $q_{23}=q_{32}^{-1}$, we finally obtain $q_{23}\equiv1$. Then the claimed equalities follow from the definition of $q_{ij}$ in \eqref{eq:def-Bij-qij}.
\end{proof}

By the asymptotically $\mathbb{S}^3/\Ga$-cylindrical in Section~5, we can get the following asymptotic behavior of the warping functions $f_i$.

\begin{lem}\label{lemma: asy warping functions}
    Let $(M,g,f)$ be a four dimensional complete $\kappa$-noncollapsed non-Ricci-flat steady gradient Ricci soliton. Suppose $(M,g)$ admits an $SU(2)$ cohomogeneity one action with only one singular orbit and the metric $g$ on the principal part of the form \eqref{su(2) invariant metric} with the normalization \eqref{eq:milnor}. Then 
    \begin{align}\label{eq:f-ratio-round-infinity}
        \frac{f_i}{f_j}(t)\to 1,\quad (t\to+\infty)
    \end{align}
    for every $i,j\in \lbrace 1,2,3 \rbrace$.
\end{lem}

\begin{proof}
     Let $g_t=f_1^2(t)\omega_1^2+f_2^2(t)\omega_2^2+f_3^2(t)\omega_3^2$ be the induced metrics on principal orbits(i.e., regular level sets). Using the normalization of the Milnor frame \eqref{eq:milnor} and the curvature formulas of homogeneous metrics(cf. \cite[Chapter 7]{B2007}),
    \begin{align}\label{Ricci formula orbits}
        \mathrm{Ric}_{g_t}(X_i/f_i,X_i/f_i)=\frac{f_i^4-(f_j^2-f_k^2)^2}{2(f_1f_2f_3)^2},\quad \lbrace i,j,k \rbrace=\lbrace 1,2,3\rbrace.
    \end{align}
    By the asymptotically $\mathbb{S}^3/\Ga$-cylindrical property proved in Proposition \ref{prop: AC property} and Deng-Zhu's level set method, one can show that after rescaling, $g_t$ converges to the round metric on $\mathbb{S}^3/\Ga$. Then we get that
    \begin{align}\label{Ric ratio to 1}
        \frac{\mathrm{Ric}_{g_t}(X_i/f_i,X_i/f_i)}{\mathrm{Ric}_{g_t}(X_j/f_j,X_j/f_j)}\to 1,\quad (t\to+\infty)
    \end{align}
    for any $i,j\in\lbrace 1,2,3 \rbrace$. Then the lemma follows from \eqref{Ricci formula orbits} and \eqref{Ric ratio to 1}. Indeed, let
    \begin{align*}
        f_{\max}(t):=\max\lbrace f_i(t),1\leq i\leq 3 \rbrace,\qquad f_{\min}(t):=\min \lbrace f_i(t),1\leq i\leq 3 \rbrace.
    \end{align*}
    Suppose \eqref{eq:f-ratio-round-infinity} is false, then there exist a sequence $t_m$ tending to infinity and a constant $\ell>1$ such that
    \begin{align*}
        \frac{f_{\max}(t_m)}{f_{\min}(t_m)}\geq \ell,\qquad \forall\,m\geq 1.
    \end{align*}
    After passing to a subsequence and reordering the indices, we may assume that
    \begin{align*}
        f_1(t_m)=f_{\max}(t_m),\quad f_3(t_m)=f_{\min}(t_m),\qquad \forall\, m\geq 1.
    \end{align*}
    Then we have
    \begin{align}\label{eq:contra-f1f3}
        \frac{f_1(t_m)}{f_3(t_m)}\geq \ell>1 ,\qquad \forall\, m\geq 1.
    \end{align}
    Note that at $t=t_m$,
    \begin{align}\label{eq:Ric ratio}
        1-\frac{\mathrm{Ric}_{g_t}(X_3/f_3,X_3/f_3)}{\mathrm{Ric}_{g_t}(X_1/f_1,X_1/f_1)}=\frac{2(f_1^2-f_3^2)}{f_1^2+f_2^2-f_3^2}\geq \frac{f_1^2-f_3^2}{f_1^2}
    \end{align}
    Since
    \begin{align*}
        \frac{\mathrm{Ric}_{g_t}(X_3/f_3,X_3/f_3)}{\mathrm{Ric}_{g_t}(X_1/f_1,X_1/f_1)}\to 1,\quad (t\to+\infty),
    \end{align*}
    \eqref{eq:Ric ratio} implies that
    \begin{align*}
        \liminf_{m\to+\infty}\frac{f_3(t_m)^2}{f_1(t_m)^2}\geq 1,
    \end{align*}
    which contradicts \eqref{eq:contra-f1f3}.
\end{proof}

Now we get down to the proof of Theorem \ref{thm: four}. By Lemma \ref{lem:four two cases}, after passing to a finite cover and finite ineffective quotients, it suffices to consider the two possibilities according to the singular isotropy group.

First, we assume that the singular orbit is a point.  The slice representation is the standard four-dimensional real representation of $SU(2)$, up to finite ineffective kernel.  Smoothness conditions at the origin yield the following~(cf.~\cite{VZ2022}).

\begin{lem}\label{lem:nut-smooth}
Under the hypothesis of Lemma \ref{lemma: asy warping functions}, if the singular orbit is a point, then there are smooth functions $\phi_i$ near $0$ such that
\begin{equation}\label{eq:nut-expansion}
  f_i(t)^2=\frac{t^2}{4}+t^4\phi_i(t^2),
  \qquad i=1,2,3.
\end{equation}
In particular,
\[
  \lim_{t\to0}\frac{f_i(t)}{f_j(t)}=1
  \qquad\text{for all }i,j\in\lbrace 1,2,3 \rbrace.
\]
\end{lem}

Since $q_{ij}=f_i^2/f_j^2$, Lemma \ref{lemma: asy warping functions} and Lemma \ref{lem:nut-smooth} show that the hypotheses of Lemma \ref{lem:two-cycle} are satisfied.  We obtain the following.

\begin{theo}\label{thm:fixed}
Let $(M,g,f)$ be a four dimensional complete $\kappa$-noncollapsed non-Ricci-flat steady gradient Ricci soliton. Suppose $(M,g)$ admits an $SU(2)$ cohomogeneity one action with the singular orbit a point, then the metric is rotationally symmetric and hence isometric to four-dimensional Bryant soliton up to scaling.
\end{theo}

\begin{proof}
By the discussion in the beginning of this section, we can choose a suitable Milnor frame $X_1,X_2,X_3$ with the normalization \eqref{eq:milnor} such that the metric $g$ on the principal part can be written as \eqref{su(2) invariant metric}. By Lemma \ref{lemma: asy warping functions} and Lemma \ref{lem:nut-smooth} and the definition of $q_{ij}$ in \eqref{eq:def-Bij-qij}, we can get
\begin{align*}
    \lim_{t\to 0}q_{ij}(t)=\lim_{t\to\infty}q_{ij}(t)=1,\qquad \forall\,i,j\in \lbrace 1,2,3 \rbrace.
\end{align*}
Then by Lemma \ref{lem:two-cycle}, for any $t>0$,
\begin{align*}
    q_{ij}(t)\leq 1,\qquad \forall\,i,j\in\lbrace 1,2,3 \rbrace.
\end{align*}
Note that $q_{ij}=q_{ji}^{-1}$. This implies that $q_{ij}\equiv 1$ for any $i,j\in \lbrace 1,2,3 \rbrace$ and hence $f_1\equiv f_2\equiv f_3$.  Thus $(M,g,f)$ is a four dimensional rotationally symmetric steady gradient Ricci soliton, which is isometric to the Bryant soliton~\cite{Bryant} up to scaling.
\end{proof}

Next, we consider the case that the singular orbit has codimension two. By Lemma \ref{lem:four two cases}, after passing to a finite cover, its group diagram is
\[
  H\subset K\cong S^1\subset SU(2).
\]
Since the singular orbit has codimension two, $K/H\cong S^1$, which implies that $H$ is a finite group. Let $X_1,X_2,X_3$ be the Milnor frame such that the metric on the principal part has the form \eqref{su(2) invariant metric} with the normalization \eqref{eq:milnor}. After a rotation, we assume that $X_1$ generates $\mathrm{Lie}(K)$. Note that the $X_1$ direction collapses at $t=0$, while the $X_2,X_3$ directions are tangent to the two-dimensional singular orbit.

Smoothness conditions around singular orbit give us the following expansions near $t=0$~(cf.~\cite{VZ2022} or \cite[Theorem 2]{B2025}).

\begin{lem}\label{lem:bolt-smooth}
Under the hypothesis of Lemma \ref{lemma: asy warping functions}, if the singular orbit has codimension 2 and the $X_1$ direction of the Milnor frame collapses on the singular orbit, then there are positive constants $b$ and $c$ such that
\begin{equation}\label{eq:bolt-expansion}
  f_1(t)^2=c^2t^2+O(t^4),
  \qquad
  f_2(t)^2=b^2+o(1),
  \qquad
  f_3(t)^2=b^2+o(1)
\end{equation}
as $t\to0$.  Consequently,
\begin{equation}\label{eq:bolt-ratio-limits}
  \frac{f_2}{f_3}\longrightarrow1,
  \qquad
  \frac{f_1}{f_2}\longrightarrow0,
  \qquad
  \frac{f_1}{f_3}\longrightarrow0
\end{equation}
as $t\to 0$.
\end{lem}

Note that by Lemma \ref{lemma: asy warping functions} and Lemma \ref{lem:bolt-smooth}, the hypothesis of Lemma \ref{lem:bolt-rigidity} is satisfied. Then we can get the following theorem.

\begin{theo}\label{thm:bolt}
Let $(M,g,f)$ be a four dimensional complete $\kappa$-noncollapsed non-Ricci-flat steady gradient Ricci soliton. Suppose $(M,g)$ admits an $SU(2)$ cohomogeneity one action with the singular orbit of codimension two. Then after passing to a finite cover, the isometric $SU(2)$ action extends to an isometric action locally isomorphic to $U(2)$.  If the original singular isotropy group $K$ is connected, no covering is needed.
\end{theo}

\begin{proof}
By the above discussion, after passing to a finite cover and choosing a suitable Milnor frame with the normalization \eqref{eq:milnor}, the metric $g$ on the principal part has the form \eqref{su(2) invariant metric} and the $X_1$  direction collapses at the singular orbit. Then by Lemma \ref{lem:bolt-smooth},
\[
  \lim_{t\to0}q_{23}=1,
  \qquad
    \lim_{t\to0} q_{12}=\lim_{t\to 0}q_{13}=0.
\]
By Lemma \ref{lemma: asy warping functions}, all $q_{ij}\to1$ as $t\to\infty$. Thus Lemma \ref{lem:bolt-rigidity} gives $f_2=f_3$. This implies that the orbit metrics are biaxial with symmetry axis generated by $X_1$. Now the theorem follows directly from Lemma \ref{lem:connected-cover}.
\end{proof}

Now the final step of proving Theorem \ref{thm: four} is to identify the $U(2)$-invariant steady soliton in Theorem \ref{thm:bolt} with the noncollapsed 4d Appleton soliton on $\mathcal O(-k)$ for some $k\geq 3$, up to a scaling and a finite covering. First, we show that the universal covering of $M$ in Theorem \ref{thm:bolt} is diffeomorphic to the total space of complex line bundle $\mathcal{O}(-k)\to \mathbb{CP}^1$ for some $k\geq 3$. 

\begin{lem}\label{lem:diffeo line bundle}
Under the hypothesis of Theorem \ref{thm:bolt}, there is a finite covering  $\widehat M$ of $M$ such that $\widehat M$ is diffeomorphic to $\mathcal O(-k)$ for some $k\geq 3$.
\end{lem}

\begin{proof}
By Lemma \ref{lem:four two cases}, there is a finite covering $(\widehat M,\widehat g,\widehat f)$ such that it is an $SU(2)$-invariant cohomogeneity one steady gradient Ricci soliton with singular isotropy group $K\cong S^1$ and principal isotropy group $H$ such that $K/H\cong S^1$. This implies that $H$ is finite. Since every finite subgroup of a circle is cyclic, we have
\[
  H\cong\Z_k
\]
for a unique $k\ge1$.  The slice representation
\[
  \rho:K\longrightarrow SO(2)
\]
is transitive on the unit normal circle and has kernel $H$.  After choosing
the orientation of the normal complex line, it is therefore the weight
$-k$ representation
\[
  z\cdot w=z^{-k}w,
  \qquad z\in S^1,\quad w\in\mathbb{C}.
\]
Consequently, we have
\begin{equation}\label{eq:associated-line}
  \widehat M
  =SU(2)\times_{S^1}\mathbb{C}_{-k}
  \cong\mathcal O(-k).
\end{equation}

Next, we need to show that it is impossible that $k=1,2$. For this, by Theorem \ref{thm:bolt}, we can assume $(\widehat M,\widehat g)$ is $U(2)$-invariant. To match the normalization in Section 2, we write $\sigma_i=\omega_i/2$ and set $a=2f_1$ and $b=2f_2=2f_3$.
Then $\theta=\sigma_1$ and $\pi^*\widehat g_B=\sigma_2^2+\sigma_3^2$, where $\operatorname{Ric}_{\widehat g_B}=4\widehat g_B$. Thus $\widehat g$ has the form \eqref{eq:ansatz} with $n=1$ and $p=2$. After scaling and adding a constant to the potential, we assume the the intial condition \eqref{eq:initial-data}. Now we consider the function
    \begin{align*}
        D=2(1-Q^2)-bQ^\pr.
    \end{align*}
    By a direct computation we can show that $D$ satisfies
    \begin{align}\label{eq:prime D}
        D^\pr=\frac{2Q}{b}D+(-bf^{\pr}+2b^{\pr}-2Q)Q^{\pr}.
    \end{align}
    We claim that
    \begin{claim}\label{claim:sign of D}
        For any $t>0$, $D^\pr(t)>0$ whenever $D(t)\geq 0$.
    \end{claim}
    \begin{proof}[Proof of Claim \ref{claim:sign of D}]
        We have the following formula of scalar curvature~(cf.~\cite[(9.1)]{Appleton})
        \begin{align*}
            R=-2\frac{a^{\pr\pr}}{a}-4\frac{b^{\pr\pr}}{b}-4\frac{a^\pr b^\pr}{ab}-2\Big(\frac{b^\pr}{b}\Big)^2-2\frac{a^2}{b^4}+\frac{8}{b^2}.
        \end{align*}
        By the third identity in \eqref{eq:soliton-identities}, we can get
        \begin{align*}
            -2\frac{a^{\pr\pr}}{a}-4\frac{b^{\pr\pr}}{b}-4\frac{a^\pr b^\pr}{ab}-2\Big(\frac{b^\pr}{b}\Big)^2-2\frac{a^2}{b^4}+\frac{8}{b^2}+(f^\pr)^2=C_0>0.
        \end{align*}
        By the soliton equation \eqref{eq:prelim-ode}, this implies that
        \begin{align*}
            \Big( -f^\pr+\frac{a^\pr}{a}+2\frac{b^\pr}{b} \Big)^2=C_0+\Big(\frac{a^\pr}{a}\Big)^2+2\Big( \frac{b^\pr}{b} \Big)^2+\frac{8-2Q^2}{b^2}.
        \end{align*}
        Multiplying this by $a^2$ then gives
        \begin{align*}
            (-af^\pr+a^\pr+2Qb^\pr)^2-(a^\pr+2Q^2)^2=C_0a^2+2Q^2(b^\pr-Q)^2+4Q^2D.
        \end{align*}
        If $t>0$ and $D(t)\geq 0$, then the right-hand side is positive. By Theorem \ref{thm:cohom-est}, $f^{\pr}<0$ and hence
        \begin{align*}
            -af^\pr+a^\pr+2Q b^\pr>a^\pr+2Q^2.
        \end{align*}
        Since $Q=a/b>0$ for all $t>0$, this implies that
        \begin{align*}
            -bf^\pr+2b^\pr-2Q>0.
        \end{align*}
        Recall $Q^{\prime}>0$ by Corollary \ref{cor:prelim-ratio}. By the expression of $D^\pr$ in \eqref{eq:prime D}, we can see that $D^\pr>0$ and the claim is proved.
    \end{proof}
    Now we can show that $D(t)<0$ for all $t>0$. Suppose there exists $t_0>0$ such that $D(t_0)\geq 0$. Then by Claim \ref{claim:sign of D}, $D(t)>0$ for any $t\geq t_0$ and $D$ has a positive lower bound for large $t$. However, this contradicts
    \begin{align*}
        D(t)<2(1-Q^2)\longrightarrow 0
    \end{align*}
    as $t\to+\infty$, where we have used $\lim_{t\to\infty} Q(t)=1$ by Lemma~\ref{lemma: asy warping functions}. This shows that $D(t)<0$ for any $t>0$.

     Finally, by \eqref{eq:initial-data} with $n=1$ and $p=2$,
    \[
        D(0)=2-b(0)Q'(0)=2-k.
    \] 
    Thus $k=1$ is excluded. Suppose $k=2$, then $D(0)=D^\pr(0)=0$. For a non-Ricci-flat complete steady soliton, we have $R+\Del f=0$. In our setting, this implies at $t=0$, 
    \begin{align*}
        0<C_0=R(0)=\lim_{t\to 0^+}\bigg(- f^{\pr\pr}-\Big( \frac{a^\pr}{a}+2\frac{b^\pr}{b} \Big)f^\pr \bigg)=-2f^{\pr\pr}(0).
    \end{align*}
    Therefore, by the Taylor expansion \eqref{eq:prelim-taylor} we can get
    \begin{align*}
        D^{\pr\pr}(0)=&(-f^{\pr\pr}(0)+2b^{\pr\pr}(0)-2Q^{\pr}(0))Q^\pr(0)\\
        =&-2f^{\pr\pr}(0)>0.
    \end{align*}
    Consequently, $D(t)>0$ for sufficiently small $t$, which is a contradiction.
\end{proof}

Finally, we classify all possibilities of quotients of four dimensional cohomogeneity one steady solitons if the singular orbit has codimension 2.

\begin{lem}\label{lem:possible quotients}
    Let $(M,g,f)$ be a non-Ricci-flat steady gradient Ricci soliton and $(\widehat M,\widehat g, \widehat f)$ be the noncollapsed 4d Appleton soliton on $\mathcal O(-k)$ for some $k\geq 3$. If there exists a map $\pi:(\widehat M,\widehat g)\to (M,g)$ such that $\pi$ is a finite Riemannian covering and $\widehat{f}$ is a lift of $f$. Then there are only two possibilities:
    \begin{enumerate}[label=\textup{(\roman*)}]
        \item If $k$ is odd, then $(M,g)$ is isometric to $(\widehat M,\widehat g)$;
        \item If $k$ is even, then $(M,g)$ is isometric to $(\widehat M,\widehat g)$ or the canonical free $\mathbb{Z}_2$ quotient of $(\widehat M,\widehat g)$.
    \end{enumerate}
\end{lem}

\begin{rem}
    The explicit form of the canonical free $\mathbb{Z}_2$ action on $(\widehat M,\widehat g)$ can be seen in the proof.
\end{rem}

\begin{proof}
    Since $\widehat{M}\cong \mathcal{O}(-k)$, which is simply connected, $\pi$ is the universal cover. Let
    \begin{align*}
        M=\widehat M/\Ga,
    \end{align*}
    where $\Ga$ is its finite deck group which acts isometrically and freely on $(\widehat M,\widehat g)$. In the following, we use hatted objects to denote the quantities in the covering manifold and unhatted objects to denote quantities in the base manifold.

    By the smoothness conditions, $f^\pr(0)=0$. By the soliton identity \eqref{eq:soliton-identities} we can get that $R(0)=C_0$. Then the singular orbit $\Sph^2$, as a submanifold in $(\widehat M,\widehat g)$ satisfies $\Sph^2\subset \lbrace R=C_0 \rbrace$. By Theorem \ref{thm:cohom-est}, we know that the converse is also true, hence
    \begin{align*}
        \Sph^2=\lbrace R=C_0 \rbrace.
    \end{align*}
    Since isometries preserve scalar curvature, this implies that they also preserve the singular orbit $\Sph^2$. The induced metric on $\Sph^2$ is a round metric. Hence, the restriction action of $\Ga$ on $\Sph^2$ gives us a homomorphism
    \begin{align*}
        \Ga \longrightarrow  \mathrm{Isom}(\Sph^2)\cong O(3).
    \end{align*}
    This action on $\Sph^2$ is free since $\Ga$ acts on $(\widehat M,\widehat g)$ freely. Thus, $\Sph^2/\Ga$ is also a quotient Riemannian manifold. By the classical result of classification of spherical space forms~(cf. for example, \cite{Wolfbook}), there are only two possibilities: $\Ga$ is trivial or $\Ga\cong \mathbb{Z}_2$. If $\Ga$ is trivial, then $(M,g)$ is isometric to $(\widehat M,\widehat g)$. It suffices to consider $\Ga\cong \mathbb{Z}_2$. 

     If $\Ga\cong\mathbb{Z}_2$, let $\tau$ be its generator. Then by the argument above, we have $\tau|_{\Sph^2}=-\Id_{\Sph^2}$. As in the proof of Lemma \ref{lem:diffeo line bundle}, we have the following identification of $\widehat M$ with an associated complex line bundle:
    \begin{align*}
        \widehat M=SU(2)\times _{\Sph^1}\mathbb{C}_{-k}\cong \mathcal{O}(-k),
    \end{align*}
    where $SU(2)$ is identified with the unit quaternions and $[st,z]=[s,t^{-k}z]$ for any $s\in SU(2)$ and $t\in S^1$. Let $j$ be the usual quaternion satisfying
    \begin{align}\label{eq:property j}
        j^2=-1,\qquad jtj^{-1}=t^{-1},\quad \forall\,t\in \Sph^1.
    \end{align}
    Let $P_\theta([s,z])=[s,e^{i\theta}z]$ be the isometric fiber rotations and $\Sigma\subset\widehat M$ denote the singular orbit, identified
with the zero section of $\mathcal O(-k)$, so that $\Sigma\cong S^2$.

\begin{claim}\label{claim:isometry rotation}
Every isometry of $(\widehat M,\widehat g)$ fixing $\Sigma$
pointwise is equal to $P_\theta$ for some
$\theta\in\mathbb R/2\pi\mathbb Z$.
\end{claim}

\begin{proof}
Let $\psi$ be such an isometry. Its differential induces an
orthogonal automorphism of the normal bundle
\[
    \nu\Sigma\cong\mathcal O(-k)_{\mathbb R}
\]
covering the identity on $\Sigma$.
Since $\Sigma$ is connected, this automorphism either preserves
the normal orientation everywhere or reverses it everywhere.
The latter possibility would imply
\[
    e(\nu\Sigma)=-e(\nu\Sigma).
\]
This is impossible because $H^2(\Sigma;\mathbb Z)\cong\mathbb Z$
is torsion-free and
\[
    \langle e(\nu\Sigma),[\Sigma]\rangle=-k\ne0.
\]
Thus $d\psi$ restricts to a rotation on each normal plane.

Fix $p\in\Sigma$. Since the fiber rotations $P_\theta$ induce
all rotations of $\nu_p\Sigma$, we may choose $\theta$ such that
\[
    dP_{\theta,p}|_{\nu_p\Sigma}
    =d\psi_p|_{\nu_p\Sigma}.
\]
Both maps fix $\Sigma$ pointwise, hence their differentials also agree on $T_p\Sigma$. Consequently,
\[
    P_\theta(p)=\psi(p)=p,
    \qquad dP_{\theta,p}=d\psi_p.
\]
Since an isometry of a connected Riemannian manifold is determined
by its value and differential at one point, we get that $\psi=P_\theta$.
\end{proof}
    We define $J_k:\widehat M\to \widehat M$ by
    \begin{align}\label{eq:defi Jk}
        J_k([s,z])=[sj,\bar{z}].
    \end{align}
     Using the property of $j$ we can show that $J_k$ is a well-defined smooth map. Let $i=\sqrt{-1}=e^{\pi i/2}\in S^1$. Let $T=S^1$. Then under the identification
    \begin{align*}
        SU(2)/\Sph^1\to \Sph^2,\qquad sT\mapsto sis^{-1},
    \end{align*}
    $J_k$ acts on the base $SU(2)/\Sph^1=\Sph^2$ as
    \begin{align*}
        sis^{-1}\mapsto sjij^{-1}s^{-1}=-sis^{-1}.
    \end{align*}
    Thus $J_k$ covers the antipodal map and has no fixed points. On the regular part, $J_k$ is induced by right multiplication by $j$ on $SU(2)/\mathbb{Z}_k$. It reverses the Hopf direction $\si$ and preserves the horizontal metric $g_{FS}$ as the following.
    \begin{align*}
        J_k^*\sigma=-\sigma,\qquad J_k^*(\pi^*g_{FS})=\pi^*g_{FS}.
    \end{align*}
    It also preserves the radial $t$. Consequently,
    \begin{align*}
        J_k^*\widehat g=\widehat g,\qquad J_k^*\widehat f=\widehat f.
    \end{align*}
    This shows that $J_k$ is an isometry acting on the singular orbit $\Sph^2$ as an antipodal map. Hence, $\tau J_k^{-1}$ fixes $\Sph^2$ pointwise. By Claim \ref{claim:isometry rotation}, there exists a constant $\theta$ such that
    \begin{align*}
        \tau=P_\theta J_k.
    \end{align*}
    By \eqref{eq:defi Jk}, $J_kP_\theta J_k^{-1}=P_{-\theta}$. It follows that
    \begin{align*}
        \tau^2=(P_\theta J_k)^2=P_\theta P_{-\theta}J_k^2=J_k^2.
    \end{align*}
    Since $\tau^2=\Id$, $J_k^2=\Id$. By a direct computation,
    \begin{align*}
        [s,z]=J_k^2([s,z])=[sj^2,z]=[-s,z]=[s,(-1)^kz]
    \end{align*}
    for any $[s,z]\in SU(2)\times_{\Sph^1}\mathbb{C}_{-k}$. Hence, $k$ must be an even number. Finally, note that
    \begin{align*}
        P_{\theta/2}J_kP_{-\theta/2}=P_\theta J_k=\tau.
    \end{align*}
    Thus every possible deck involution is conjugate to $J_k$ by an isometry. In particular, all the resulting quotients of $(\widehat M,\widehat g)$ are isometric.
\end{proof}

Now we can complete the proof of Theorem \ref{thm: four}.

\begin{proof}[Proof of Theorem \ref{thm: four}]
By Lemma \ref{lem:average}, we can choose a $G$-invariant potential.  By Corollary \ref{cor:no-exceptional}, the orbit space $M/G$ is a closed half line and there is only one singular orbit of codimension at least 2. 

Using the asymptotic geometry, Corollary \ref{coro: principal orbit} shows that the principal orbit is a spherical space form. Then Lemma \ref{lem:transitive} tells us that after passing to a finite cover and a finite ineffective kernel, $G$ is isomorphic to $SU(2)$, $U(2)$ or $SO(4)$. In each of the three cases, the acting group contains a connected subgroup locally isomorphic to $SU(2)$ that acts transitively on the principal orbits. Restricting to this subgroup and lifting through its covering group $SU(2)$ by using Lemma \ref{lemma: lift action}, we obtain an almost effective isometric $SU(2)$-action with the same orbit decomposition. Hence it suffices to prove the theorem for an $SU(2)$-action. 

By Lemma \ref{lem:four two cases}, after passing to a finite covering, the singular isotropy group is either all of $SU(2)$, giving a fixed point, or is a circle on the connected-isotropy cover, giving a codimension-two orbit.  The fixed-point case is handled by Theorem \ref{thm:fixed} and the codimension-two case is handled by Theorem \ref{thm:bolt}.  Hence a finite Riemannian cover has an isometric action locally isomorphic to $U(2)$ and the original metric is locally $U(2)$-invariant. Moreover, if the singular orbit has codimension 2, Corollary \ref{cor:unique four noncollapsed soliton} and Lemma \ref{lem:diffeo line bundle} show that the universal cover of $(M,g)$ must be the noncollapsed 4d Appleton soliton on $\mathcal{O}(-k)$ for some $k\geq 3$. Lemma \ref{lem:possible quotients} classifies all possible quotients, hence the proof of Theorem \ref{thm: four} is done.
\end{proof}

\section{The five-dimensional topological trichotomy}\label{sec:five-topology}

Let $(M^5,g,f)$ be a five dimensional steady gradient Ricci soliton satisfying the hypotheses of Theorem \ref{thm:five}.  By Proposition
\ref{prop:kappa-entropy}, $\pi_1(M)$ is finite. By Lemma \ref{lemma: lift action}, after passing to a finite universal
Riemannian cover and lifting the action, we assume that $M$ is simply connected. For simplicity, we still use $G$ to denote the lifted Lie group. By Corollary \ref{cor:no-exceptional}, the orbit space remains $[0,\infty)$.

For the cohomogeneity one $G$-action on $M$, let
\[
  H\subset K\subset G
\]
be the group diagram. It follows that $P=G/H$ is the principal and $Q=G/K$ is the singular orbit. Corollary \ref{cor:no-exceptional} also tells us that
\begin{align*}
    K/H=\Sph^{\ell},\qquad1\leq \ell\leq 4.
\end{align*}
Moreover, we can reconstruct $M$ by
\[
  M\cong G\times_K\mathbb{R}^{\ell+1}.
\]
In particular, this shows that $M$ deformation retracts onto $Q$. Since we assume that $M$ is simply connected, $Q$ is also simply connected. By $K/H\cong \mathbb{S}^l$, we can get the codimension of $Q$ is $1+\ell$, which implies that
\[
  \dim Q=4-\ell,\qquad 1\le\ell\le4.
\]

We have the following topological trichotomy.

\begin{prop}\label{prop:five-trichotomy}
Let $(M^5,g,f)$ be a simply connected complete non-Ricci-flat cohomogeneity one steady gradient Ricci soliton invariant under a compact connected Lie group $G$-action. Let $K,H,P,Q,\ell$ be defined as above. Then exactly one of the following occurs:
\[
\begin{array}{c|c|c|c}
\ell &Q&P&M\\ \hline
4&\{\mathrm{pt}\}&\Sph^4&\mathbb{R}^5\\
2&\Sph^2&\Sph^2\times\Sph^2&\mathbb{R}^3\times\Sph^2\\
1&\Sph^3&\Sph^3\times\Sph^1&\Sph^3\times\mathbb{R}^2.
\end{array}
\]
In particular, the case $\ell=3$ is impossible.
\end{prop}

\begin{proof}
By the above discussion, $1\leq \ell\leq 4$. If $\ell=4$, then $Q$ is a point and $G=K$,  $P=K/H\cong\Sph^4$. Thus, $M\cong G\times_K\mathbb R^{\ell+1}=K\times _K\mathbb R^{5}=\mathbb R^5$.

If $\ell=3$, then $Q$ would be a compact connected simply connected
one-dimensional manifold, which does not exist.

Suppose $\ell=2$.  Then $Q$ is a compact simply connected homogeneous
surface, hence $Q\cong\Sph^2$, and
\[
  \Sph^2=K/H\longrightarrow P=G/H\longrightarrow G/K=\Sph^2
\]
is a fiber bundle.  Thus $P$ is simply connected and
\[
  \chi(P)=\chi(\Sph^2)^2=4.
\]
By the classification result in \cite{Gorbatsevich1977,Gorbatsevich2012}(cf. Theorem \ref{thm:homogeneous-four}), a simply
connected compact homogeneous four-manifold is diffeomorphic to
$\Sph^4$, $\mathbb{C} P^2$, or $\Sph^2\times\Sph^2$.  Their Euler characteristics
are $2$, $3$, and $4$, respectively.  Hence
\[
  P\cong\Sph^2\times\Sph^2.
\]
It remains to determine the normal bundle of $Q$.  Let
$\nu_Q\to\Sph^2$ be its oriented rank-three normal bundle.  The principal orbit is its unit sphere bundle as the following.
\[
  S(\nu_Q)=G\times_K(K/H)\cong G/H=P.
\]
Any oriented rank-three bundle over $\Sph^2$ admits a nonzero section.  If
$\pi:S(\nu_Q)\to\Sph^2$ is the sphere bundle and $V$ its vertical tangent
bundle, then
\[
  V\oplus\mathbb{R}\cong\pi^*\nu_Q,\qquad
  TS(\nu_Q)\cong\pi^*T\Sph^2\oplus V.
\]
Therefore,
\[
  w_2\bigl(TS(\nu_Q)\bigr)
  =\pi^*\bigl(w_2(T\Sph^2)+w_2(\nu_Q)\bigr)
  =\pi^*w_2(\nu_Q).
\]
The total space $S(\nu_Q)\cong\Sph^2\times\Sph^2$ is spin.  Since the sphere
bundle has a section, $\pi^*$ is injective.  Thus $w_2(\nu_Q)=0$.
Oriented rank-three bundles over $\Sph^2$ are classified by
$\pi_1(SO(3))\cong\Z_2$, with the class detected by $w_2$.
Thus $\nu_Q$ is trivial and consequently
\[
  M\cong E(\nu_Q)\cong\Sph^2\times\mathbb{R}^3.
\]

Finally, suppose $\ell=1$.  Then $Q$ is a compact simply connected
homogeneous three-manifold, which implies that $Q\cong\Sph^3$ by Poincare theorem. Now we can get the fibration
\[
  \Sph^1\longrightarrow P\longrightarrow\Sph^3
\]
is trivial since its Euler class lies in $H^2(\Sph^3;\Z)=0$.  This gives us that
$P\cong\Sph^3\times\Sph^1$.  The oriented normal two-plane bundle of $Q$ is
also classified by $H^2(\Sph^3;\Z)$ and is trivial.  Therefore
$M\cong\Sph^3\times\mathbb{R}^2$.
\end{proof}

By Proposition \ref{prop:five-trichotomy}, it suffices to consider three cases. The cases $\ell=4$ or 2, we will see, are relatively easy to handle. 

\subsection{The fixed point branch \texorpdfstring{$\ell=4$}{}}

First, if the singular orbit is a point, i.e., $\ell=4$, the corresponding steady soliton is the five dimensional Bryant soliton.

\begin{prop}\label{prop: 5d Bryant}
    If $\ell=4$ in Proposition \ref{prop:five-trichotomy}, then $M$ is isometric to the five dimensional Bryant soliton up to scaling.
\end{prop}

\begin{proof}
    If $\ell=4$, then the singular orbit is a point.  Hence, the principal orbit is
$\Sph^4$.  By the classification of compact connected transitive groups on
spheres \cite{MS1943}, the isotropy representation of the lifted
effective action on $\Sph^4$ is irreducible.  This implies that every invariant orbit
metric is round up to scale and the soliton metric $g$ has the following form.
\[
  g=\md t^2+a(t)^2g_{\Sph^4(1)}.
\]
By smoothness conditions on the singular orbit,
\[
  a(0)=0,\qquad  a^{\prime}(0)=1,\qquad  f^{\prime}(0)=0.
\]
This shows that the soliton is rotationally symmetric. Hence, it is isometric to the five-dimensional Bryant soliton up to scaling~\cite{Bryant}.
\end{proof}

\subsection{The \texorpdfstring{$\Sph^2$}--bolt branch \texorpdfstring{$\ell=2$}{}}
Suppose the second row of Proposition \ref{prop:five-trichotomy} occurs.  The proof of Theorem
 \ref{thm:homogeneous-four} shows that
\[
  G=SU(2)\times SU(2),\qquad H=\Sph^1\times\Sph^1
\]
after removing ineffective kernels and passing to finite covers. The connected subgroup $K$ has
codimension two in $G$ and $K/H\cong\Sph^2$.  At the level of roots it is obtained by adjoining
the root spaces of exactly one $A_1$ factor.  After interchanging the two
factors, we have the following group diagram of this cohomogeneity one action.
\[
  \Sph^1\times\Sph^1
  \ \subset\
  SU(2)\times\Sph^1
  \ \subset\
  SU(2)\times SU(2).
\]
Thus the first $\Sph^2$ factor collapses at the singular orbit and the second
one is the bolt.

By the proof of Theorem \ref{thm:homogeneous-four}, we can see the isotropy representation on the principal orbit in this case splits into two
inequivalent real two-dimensional modules,
\[
  \mathfrak m=\mathfrak m_1\oplus\mathfrak m_2.
\]
The first circle of $H$ rotates $\mathfrak m_1$ and acts trivially on
$\mathfrak m_2$, while the second circle has the reverse action.  Schur's
lemma excludes mixed terms and shows that every invariant metric has the
form
\begin{equation}\label{eq:idw-metric}
  g=\md t^2+a(t)^2g_{\Sph^2(1)}+b(t)^2g_{\Sph^2(1)}.
\end{equation}
The bolt conditions are
\begin{equation}\label{eq:idw-bolt}
  a(0)=0,\quad  a^\pr(0)=1,\qquad
  b(0)=b_0>0,\quad  b^\pr(0)=0,\qquad
   f^\pr(0)=0.
\end{equation}

Let $\xi$ be defined as in \eqref{eq:def-LHxi}. Then with respect to the metric \eqref{eq:idw-metric}, $\xi$ can be represented as
\begin{align}\label{eq:xi-explicit-form}
    \xi=2\frac{a^{\pr}}{a}+2\frac{b^\pr}{b}-f^{\pr}.
\end{align}
As in Dancer-Wang\cite{DancerWangExamples}, we can make the following change of variables.
\begin{equation}\label{eq:idw-XY}
  X_1=\frac{\sqrt2 a^\pr}{\xi a},\qquad
  X_2=\frac{\sqrt2   b^\pr}{\xi b},\qquad
  Y_1=\frac{\sqrt2}{a\xi},\qquad
  Y_2=\frac{\sqrt2}{b\xi}.
\end{equation}
Moreover, we define a new variable $s$ by
\begin{equation}\label{eq:idw-s}
  \frac{\md}{\md s}=\frac1{\xi}\frac{\md}{\md t}.
\end{equation}
Then the soliton equations \eqref{eq:cohom-radial}-\eqref{eq:cohom-tangent} descend to the following system of $X_i,Y_i$~(cf.~\cite[Section~2]{DancerWangExamples}).
\begin{align}
  X_i'&=X_i(X_1^2+X_2^2-1)+\frac{Y_i^2}{\sqrt2},\label{eq:idw-X}\\
  Y_i'&=Y_i\left(X_1^2+X_2^2-\frac{X_i}{\sqrt2}\right),
  \qquad i=1,2,\label{eq:idw-Y}
\end{align}
where prime now denotes derivatives with respect to $s$. Let
\begin{align}\label{eq:def-L}
    \mathcal{L}=X_1^2+X_2^2+Y_1^2+Y_2^2-1.
\end{align}
If we assume $C_0=1$ in the soliton identity \eqref{eq:soliton-identities}, then this identity implies~(cf.~\cite[equation~(2.4)]{DancerWangExamples})
\begin{equation}\label{eq:idw-L}
  \mathcal L=X_1^2+X_2^2+Y_1^2+Y_2^2-1
  =-\frac{1}{\xi^2}<0
\end{equation}
Using the system \eqref{eq:idw-X}--\eqref{eq:idw-Y}, we can compute
\begin{equation}\label{eq:idw-Lprime}
  \mathcal L'=2(X_1^2+X_2^2)\mathcal L.
\end{equation}
These are exactly equations (2.7)--(2.9) of
Dancer--Wang for $r=2$, $d_1=d_2=2$ and Einstein constants
$\lambda_1=\lambda_2=1$ \cite{DancerWangExamples}. In this case, Dancer--Wang proved the following.

\begin{theo}[{\cite[Section~3 and Theorem~4.17]{DancerWangExamples}}]\label{thm:DW-construction}
Let $(X_1,X_2,Y_1,Y_2)$ satisfy the system \eqref{eq:idw-X}--\eqref{eq:idw-Y} and
\[
  P_0=\left(\frac1{\sqrt2},0,\frac1{\sqrt2},0\right),
  \qquad
  \mathcal L=X_1^2+X_2^2+Y_1^2+Y_2^2-1.
\]
Then the trajectories of this system satisfying
\[
  \lim_{s\to-\infty}(X_1,X_2,Y_1,Y_2)=P_0,
  \qquad
  Y_1,Y_2>0,\qquad \mathcal L<0
\]
form a one-parameter family modulo translation in $s$. For every such trajectory, if we define 
\begin{equation}\label{eq:reconstruction-abf}
    \begin{aligned}
       t(s)=\int_{-\infty}^{s}\sqrt{-\mathcal{L}(\sigma)}\,d\sigma,\qquad
  a(s)=\frac{\sqrt{-2\mathcal{L}(s)}}{Y_1(s)},\\
  b(s)=\frac{\sqrt{-2\mathcal{L}(s)}}{Y_2(s)},\qquad
  \frac{df}{ds}=\sqrt2(X_1+X_2)-1, 
    \end{aligned}
\end{equation}
then they define a cohomogeneity one metric on the principal part of $\mathbb R^3\times \Sph^2$ under $SU(2)\times SU(2)$ action by 
\[
  g=dt^2+a(t)^2g_{S^2(1)}+b(t)^2g_{S^2(1)}.
\]
Moreover, this metric extends smoothly across an $\Sph^2$ singular orbit at $t=0$ and is a complete steady soliton on $\mathbb R^3\times S^2$ with the normalization $R+|\nabla f|^2=1$. These solitons constitute the normalized Ivey--Dancer--Wang family.
\end{theo}

\begin{proof}
    Specializing to $r=2$ and $d_1=d_2=2$, the unstable manifold discussion immediately before Remark~3.5 in \cite{DancerWangExamples} gives the family of trajectories. Completeness follows from their Lemma~3.6, while smooth extension across the singular orbit is established in the proof of their Theorem~4.17.
\end{proof}

To obtain the classification in this case, it remains to show that an arbitrary soliton in the $\Sph^2$-bolt branch determines one of these trajectories.

\begin{prop}\label{prop:S2-bolt-classification}
Let $(M^5,g,f)$ be a complete non-Ricci-flat simply connected steady gradient
Ricci soliton invariant under a cohomogeneity one action with a singular orbit of codimension 3. Suppose $f$ is also invariant under this action. Then the soliton belongs to the Ivey--Dancer--Wang family $\IDW$ up to homothety and addition of a constant to $f$.
\end{prop}

\begin{proof}
Since $(M^5,g,f)$ is a complete non-Ricci-flat simply connected cohomogeneity one steady gradient Ricci soliton with a singular orbit of codimension 3, we can use Proposition \ref{prop:five-trichotomy} to get that $M\cong \mathbb R^3 \times \Sph^2$ with principal orbit $\Sph^2\times \Sph^2$ and singular orbit $\Sph^2$. Then by the above discussion, we can see the metric $g$ on the principal part has the form
\begin{align*}
    g=dt^2+a(t)^2g_{\mathbb S^2}+b(t)^2g_{\mathbb S^2}.
\end{align*}
By the smoothness conditions at the singular orbit, as $t\to 0$,
\[
  a(t)=t+O(t^3),\qquad
  b(t)=b_0+O(t^2),\qquad
  f'(t)=O(t),
  \qquad b_0>0.
\]
Let $\xi$ be defined as in \eqref{eq:xi-explicit-form}. Then by Theorem \ref{thm:cohom-est}, $\xi(t)>0$ for any $t>0$. Hence, the change of variables \eqref{eq:idw-XY}--\eqref{eq:idw-s} of $X_i,Y_i,s$ is well defined near the singular orbit and the steady soliton equation descends to the system \eqref{eq:idw-X}--\eqref{eq:idw-Y}. Moreover, by a direct computation, as $t\to0$ we have
\begin{align}\label{eq:xi-expansion}
    \xi=2\frac{a'}a+2\frac{b'}b-f'
      =\frac2t+O(t).
\end{align}
Then by the definition of $s$ and \eqref{eq:xi-expansion}, we can see as $t\to 0$,
\[
  s=2\log t+c+O(t^2).
\]
Thus $t\downarrow0$ corresponds to $s\to-\infty$. By the definitions of $X_i,Y_i$, we can also compute the following expansions of them near $t=0$.
\[
  X_1=\frac1{\sqrt2}+O(t^2),\qquad
  X_2=O(t^2),\qquad
  Y_1=\frac1{\sqrt2}+O(t^2),\qquad
  Y_2=\frac{t}{\sqrt2\,b_0}+O(t^3).
\]
This implies that the phase trajectory converges to $P_0$ as $s\to-\infty$ and $Y_1(t),Y_2(t)>0$ for sufficiently small $t$. By the system \eqref{eq:idw-Y} satisfied by $Y_1,Y_2$, the positivity of $Y_1,Y_2$ is preserved along the trajectory, hence
\[
  Y_1(t),\,Y_2(t)>0,\qquad \forall\,t>0. 
\]
Let $\mathcal L$ be defined as in \eqref{eq:def-L}. We need to show $\mathcal{L}<0$. For this, after a homothety, we assume $C_0=1$ in the soliton identity \eqref{eq:soliton-identities}. Then by \eqref{eq:idw-L}, this soliton identity implies that
\begin{align*}
    \mathcal{L}=-\frac{1}{\xi^2}<0.
\end{align*}
Now using Theorem~\ref{thm:DW-construction}, we can get that this trajectory reconstructs a member of the normalized Ivey--Dancer--Wang family.

It remains to identify this reconstruction with the given soliton. Indeed, since $\mathcal{L}=-\xi^{-2}$ and $X_i,Y_i$ are defined by \eqref{eq:idw-XY}, the reconstructed formula \eqref{eq:reconstruction-abf} gives us the same $a,b,f^{\pr}$. This implies that the two solitons are identical up to addition of a constant to $f$.
\end{proof}

As in Proposition \ref{prop:S2-bolt-classification}, we use the notation $\IDW$ for this full admissible phase-space family.  Its
members have asymptotically paraboloidal geometry
\begin{equation}\label{eq:idw-asymptotic}
  g=\md t^2+t\,h_\infty+o(t),
\end{equation}
where $h_\infty$ is homothetic to the product Einstein metric on
$\Sph^2\times\Sph^2$ \cite[Theorem 3.10]{DancerWangExamples}. From the asymptotically paraboloidal property, we can get steady solitons in $\IDW$ are $\kappa$-noncollapsed.

\begin{prop}\label{prop:idw-noncollapse}
Every member of $\IDW$ is $\kappa$-noncollapsed on all scales for some
$\kappa>0$ depending on the member.
\end{prop}

\section{Proof of Theorem \ref{thm:five}}\label{sec:s3r2}

To complete the proof of Theorem \ref{thm:five}, it suffices to exclude the $\mathbb{S}^3\times \mathbb{R}^2$ branch, i.e., the third row of Proposition \ref{prop:five-trichotomy}. Indeed, we will show that the steady solitons in this branch are collapsed.

We now assume that the third row of Proposition \ref{prop:five-trichotomy} occurs. In this case, $M\cong \mathbb{S}^3\times \mathbb R^2$ with the singular orbit $\Sph^3$ and the principal orbit $\Sph^3\times\Sph^1$. 

In this case, the isotropy representation may contain equivalent irreducible summands, hence invariant orbit metrics can have mixed terms between these summands. To tackle this difficulty, we need the following lemma.

\begin{lem}\label{lem:orbit-lie-group}
If $\ell=1$ in Proposition \ref{prop:five-trichotomy}, then after passing to a finite cover of the principal orbit, there is a locally
free~(i.e., the isotropy group is discrete) transitive isometric action on the principal orbit whose Lie algebra is
\[
  \mathfrak l=\mathfrak{su}(2)\oplus\mathbb{R} Z,\qquad [Z,\mathfrak{su}(2)]=0.
\]
Equivalently, the lifted orbit metric is an arbitrary left-invariant metric
on $SU(2)\times\mathbb{R}$, invariant under a fixed cocompact cyclic deck group.
\end{lem}

\begin{proof}
If $\ell=1$ in Proposition \ref{prop:five-trichotomy}, then we have group diagram $H\subset K\subset G$ with $G/K\cong\Sph^3$ and $K/H\cong\Sph^1$. Let $\mathfrak g,\mathfrak k,\mathfrak h$ be the Lie algebras of $G,K,H$ respectively.
Let $\mathfrak n$ be the Lie algebra of the kernel of the
$G$-action on $G/K$. By the classification of connected
compact transitive groups on $\Sph^3$ \cite{MS1943}, we obtain
a simple ideal
$\overline{\mathfrak s}\cong\mathfrak{su}(2)$ in
$\mathfrak g/\mathfrak n$ whose connected subgroup acts
transitively on $G/K$.
Since $\mathfrak g$ is compact, $\mathfrak n$ admits an ideal
complement. Thus $\overline{\mathfrak s}$ lifts to a simple
ideal $\mathfrak s\subset\mathfrak g$ whose connected subgroup
acts transitively on $G/K$ with finite isotropy. Consequently,
\[
\mathfrak g=\mathfrak s\oplus\mathfrak k
\qquad\text{as vector spaces}.
\]
Let $Z_0\in\mathfrak k$ whose image spans
$\mathfrak k/\mathfrak h$.  Since $\mathfrak s$ is an ideal,
$\ad_{Z_0}|_{\mathfrak s}$ is a derivation of $\mathfrak{su}(2)$. On a compact simple Lie algebra, every derivation is inner, which implies that there is $Y\in\mathfrak s$ with
\[
  \ad_{Z_0}|_{\mathfrak s}=\ad_Y.
\]
Let $Z=Z_0-Y$.  Then $[Z,\mathfrak s]=0$. If we define $\mathfrak l$ by
\[
  \mathfrak l:=\mathfrak s\oplus\mathbb{R} Z,
\]
then it satisfies $\mathfrak l+\mathfrak h=\mathfrak g$ and
$\mathfrak l\cap\mathfrak h=0$.  Thus its simply connected group
$\widetilde L=SU(2)\times\mathbb{R}$ acts locally freely and transitively on
$G/H$.

We now identify the global quotient underlying the orbit metrics.
Fix $o=eH\in G/H$ and let
$\pi:\widetilde L=SU(2)\times\mathbb R\to G/H$ be the orbit map.
Since $\mathfrak l\cap\mathfrak h=\{0\}$ and
$\mathfrak l+\mathfrak h=\mathfrak g$, this map is a local
diffeomorphism. Its equivariance implies that $\pi^*g_t$ is
left-invariant and hence complete. Therefore $\pi$ is a covering map. Moreover, since $\widetilde L$ is simply connected, it is the universal covering.

Let $\Gamma$ be the stabilizer of $o$, we consequently have
$G/H\cong\widetilde L/\Gamma$ with deck transformations given
by right translations. The compactness of $G/H$ makes $\Gamma$
cocompact, while $G/H\cong S^3\times S^1$ gives
$\Gamma\cong\pi_1(G/H)\cong\mathbb Z$.
Finally, since $\pi$ and $\Gamma$ depend only on the fixed group action, the same cyclic deck group preserves every lifted metric
$\pi^*g_t$.
\end{proof}

Let
\[
  \mathfrak s=[\mathfrak l,\mathfrak l]=\mathfrak{su}(2).
\]
Fix the left-invariant one-form $\theta\in\mathfrak l^*$ satisfying
\begin{align}\label{eq:def-theta}
    \theta|_{\mathfrak s}=0,\qquad \theta(Z)=1.
\end{align}
Since $\theta$ annihilates the commutator algebra,
\begin{equation}\label{eq:theta-closed}
  \md\theta(X,Y)=-\theta([X,Y])=0.
\end{equation}
For each lifted orbit metric $g_t$, we define
\begin{equation}\label{eq:ell-def}
  \ell(t)=|\theta|^{-1}_{g_t^{-1}},\qquad
  W(t)=\ell(t)\theta^{\sharp_{g_t}}.
\end{equation}
Then
\[
  |W|_{g_t}=1,\qquad W\perp_{g_t}\mathfrak s.
\]
Let 
\begin{align}\label{def-c-fiveclass}
    c(t)=\frac{\dot\ell}{\ell},
\end{align}
then we have the following properties of $c(t)$.

\begin{lem}\label{lem:ell-shape}
The function $c(t)$ defined by \eqref{def-c-fiveclass} satisfies
\begin{equation}\label{eq:c-L}
  c=\langle LW,W\rangle,\qquad
  \dot W=-2LW+cW,
\end{equation}
where $L$ is the shape operator defined in \eqref{eq:def-LHxi}.
\end{lem}

\begin{proof}
By the definition of $L$ in \eqref{eq:def-LHxi}, $\dot g_t=2g_tL$. Then we can get
\[
  \frac{\md}{\md t} |\theta|^2_{g_t^{-1}}=-2\langle L\theta^\sharp,\theta^\sharp\rangle.
\]
Note that $\theta^\sharp=\ell^{-1}W$, this gives the first
identity of \eqref{eq:c-L}. Differentiating the identity
$g_t(\theta^\sharp,\cdot)=\theta$ gives
$\frac{\md}{\md t}\theta^\sharp=-2L\theta^\sharp$, from which the second identity follows.
\end{proof}

We need the following formula of Ricci curvature and scalar curvature about a unimodular Lie algebra and its codimension one derived ideal.

\begin{lem}\label{lem:codim-one-ricci}
Let $(\mathfrak l,\langle\cdot,\cdot\rangle)$ be a unimodular four dimensional metric Lie
algebra whose derived ideal $\mathfrak s=[\mathfrak l,\mathfrak l]$ has
codimension one.  If $X$ is a unit vector orthogonal to $\mathfrak s$, then
\begin{align}
  \Ric^{\mathfrak l}(X,X)&=-|\sym (\ad_X|_{\mathfrak s})|^2\le0,
  \label{eq:ric-W}\\
  \Scal(\mathfrak l)&=\Scal(\mathfrak s)-|\sym (\ad_X|_{\mathfrak s})|^2,
  \label{eq:scal-split}
\end{align}
where $\sym(\ad_X|_{\mathfrak s})$ denotes the symmetric part of the linear map $\ad_{X}|_\mathfrak s$ and $\mathfrak s$ carries the restricted metric.
\end{lem}

\begin{proof}
Let $e_1,e_2,e_3$ be an orthonormal basis of $\mathfrak s$.  The Ricci
formula for a unimodular metric Lie algebra is
\[
  \Ric(X,X)
  =-\frac12\sum_\alpha|[X,e_\alpha]|^2
   -\frac12B(X,X)
   +\frac14\sum_{\alpha,\beta}
     \langle[e_\alpha,e_\beta],X\rangle^2,
\]
where $B$ is the Killing form~(cf.~\cite{Milnor,B2007}). Note that the last term vanishes since all brackets lie in $\mathfrak s$. Since $\ad_X(e_\alpha)=[X,e_\al]$,
\begin{align*}
    \sum_\alpha|[X,e_\alpha]|^2=|\ad_X|_{\mathfrak s}|^2.
\end{align*}
By the definition of Killing form, $B(X,X)=\tr((\ad_X|_{\mathfrak s})^2)$.  Writing $\ad_{X}|_{\mathfrak s}=S+A$ with $S$ self-adjoint and $A$
skew-adjoint, then we have
\[
  |\ad_X|_{\mathfrak s}|^2=|S|^2+|A|^2,\qquad
  \tr((\ad_X|_{\mathfrak s})^2)=|S|^2-|A|^2.
\]
This proves \eqref{eq:ric-W}.

The scalar-curvature formula for a unimodular metric Lie algebra is
\[
  \Scal=-\frac14\sum_{i,j}|[e_i,e_j]|^2
        -\frac12\sum_iB(e_i,e_i).
\]
Comparing the formula for the orthonormal basis
$W,e_1,e_2,e_3$ with the formula on $\mathfrak s$, the extra terms are
\[
  -\frac12|\ad_W|_{\mathfrak s}|^2-\frac12\tr((\ad_W|_{\mathfrak s})^2)=-|\sym (\ad_W|_{\mathfrak s})|^2,
\]
which gives \eqref{eq:scal-split}.
\end{proof}

The following uniform boundedness of the abelianized radius $\ell(t)$ plays an important role in the proof of the collapsing property.

\begin{prop}\label{prop:ell-bounded}
Let the function $\ell$ be defined by \eqref{eq:ell-def}. Then $\ell(t)$ is uniformly bounded above on $[0,\infty)$.
\end{prop}

\begin{proof}
Differentiate $c=\langle LW,W\rangle$(note that both the metric and
$W$ depend on $t$) and using \eqref{eq:c-L} we can get
\[
\begin{aligned}
  \dot c
  &=2\langle L^2W,W\rangle+\langle\dot LW,W\rangle
    +2\langle L\dot W,W\rangle\\
  &=\langle\dot LW,W\rangle-2|LW|^2+2c^2.
\end{aligned}
\]
The tangential soliton equation \eqref{eq:cohom-tangent} and Lemma
\ref{lem:codim-one-ricci} yield
\begin{equation}\label{eq:c-damping}
  \dot c+\xi c
  =-|\sym(\ad_W|_{\mathfrak s})|^2
   -2|LW-cW|^2
  \le0.
\end{equation}
After a normalization, we assume $C_0=1$. Then by Theorem \ref{thm:cohom-est}, $\xi(t)\to 1>0$.  Hence there are
$T,\delta>0$ such that $\xi\ge\delta$ on $[T,\infty)$. The integrating
factor in \eqref{eq:c-damping} then gives
\[
  c_+(t)\le c_+(T)e^{-\delta(t-T)},
\]
where $c_+(t):= \max\lbrace c(t),0 \rbrace $. Therefore,
\[
  \log\frac{\ell(t)}{\ell(T)}
  =\int_T^tc(s)\,\md s
  \le\int_T^\infty c_+(s)\,\md s<\infty.
\]
It remains only to control the singular end.  In the notation of Lemma \ref{lem:orbit-lie-group}, the isotropy direction is generated by
$Z_0=Z+Y$ for some $Y\in\mathfrak s$ and
$\theta(Z_0)=1$.  By the smoothness conditions of the codimension-two singular orbit, we have
$|Z_0^\#|_{g_t}=O(t)$ as $t\downarrow0$. Using duality,
\[
  |\theta|_{g_t^{-1}}
  \ge \frac{|\theta(Z_0)|}{|Z_0^\#|_{g_t}}
  =\frac1{|Z_0^\#|_{g_t}}
\]
and hence
\[
  \ell(t)=|\theta|_{g_t^{-1}}^{-1}
  \le |Z_0^\#|_{g_t}=O(t).
\]
Thus $\ell$ is bounded near the bolt. The continuity of $\ell$ gives boundedness on
every compact subinterval of $(0,\infty)$.  Combining this with the estimate
above proves the proposition.
\end{proof}

 The following lemma gives us a type of curvature-volume estimate on the orbit.

\begin{lem}\label{lem:su2-scal-vol}
There is a constant $C$ such that every left-invariant metric $h$ on
$SU(2)$ satisfies
\begin{equation}\label{eq:su2-scal-vol}
  \Scal(h)\le C\,\Vol(SU(2),h)^{-2/3}.
\end{equation}
\end{lem}

\begin{proof}
Fix a bi-invariant metric $h_*$ on $SU(2)$ and write
\[
  V_*=\operatorname{Vol}(SU(2),h_*).
\]
For any left-invariant metric $h$, let $E_1,E_2,E_3$ be an oriented $h_*$-orthonormal basis that diagonalizes $h$. Let $\o_i$ be the dual one-form of $E_i$, then we have
\[
  h=x_1^2\omega_1^2+x_2^2\omega_2^2+x_3^2\omega_3^2,
\]
where $x_i$ are positive constants. Since the bracket of $\mathfrak{su}(2)$ is a fixed positive multiple of the cross product determined by $h_*$ and a suitable
orientation,
\[
  [E_2,E_3]=cE_1,\qquad
  [E_3,E_1]=cE_2,\qquad
  [E_1,E_2]=cE_3
\]
for a positive constant $c$ depending only on $h_*$.

For the $h$-orthonormal basis $e_i=E_i/x_i$, the corresponding
structure constants of brackets of $e_i$ are
\begin{align}\label{eq;def-lambda}
    \lambda_1=\frac{cx_1}{x_2x_3},\qquad
  \lambda_2=\frac{cx_2}{x_3x_1},\qquad
  \lambda_3=\frac{cx_3}{x_1x_2}.
\end{align}
Taking the trace of Milnor's Ricci curvature formula
\cite[Theorem~4.3]{Milnor} gives
\[
  \operatorname{Scal}(h)
  =\lambda_1\lambda_2+\lambda_2\lambda_3+\lambda_3\lambda_1
   -\frac12(\lambda_1^2+\lambda_2^2+\lambda_3^2).
\]
By the elementary inequality
\[
  ab+bc+ca-\frac12(a^2+b^2+c^2)
  \le \frac32(abc)^{2/3},
  \qquad a,b,c>0,
\]
we can get
\begin{align}\label{eq:scalar-bound}
    \operatorname{Scal}(h)\leq \frac{3}{2}(\la_1\la_2\la_3)^{2/3}
\end{align}
Finally, since the coframe $\o_1,\o_2,\o_3$ is orthonormal for $h_*$,
\begin{align}\label{eq:volume-formula}
    d\operatorname{vol}_h
  =x_1x_2x_3\,d\operatorname{vol}_{h_*},
  \qquad
  \operatorname{Vol}(SU(2),h)=V_*x_1x_2x_3.
\end{align}
By \eqref{eq;def-lambda}, \eqref{eq:scalar-bound} and \eqref{eq:volume-formula}, we can get
\[
  \operatorname{Scal}(h)\leq \frac{3}{2}(\la_1\la_2\la_3)^{2/3}
  =C{\operatorname{Vol}(SU(2),h)}^{-2/3},
\]
where $C$ is a constant independent of $h$.
\end{proof}

For the proof of collapsing in this case, let
\begin{align}\label{eq:def-V-volume}
    V(t):=\Vol(P,g_t)
\end{align}
be the volume function of the principal orbits, $h_t=g_t|_{\mathfrak s}$ and $v_h(t)$ be its volume density relative to a fixed Haar normalization on $SU(2)$. By the definition of $\ell,W$ in \eqref{eq:ell-def}, the horizontal vector $\ell W$ satisfies $\theta(\ell W)=1$. Moreover, it is orthogonal
to $\mathfrak s$ and $|\ell W|_{g_t}=\ell$.  Replacing the fixed vector $Z$ by
$\ell W=Z+X$, $X\in\mathfrak s$, has determinant one.  The fixed cocompact
deck group in Lemma \ref{lem:orbit-lie-group} therefore contributes only a
constant, and
\begin{equation}\label{eq:orbit-volume-factor}
  V(t)=C_P\,\ell(t)v_h(t).
\end{equation}
If $s_t$ denotes the intrinsic scalar curvature of the principal orbit,
\eqref{eq:scal-split}, Proposition \ref{prop:ell-bounded}, Lemma \ref{lem:su2-scal-vol}, and \eqref{eq:orbit-volume-factor} give
\begin{equation}\label{eq:orbit-scal-vol}
  s_t\le C\,V(t)^{-2/3}.
\end{equation}

Now we can get an upper volume-growth bound of the geodesic ball centered at a point $o$ of the singular orbit.

\begin{prop}\label{prop:s3-upper}
For a complete non-Ricci-flat soliton in the $\Sph^3\times\mathbb{R}^2$ branch,
\begin{equation}\label{eq:s3-upper}
  \Vol B(o,r)\le C(1+r)^{5/2}.
\end{equation}
\end{prop}

\begin{proof}
By \eqref{eq:xi-H}, \eqref{eq:Vprime} and
\eqref{eq:orbit-scal-vol},
\begin{align}\label{eq:dotHV}
    \dot H+\xi H\le C V^{-2/3},\qquad \dot V=HV.
\end{align}
Let $F=HV^{2/3}$. The estimate $H>0$ from Theorem \ref{thm:cohom-est} gives $F\ge0$. By \eqref{eq:dotHV},
\begin{align}\label{eq:ODE F}
    \dot F
  =V^{2/3}\left(\dot H+\frac23H^2\right)
  \le-\left(\xi-\frac23H\right)F+C.
\end{align}
After a normalization, we assume $C_0=1$ in \eqref{eq:soliton-identities}. Then Theorem \ref{thm:cohom-est} tells us that $\xi\to 1$ and $H\to0$. Hence, $\xi-2H/3$ is bounded below by a positive constant for large $t$.  Integrating the inequality \eqref{eq:ODE F} gives $F\le C$(note that the constant $C$ may change from line to line). Hence, by \eqref{eq:dotHV},
\[
  \frac{\md}{\md t}V^{2/3}=\frac23HV^{2/3}\le C
\]
and
\[
  V(t)\le C(1+t)^{3/2}.
\]
Now we choose a point $o$ on the singular orbit. Since $t(x)=d(x,Q)\le d(x,o)$,
\[
  B(o,r)\subset\{0\le t\le r\}.
\]
Then the proposition follows from integration of the orbit volumes:
\[
  \Vol B(o,r)\le\int_0^rV(t)\,\md t\le C(1+r)^{5/2}.
\]
\end{proof}

On the other hand, $\kappa$-noncollapsing condition and linear curvature decay in Proposition \ref{prop:full-linear} give us a lower bound of the geodesic ball $B(o,r)$, due to the result of Bamler-Chan-Ma-Zhang.

\begin{prop}\label{prop:packing-lower}
Let $(M^5,g,f)$ be a $\kappa$-noncollapsed non-Ricci-flat cohomogeneity one steady soliton. Let $o\in M$ be a fixed basepoint. Then there exists a positive constant $c$ such that
\begin{equation}\label{eq:packing-lower}
  \Vol B(o,r)\ge cr^3
\end{equation}
for all sufficiently large $r$.
\end{prop}

\begin{proof}
This proposition directly follows from Proposition \ref{prop:full-linear}, Proposition \ref{prop:kappa-entropy} and \cite[Theorem 1.1]{BamlerChanMaZhang}
\end{proof}

Combining Proposition \ref{prop:s3-upper} and Proposition \ref{prop:packing-lower} yields the following theorem.

\begin{theo}\label{thm:s3-exclusion}
There is no complete non-Ricci-flat $\kappa$-noncollapsed
five-dimensional cohomogeneity one steady gradient Ricci soliton in the
$\Sph^3\times\mathbb{R}^2$ branch.
\end{theo}

At this stage, we can complete the proof of Theorem \ref{thm:five}.

\begin{proof}[Proof of Theorem \ref{thm:five}]
By Lemma \ref{lem:average}, we can assume that the potential function is also invariant under this cohomogeneity one action. Proposition \ref{prop:kappa-entropy} and Lemma \ref{lemma: lift action} allow us to lift the action to a finite universal cover. Proposition \ref{prop:five-trichotomy} gives exactly three branches.
The fixed-point branch is the five-dimensional Bryant soliton by Proposition \ref{prop: 5d Bryant}.  The $\Sph^2$-bolt branch belongs to $\IDW$ by Proposition \ref{prop:S2-bolt-classification}.  The remaining branch is excluded by Theorem \ref{thm:s3-exclusion}.  Finally, Proposition \ref{prop:idw-noncollapse} and the standard asymptotics of the Bryant soliton show that the listed models are
indeed noncollapsed.
\end{proof}

\appendix

\section{Calculations of ODEs about \texorpdfstring{$m,\mathcal B_n$}{1}}

\renewcommand{\theequation}{A.\arabic{equation}}
\setcounter{equation}{0}

In this appendix, we carry out some computations of differential equations of $m,\mathcal{B}_n$ used in Section~3 in the proof of uniqueness of Appleton's noncollapsed steady solitons.

Recall that
\begin{align*}
    m=\frac{b^\pr}{Q}=\frac{bb^\pr}{a}.
\end{align*}

The following lemma gives us the first and second-order differential equations satisfied by $m$.

\begin{lem}\label{lem:comparison-m-evolution}
The quantity $m$ satisfies
\begin{equation}\label{eq:m-first-appendix}
 m'=\frac{2(n+1-Q^2)}{bQ}
 -\left(2n\frac{b'}b-f'+2\frac{Q'}Q\right)m
\end{equation}
and
\begin{align}\label{eq:m-second-appendix}
 m''+\left((2n-1)\frac{b'}b-f'+3\frac{Q'}Q\right)m'
 =\frac{2(m^2-1)}{b^2}\left(nmQ^2+2bQ'\right).
\end{align}
Moreover, as $t\to 0$,
\begin{align}\label{eq:m-bolt-appendix}
 m(t)=\frac{n+1}{\alpha}
 +\frac{\alpha}{2}\left(\frac{(n+1)^2}{\alpha^2}-1\right)t^2
 +O(t^4).
\end{align}
\end{lem}

\begin{proof}
Substituting $a=Qb$ into the equation for $b$ in \eqref{eq:prelim-ode} gives
\[
 b''=\frac{2(n+1-Q^2)}b-2n\frac{(b')^2}b
 +\left(f'-\frac{Q'}Q\right)b'.
\]
Then \eqref{eq:m-first-appendix} follows from differentiating $m=b'/Q$. Differentiating once more, while using the soliton equations and the quotient equation \eqref{eq:prelim-qode} to eliminate $f''$, $b''$, and $Q''$, gives \eqref{eq:m-second-appendix}.

To determine the expansion at the zero section, note that $m$ is even, hence $m'(0)=0$. Using $Q'/Q=t^{-1}+O(t)$ and taking the limit in \eqref{eq:m-second-appendix} as $t \downarrow 0$ gives
\[
 4m''(0)=4\alpha\bigl(m(0)^2-1\bigr).
\]
Together with the evenness of $m$, this proves \eqref{eq:m-bolt-appendix}.
\end{proof}

The next quantity $\mathcal{B}_n$ is defined by
\begin{align*}
    \mathcal B_n
 :=-\frac{2(n+1)(1-Q^2)}{Q}(b^2Q')'
   +\bigl(nm^2-3n-2\bigr)b^2(Q')^2,
\end{align*}
which is equivalent to
\begin{equation}\label{eq:Bn-diff-appendix}
    \begin{aligned}
        \mathcal B_n
 ={}&2(n+1)(1-Q^2)
 \left[\frac{bQ'}Q\bigl((2n-1)b'-bf'\bigr)
       +2(n+1)(1-Q^2)\right]\\
 &+\bigl(nm^2-3n-2\bigr)b^2(Q')^2.
    \end{aligned}
\end{equation}
by using the equation \eqref{eq:prelim-qode}. We give a detailed calculation of the differential equation satisfied by $\mathcal{B}_n$ at a zero of it.

\begin{lem}\label{lem:Bn-zero-identities}
Let
\[
 \vartheta(x)=\frac{bb'Q'}{Q(1-Q^2)},\qquad
 \zeta(x)=-\frac{b^2f'Q'}{Q(1-Q^2)},
\]
where $x=Q^2$. If $\mathcal B_n(x_0)=0$ for some $x_0\in(0,1)$, then at $x=x_0$,
\begin{equation}\label{eq:Bn-zero-constraint}
 \frac{(3n+2)b^2(Q')^2}{(1-x)^2}
 =n\vartheta^2+2(n+1)(2n-1)\vartheta
   +2(n+1)\zeta+4(n+1)^2.
\end{equation}
Moreover, throughout $(0,1)$,
\begin{equation}\label{eq:m-x-appendix}
 (bQ')m_x
 =\frac{n+(1-x)(1-\vartheta)}{x}
 -\frac{(1-x)^2\vartheta(2n\vartheta+\zeta)}{2b^2(Q')^2},
\end{equation}
and at every zero $x_0$ of $\mathcal B_n$,
\begin{equation}\label{eq:Bn-x-appendix}
\begin{aligned}
 (\mathcal B_n)_x
 ={}&2(1-x)(n\vartheta-n-1)
       \frac{n+(1-x)(1-\vartheta)}{x}\\
 &-\frac{(1-x)^3}{(3n+2)b^2(Q')^2}
   \Bigl[(n+1)\mathscr P_n(\vartheta,\zeta)
       +n(3n+2)\vartheta^2(2n\vartheta+\zeta)\Bigr],
\end{aligned}
\end{equation}
where $\mathscr P_n$ is the polynomial defined in~\eqref{eq:comparison-polynomial}.
\end{lem}

\begin{proof}
For simplicity, let
\[
 w=b^2(Q')^2,\qquad
 H=n\vartheta^2+2(n+1)(2n-1)\vartheta
      +2(n+1)\zeta+4(n+1)^2.
\]
The definitions give
\[
 b'=\frac{\vartheta\sqrt{x}(1-x)}{\sqrt w},\qquad
 bf'=-\frac{\zeta\sqrt{x}(1-x)}{\sqrt w},\qquad
 m=\frac{\vartheta(1-x)}{\sqrt w}.
\]
Substituting these identities into \eqref{eq:Bn-diff-appendix} yields
\begin{align}\label{eq:another-Bn-appendix}
    \mathcal B_n=(1-x)^2H-(3n+2)w
\end{align}
on the entire interval $(0,1)$.
Its value at a zero of $\mathcal B_n$ gives
\eqref{eq:Bn-zero-constraint}.

We first rewrite the equation for $m$. By \eqref{eq:comparison-m-first-order} and
$\partial_x=(2QQ')^{-1}\partial_t$,
\[
\begin{aligned}
 (bQ')m_x
 &=\frac{n+1-x}{x}
   -\frac{b}{2Q}\left(2n\frac{b'}b-f'
                         +2\frac{Q'}Q\right)m\\
 &=\frac{n+1-x-(1-x)\vartheta}{x}
   -\frac{(1-x)^2\vartheta(2n\vartheta+\zeta)}{2w}.
\end{aligned}
\]
This proves \eqref{eq:m-x-appendix}.

We next compute the derivatives needed to differentiate $\mathcal B_n$.
The quotient equation \eqref{eq:prelim-qode} gives
\[
\begin{aligned}
 w_x
 &=\frac{bb'Q'+b^2Q''}{Q}\\
 &=\frac{bQ'}Q(bf'-2nb')-2(n+1)(1-x)\\
 &=-(1-x)(2n\vartheta+\zeta+2n+2).
\end{aligned}
\]
Differentiating $\vartheta=m\sqrt w/(1-x)$ and using \eqref{eq:m-x-appendix}, we obtain
\[
 \vartheta_x
 =\frac{n+1-x+(2x-1)\vartheta}{x(1-x)}
  -\frac{(1-x)\vartheta(2n\vartheta+\zeta+n+1)}{w}.
\]
For the derivative of $\zeta$, the soliton equations \eqref{eq:prelim-ode} imply
\[
\begin{aligned}
 b^2f''
 ={}&2n(2n+2-x)-(1-x)(4n\vartheta+\zeta)\\
 &-\frac{(2n+1)x(1-x)^2\vartheta(2n\vartheta+\zeta)}{w}.
\end{aligned}
\]
Differentiating $\zeta=-bf'\sqrt w/(\sqrt x(1-x))$ now yields
\[
\begin{aligned}
 \zeta_x
 ={}&\frac{x\zeta+2n(1-x)\vartheta-n(2n+2-x)}{x(1-x)}\\
 &+\frac{1-x}{2w}
 \Bigl[2n(2n+1)\vartheta^2+2\vartheta\zeta
                  -\zeta^2-2(n+1)\zeta\Bigr].
\end{aligned}
\]
By \eqref{eq:another-Bn-appendix},
\[
\begin{aligned}
 (\mathcal B_n)_x
 ={}&-2(1-x)H
 +2(1-x)^2\bigl[n\vartheta+(n+1)(2n-1)\bigr]\vartheta_x\\
 &+2(n+1)(1-x)^2\zeta_x-(3n+2)w_x.
\end{aligned}
\]
Let $x_0$ be a zero of $\mathcal{B}_n$. If we evaluate at $x_0$, then $(3n+2)w=(1-x)^2H$. Now using the derivative formulas of $w_x,\,\vartheta_x,$ and $\zeta_x$, we can get at $x_0$,
\[
\begin{aligned}
 &\frac{(3n+2)w}{(1-x)^3}
 \left[(\mathcal B_n)_x
 -2(1-x)(n\vartheta-n-1)
       \frac{n+(1-x)(1-\vartheta)}{x}\right]\\
 &\quad=\bigl[2(n^2+n+1)\vartheta+n\zeta-4(n+1)^2\bigr]H\\
 &\qquad\quad
 -2(3n+2)\bigl[n\vartheta+(n+1)(2n-1)\bigr]
       \vartheta(2n\vartheta+\zeta+n+1)\\
 &\qquad\quad
 +(n+1)(3n+2)
       \bigl[2n(2n+1)\vartheta^2+2\vartheta\zeta
                       -\zeta^2-2(n+1)\zeta\bigr]\\
 &\quad=-(n+1)\mathscr P_n(\vartheta,\zeta)
        -n(3n+2)\vartheta^2(2n\vartheta+\zeta),
\end{aligned}
\]
where the last equality follows by expanding $H$ and grouping the powers of
$\zeta$ and $\vartheta$ and $\mathscr P_n$ is the polynomial defined in~\eqref{eq:comparison-polynomial}. Rearranging the above identity proves \eqref{eq:Bn-x-appendix}.
\end{proof}

\section{Simply connected compact homogeneous four-manifolds}
\label{app:homogeneous-four}

\renewcommand{\theequation}{B.\arabic{equation}}
\setcounter{equation}{0}

For completeness we prove the low-dimensional homogeneous-space
classification used in Proposition \ref{prop:five-trichotomy}.

\begin{theo}\label{thm:homogeneous-four}
Let $M^4=G/H$ be a compact and simply connected homogeneous space under a transitive compact and connected Lie group $G$ action. Then $M$ is diffeomorphic to one of the following three manifolds
\[
  \Sph^4,\qquad \mathbb{CP} ^2,\qquad
  \Sph^2\times\Sph^2.
\]
After removing ineffective kernels and passing to finite covers of the
acting groups, the corresponding homogeneous presentations are
\[
  \frac{Sp(2)}{Sp(1)\times Sp(1)},\qquad
  \frac{SU(3)}{S(U(2)\times U(1))},\qquad
  \frac{SU(2)\times SU(2)}{\Sph^1\times\Sph^1}.
\]
\end{theo}

\begin{proof}
After quotient by the ineffective kernel, we assume that the action of $G$ is effective.  The fibration
$H\to G\to G/H$ and simple connectivity of $M$ show that $H$ is connected.
Since a closed simply connected four-manifold is orientable and has
\[
  \chi(M)=2+b_2(M)>0,
\]
the Hopf--Samelson theorem \cite{HopfSamelson1940} gives
\[
\operatorname{rank}G=\operatorname{rank}H.
\]
Then we can choose a common maximal torus $T\subset H\subset G$. Since the connected center
$Z(G)^0$ lies in $T$, we have $Z(G)^0\subset H$. As a subset of the center, $Z(G)^0$ acts trivially on $G/H$. By the effectivity of the $G$-action, $Z(G)^0=\{e\}$. This implies that $G$ is semisimple. After a finite covering we may assume it is a product of compact simply
connected simple groups.

Let $\Phi$ be the root system of $G$ relative to $T$.  Since
$\mathfrak h_\mathbb{C}$ is $T$-invariant, there is a closed root subsystem
$\Psi\subset\Phi$ such that
\[
  \mathfrak h_\mathbb{C}
  =\mathfrak t_\mathbb{C}\oplus\bigoplus_{\alpha\in\Psi}\mathfrak g_\alpha.
\]
After choosing compatible positive roots, we have
\[
  4=\dim G/H=2\bigl(|\Phi^+|-|\Psi^+|\bigr).
\]
This implies that
\begin{equation}\label{eq:root-difference}
  |\Phi^+|-|\Psi^+|=2.
\end{equation}
Now we decompose $\Phi=\bigsqcup\Phi_i$ according to the simple factors and put
$d_i=|\Phi_i^+|-|\Psi_i^+|$.  If $d_i=0$, then the whole simple factor would lie in $H$ and hence in the ineffective kernel. By the effectivity of the action, we must have $d_i>0$.
Thus \eqref{eq:root-difference} tells us that there is one simple factor with $d_1=2$ or two simple factors with $d_1=d_2=1$.

For an irreducible root system, using the Borel--de Siebenthal classification of connected maximal-rank subgroups \cite{BorelDeSiebenthal1949} and the bound $|\Phi^+|-|\Psi^+|\leq2$, we obtain the following possibilities:
\begin{align}\label{eq:three-poss-appendix}
    (A_1,\varnothing),\qquad
  (A_2,A_1),\qquad
  (B_2,A_1\sqcup A_1).
\end{align}
If there is one simple factor, then $d_1 = 2$. Hence, $G$ is simple and $\dim(G/H) = 4$. By the possibilities in \eqref{eq:three-poss-appendix}, there are the following two cases.
\begin{align*}
    (\mathfrak g,\mathfrak h)
  =(\mathfrak{su}(3),\mathfrak{su}(2)\oplus\mathbb{R}),\qquad G/H\cong \mathbb {CP}^2,\\
  (\mathfrak g,\mathfrak h)
  =(\mathfrak{sp}(2),\mathfrak{sp}(1)\oplus\mathfrak{sp}(1)),\qquad G/H\cong \Sph^4.
\end{align*}
If there are two simple factors, then $d_1=d_2=1$. By \eqref{eq:three-poss-appendix}, each must have
root system $A_1$ and the isotropy has no roots.  Hence
\[
  (\mathfrak g,\mathfrak h)
  =(\mathfrak{su}(2)\oplus\mathfrak{su}(2),\mathbb{R}\oplus\mathbb{R}),\qquad G/H\cong \Sph^2\times \Sph^2.
\]
These cases exhaust \eqref{eq:root-difference} and the theorem is proved.
\end{proof}

\begin{rem}
The classification above is classical~(see Gorbatsevich~\cite{Gorbatsevich1977,Gorbatsevich2012}). The proof was included to make
clear why the Euler characteristic $4$ in the $\Sph^2$-bolt branch forces $\Sph^2\times \Sph^2$ rather than the twisted $\Sph^2$-bundle over $\Sph^2$.
\end{rem}

\end{document}